\documentclass[11pt]{article}
\usepackage{graphicx}
\usepackage{natbib,amssymb}
\usepackage{amsmath,multirow}
\usepackage{amsthm}
\usepackage{geometry}
\usepackage{url}

\newtheorem{prop}{Proposition}[section]
\newtheorem{theo}[prop]{Theorem}
\newtheorem{lem}[prop]{Lemma}
\newtheorem{cond}[prop]{Condition}
\newtheorem{cor}[prop]{Corollary}

\newtheorem{rem}[prop]{Remark}
\newtheorem{alg}[prop]{Algorithm}
\input{tcilatex}

\begin{document}

\title{Valid confidence intervals for post-model-selection predictors}
\author{Fran\c{c}ois Bachoc*, Hannes Leeb**, and Benedikt M. P\"{o}tscher** 
\\
*Department of Mathematics, University Paul Sabatier\\
**Department of Statistics, University of Vienna}
\date{First version: December 2014\\
This version: December 2016}
\maketitle

\begin{abstract}
We consider inference post-model-selection in linear regression. In this
setting, \cite{berk13valid} recently introduced a class of confidence sets,
the so-called PoSI intervals, that cover a certain non-standard quantity of
interest with a user-specified minimal coverage probability, irrespective of
the model selection procedure that is being used. In this paper, we
generalize the PoSI intervals to confidence intervals for
post-model-selection predictors.

\medskip

{\footnotesize \textsc{AMS Mathematics Subject Classification 2010:} 62F25,
62J05.\vspace{0.2cm}}


{\footnotesize \ }
\end{abstract}

\section{Introduction and overview\label{section:intro}}

In statistical practice, the model used for analysis is very often chosen
after the data have been observed, either by ad-hoc methods or by more
sophisticated model selection procedures. Inference following such a model
selection step (inference post-model-selection) has proven to be a
challenging problem. `Naive' procedures, which ignore the presence of model
selection, are typically invalid (e.g., in the sense that the actual
coverage probability of `naive' confidence sets for the true parameter can
be dramatically smaller than the nominal one), and the construction of valid
procedures is often non-trivial; see \cite{leeb05model,leeb06estimate,LP2012}%
, \cite{kabaila06large}, \cite{Poe09a} and references therein for an
introduction to the issues involved here. In these references, inference is
focused on the true parameter of the data-generating model (or on components
thereof). Shifting the focus away from the true parameter as the target of
inference, \cite{berk13valid} recently introduced a class of confidence
sets, the so-called PoSI intervals, that guarantee a user-specified minimal
coverage probability after model selection in linear regression,
irrespective of the model selector that is being used; see also \cite%
{berk13validold} and \cite{leeb13various}. In this paper, we generalize the
PoSI intervals to intervals for post-model-selection predictors.

Prediction following model selection is obviously also of great importance.
In the case where the selected model is misspecified, parameter estimates
are typically biased or at least difficult to interpret; cf. Remark \ref%
{rem:interpretation}. But even a misspecified model may perform well for
prediction. In particular, \cite{Gree04a} derive, under appropriate sparsity
assumptions, feasible predictors that asymptotically perform as well as the
(infeasible) best candidate predictor even if the available number of
explanatory variables by far exceeds the sample size. These feasible
predictors are also covered by the results in the present paper, among
others. Like \cite{Gree04a}, our analysis does not rely on the assumption
that the true data generating model is among the candidates for model
selection. We develop confidence intervals for such predictors, that are
easy to interpret and that are optimal in an appropriate sense; cf. Remarks %
\ref{rem:interpretation}(ii), \ref{rem:predictor}, and \ref{rem:predictor_2}%
, as well as \cite{Gree04a}. A further rationale for extending the
PoSI-approach of \cite{berk13valid} to problems related to prediction is
that this framework seems to provide a more natural habitat for considering
non-standard targets; see the discussion in Remark 2.1 of \cite%
{leeb13various} as well as in Remarks \ref{rem:interpretation} and \ref%
{rem:justification} given further below.

The crucial feature of the approach of \cite{berk13valid} is that the
coverage target, i.e., the quantity for which a confidence set is desired,
is \emph{not} the standard target, i.e., the parameter in an overall model
(or components thereof), but a \emph{non-standard} quantity of interest that
depends on the selected model and thus on the data. This non-standard
quantity of interest is denoted by $\beta _{\hat{M}}^{(n)}$ throughout the
paper (cf. Section \ref{section:design:dependent} for details). Here $\hat{M}
$ stands for the (data-dependent) model chosen by the model selector and $n$
stands for sample size. The non-standard target $\beta _{\hat{M}}^{(n)}$
provides a certain vector of regression coefficients for those explanatory
variables that are `active' in the model $\hat{M}$ (more precisely, $\beta _{%
\hat{M}}^{(n)}$ represents the coefficients of the projection of the
expected value vector of the dependent variable on the space spanned by the
regressors included in $\hat{M}$); for a precise definition see eqs. (\ref%
{eq:def:betaMn}) and (\ref{eq:def:Beta:hatM:n}) in Section \ref%
{section:design:dependent}.

For a new set of explanatory variables $x_{0}$, we first extend the
PoSI-approach to obtain confidence intervals for the predictor $%
x_{0}^{\prime }[\hat{M}]\beta _{\hat{M}}^{(n)}$. Here, $x_{0}[\hat{M}]$
denotes the set of explanatory variables from $x_{0}$ that correspond to the
`active' regressors in the model $\hat{M}$. We call $x_{0}^{\prime }[\hat{M}%
]\beta _{\hat{M}}^{(n)}$ the design-dependent (non-standard) coverage
target, because different design matrices in the training data typically
result in different values of $x_{0}^{\prime }[\hat{M}]\beta _{\hat{M}%
}^{(n)} $ even if both training data sets lead to selection of the same
model $\hat{M}$. We construct PoSI confidence intervals for $x_{0}^{\prime }[%
\hat{M}]\beta _{\hat{M}}^{(n)}$ that guarantee a user-specified minimal
coverage probability, irrespective of the model selector that is being used.
The design-dependent coverage target minimizes a certain \emph{`in-sample'}
prediction error; cf. Remark \ref{rem:predictor}. However, when the goal is
to predict a new response corresponding to a new vector $x_{0}$ of
explanatory variables, this `in-sample' optimality property may have little
relevance and thus the focus on covering the design-dependent target $%
x_{0}^{\prime }[\hat{M}]\beta _{\hat{M}}^{(n)}$ may be debatable.

In view of this, we next consider an alternative coverage target that
depends on the selected model but not on the training data otherwise, and
that we denote by $x_{0}^{\prime }[\hat{M}]\beta _{\hat{M}}^{(\star )}$. We
call $x_{0}^{\prime }[\hat{M}]\beta _{\hat{M}}^{(\star )}$ the
design-independent (non-standard) coverage target. The design-independent
coverage target minimizes a certain `\emph{out-of-sample}' prediction error,
namely the mean-squared prediction error, over all (infeasible) predictors
of a future response $y_{0}$ that are of the form $x_{0}^{\prime }[\hat{M}%
]\gamma (\hat{M})$, when $x_{0}$ and the row-vectors of $X$ are sampled from
the same distribution; cf. Remark \ref{rem:predictor_2}. In particular, this
target does not suffer from the issues that plague the design-dependent
coverage target, as discussed at the end of the preceding paragraph. Certain
optimality properties of a feasible counterpart of $x_{0}^{\prime }[\hat{M}%
]\beta _{\hat{M}}^{(\star )}$ are derived in \cite{Gree04a}, for a
particular model selector $\hat{M}$ and under appropriate sparsity
assumptions; a target closely related to $x_{0}^{\prime }[\hat{M}]\beta _{%
\hat{M}}^{(\star )}$ is also studied in \cite{Lee09a}. For a large class of
model selectors, we show that the PoSI confidence intervals constructed
earlier also cover the design-independent coverage target with minimal
coverage probability not below the user-specified nominal level
asymptotically. In that sense, the PoSI confidence intervals are
approximately valid for the target $x_{0}^{\prime }[\hat{M}]\beta _{\hat{M}%
}^{(\star )}$, irrespective of the model selector $\hat{M}$ in that class.
In simulations we find that our asymptotic result is representative of the
finite-sample situation even for moderate sample sizes.

When extending the PoSI-approach to confidence intervals for both the
design-dependent and the design-independent coverage target, i.e., for both $%
x_{0}^{\prime }[\hat{M}]\beta _{\hat{M}}^{(n)}$ and $x_{0}^{\prime }[\hat{M}%
]\beta _{\hat{M}}^{(\star )}$, we find that the resulting intervals
necessarily depend not only on $x_{0}[\hat{M}]$ but also on those components
of $x_{0}$ that are `in-active' in the model $\hat{M}$. This may appear
surprising at first sight but turns out to be inherent to the PoSI-approach
(because of the need to take the maximum over all models $M$ in (\ref%
{eq:def:posi:xzero})). In any case, this is problematic in situations when,
after having selected a given model, only the `active' components of $x_{0}$
are observed, e.g., in situations where observations are costly and model
selection is carried out also with the goal of reducing cost by not having
to observe irrelevant components of $x_{0}$. To resolve this, we also
develop PoSI confidence intervals that depend on the `active' variables $%
x_{0}[\hat{M}]$ only. These intervals are obtained by maximizing over all
inactive variables and are hence larger than the intervals for the case
where $x_{0}$ is known entirely. In simulations, we find that the excess
width of these intervals is moderate. We also provide analytic results
regarding the excess width of these intervals in an asymptotic setting where
the number of regressors goes to infinity, see Section \ref{sub2.3}.

Inference post-model-selection is currently a very active area of research
and we can only give a selection of work relevant for, or related to, this
paper. Contemporary analyses of confidence sets for (components of) the true
parameter of the underlying model include \cite{And09a}, \cite%
{kabaila06large}, \cite{leeb05model}, \cite{Poe09a}, \cite{Poe10a}, and \cite%
{Sch14a}. These references also point to numerous earlier results. Also, the
work of \cite{Loc14a}, \cite{Was09a}, and \cite{Was14a} should be mentioned
here. For the LASSO, in particular, a de-sparsifying method has recently
been developed by \cite{Bel11a,Bel14a}, \cite{van14a}, and \cite{Zha14a}.
Another strand of literature that, like the PoSI approach, also focuses on $%
\beta _{\hat{M}}^{(n)}$ as the quantity of interest, is developed in \cite%
{Fit15a}, \cite{lee15exact}, \cite{LeeTay14}, \cite{Tia15a}, \cite{Tib15a},
and \cite{Lock15}: In these papers, confidence sets for $\beta _{\hat{M}%
}^{(n)}$ are considered that have a guaranteed coverage probability \emph{%
conditionally} on the event that a particular model has been selected by the
model selection procedure. In contrast to PoSI procedures, the confidence
intervals obtained in these papers are specific to the model selection
procedure used (the LASSO, in particular, being considered in these
references) and generally rely on certain geometric properties of the
specific model selection procedure under consideration. In simulation
experiments, we compare the confidence intervals proposed in these
references with the intervals developed here and observe some interesting
phenomena, see Section \ref{taylor}. As prompted by a referee, we point out
here that in the presence of a large number of regressors PoSI intervals
(including intervals considered in the present paper) typically are
computationally more burdensome than the confidence intervals proposed in 
\cite{lee15exact} for the LASSO with a fixed value for the tuning parameter;
see, however, also the discussion towards the end of Section \ref{taylor}.

The rest of the paper is organized as follows. In Section \ref%
{section:design:dependent}, we introduce the models, the model-selection
procedures, the design-dependent target, and the PoSI confidence intervals
for both the case where all explanatory variables in $x_{0}$ are observed
and the case where only the components of $x_{0}$ corresponding to the
`active' explanatory variables are available; moreover, we analyze
properties of these intervals in an asymptotic framework where the model
dimension increases; cf. Section \ref{sub2.3}. In Section \ref%
{section:distribution:dependent}, we present the design-independent target
and show that the PoSI confidence intervals introduced earlier also cover
the design-independent target, with minimal coverage probability not below
the nominal one asymptotically when sample size increases. The results of a
numerical study are reported in Section \ref{section:simulation:study}.
Conclusions are drawn in Section~\ref{conclusio}. Appendix \ref%
{Section_on_sigma} contains some comments on the assumptions made on the
errors variance. The proofs of the results in Sections \ref%
{section:design:dependent} and \ref{section:distribution:dependent} are
given in Appendices \ref{app A} and \ref{app:B}. Appendix \ref{app_various}
contains some comments on and extensions of the results in Section \ref%
{section:distribution:dependent}. In Appendix \ref%
{section:practical:algorithms} we describe algorithms for computing the PoSI
confidence intervals, that are comparable with those proposed by \cite%
{berk13valid} in terms of computational complexity. Finally, Appendix \ref%
{app D} contains details concerning the numerical calculations used for the
results in Section \ref{section:simulation:study}.

\section{Confidence intervals for the design-dependent non-standard target 
\label{section:design:dependent}}

\subsection{The framework}

Consider the model

\begin{equation}
Y=\mu +U  \label{eq:full:linear:model}
\end{equation}%
where $\mu \in \mathbb{R}^{n}$ is unknown and $U$ follows an $N\left(
0,\sigma ^{2}I_{n}\right) $-distribution; here $\sigma ^{2}$, $0<\sigma
<\infty $, is the unknown error variance and $I_{n}$ is the identity matrix
of size $n\geq 1$. An important instance of this model arises when $\mu $ is
known to reside in a lower dimensional linear subspace of $\mathbb{R}^{n}$,
but we do not make such an assumption at this point. Apart from the data $Y$%
, we are given a (real) $n\times p$ matrix $X$, not necessarily of full
column rank, the columns of which represent potential regressors. This setup
allows for $p>n$ as well as for $1\leq p\leq n$. The rank of $X$ will be
denoted by $d$. The design matrix $X$ is treated as fixed throughout Section %
\ref{section:design:dependent}.

We consider fitting (potentially misspecified) linear models with design
matrices that are obtained by deleting columns from $X$. Such a model will
be represented by $M$, a subset of $\left\{ 1,...,p\right\} $, where the
elements of $M$ index the columns of $X$ that are retained. We use the
following notation: For $M\subseteq \left\{ 1,...,p\right\} $, we write $%
M^{c}$ for the complement of $M$ in $\{1,...,p\}$. It proves useful to allow 
$M$ to be the empty set. We write $|M|$ for the cardinality of $M$. With $%
m=|M|$, let us write $M=\left\{ j_{1},...,j_{m}\right\} $ in case $m\geq 1$.
For $M\neq \varnothing $ and for an $l\times p$ matrix $T$, $l\geq 1$, let $T%
{[M]}$ be the matrix of dimension $l\times m$ obtained from $T$ by retaining
only the columns of $T$ with indices $j\in M$ and deleting all others; if $%
M=\varnothing $ we set $T{[M]}=0\in \mathbb{R}^{l}$. In abuse of notation we
shall, for a $p\times 1$ vector $v$, write $v{[M]}$ for $\left( v^{\prime }{%
[M]}\right) ^{\prime }$, i.e., $v{[M]}=(v_{j_{1}},...,v_{j_{m}})^{\prime }$
for $m\geq 1$ and $v{[M]}=0\in \mathbb{R}$ in case $M=\varnothing $. For a
given model $M$, we denote the corresponding least squares estimator by $%
\hat{\beta}_{M}$, i.e., 
\begin{equation}
\hat{\beta}_{M}=\left( X{[M]}^{\prime }{X[M]}\right) ^{-1}X{[M]}^{\prime }{Y}%
,  \label{eq:restr_OLS}
\end{equation}%
where the inverse is to be interpreted as the Moore-Penrose inverse in case $%
{X[M]}$ does not have full column rank. For any given model $M$ the
corresponding least squares estimator $\hat{\beta}_{M}$ is obviously an
unbiased estimator of 
\begin{equation}
\beta _{M}^{(n)}=\left( X{[M]}^{\prime }{X[M]}\right) ^{-1}X{[M]}^{\prime }{%
\mu .}  \label{eq:def:betaMn}
\end{equation}%
Note that $\hat{\beta}_{M}$ as well as $\beta _{M}^{(n)}$ reduce to $0$ in
case $M=\varnothing $.

As in \cite{berk13valid} we further assume that, as an estimator for $\sigma
^{2}$, we have available an (observable) random variable $\hat{\sigma}^{2}$
that is independent of $P_{X}Y$ and that is distributed as $\sigma ^{2}/r$
times a chi-square distributed random variable with $r$ degrees of freedom ($%
1\leq r<\infty $), with $P_{X}$ denoting orthogonal projection on the column
space of $X$. This assumption is always satisfied in the important special
case where one assumes that $d<n$ and $\mu \in \limfunc{span}(X)$ hold, upon
choosing for $\hat{\sigma}^{2}$ the standard residual variance estimator
obtained from regressing $Y$ on $X$ and upon setting $r=n-d$. However,
otherwise it is not an innocuous assumption at all and this is further
discussed in Appendix \ref{Section_on_sigma}. Observe that our assumption
allows for estimators $\hat{\sigma}^{2}$ that not only depend on $Y$ and $X$%
, but possibly also on other observable random variables (e.g., additional
data). The joint distribution of $Y$ and $\hat{\sigma}^{2}$ depends on $\mu $
and $\sigma $ as well as on sample size $n$ and will be denoted by $P_{n,\mu
,\sigma }$ (see also Appendix \ref{notation}).

We are furthermore given a (non-empty) collection $\mathcal{M}$ of
admissible models $M\subseteq \left\{ 1,...,p\right\} $, the `universe' of
models considered by the researcher. Without loss of generality we will
assume that any column of $X$ appears as a regressor in at least one of the
models $M$ in $\mathcal{M}$, i.e., that $\bigcup \left\{ M:M\in \mathcal{M}%
\right\} =\left\{ 1,...,p\right\} $ holds (otherwise we can just redefine $X$
by discarding all columns that do not appear in any of the models in $%
\mathcal{M}$); of course, we have excluded here the trivial and
uninteresting case $\mathcal{M}=\left\{ \varnothing \right\} $. For such a
collection $\mathcal{M}$ it is easy to see that the assumed independence of $%
\hat{\sigma}^{2}$ and $P_{X}Y$ is in fact equivalent to independence of $%
\hat{\sigma}^{2}$ from the collection $\{\hat{\beta}_{M}:M\in \mathcal{M}\}$
of least squares estimators. While not really affecting the results, it
proves useful to assume, throughout the following, that the empty model
belongs to $\mathcal{M}$. We shall furthermore always assume that any
non-empty $M\in \mathcal{M}$ is of full-rank in the sense that $\limfunc{rank%
}{X[M]}=|M|$. We point out here that our assumptions on $\mathcal{M}$ imply
that $X$ can not have a zero column, and hence $d\geq 1$ must hold. An
important instance of a collection $\mathcal{M}$ satisfying our assumptions
is the collection of all full-rank submodels of $\left\{ 1,...,p\right\} $
(enlarged by the empty model) provided that no column of $X$ is zero; of
course, there are many other examples, see, e.g., the list in Section 4.5 of 
\cite{berk13valid}.

A model selection procedure $\hat{M}$ is now a (measurable) rule that
associates with every $(X,Y,\hat{\sigma}^{2})$ a (possibly empty) model $%
\hat{M}(X,Y,\hat{\sigma}^{2})\in \mathcal{M}$. In the following we shall, in
abuse of notation, often write $\hat{M}$ for $\hat{M}(X,Y,\hat{\sigma}^{2})$%
. Allowing explicitly dependence of $\hat{M}$ on $\hat{\sigma}^{2}$ is only
relevant in case $\hat{\sigma}^{2}$ depends on extraneous data beyond $(X,Y)$
and the model selection procedure actually makes use of $\hat{\sigma}^{2}$.
[We note here that in principle we could have allowed $\hat{M}$ to depend on
further extraneous data, in which case $P_{n,\mu ,\sigma }$ would have to be
redefined as the joint distribution of $Y$, $\hat{\sigma}^{2}$, and this
further extraneous data.] The post-model-selection estimator $\hat{\beta}_{%
\hat{M}}$ corresponding to the model selection procedure is now given by (%
\ref{eq:restr_OLS}) with $M$ replaced by $\hat{M}$.

The non-standard quantity of interest studied in \cite{berk13valid} is the
random vector (with random dimension) $\beta _{\hat{M}}^{(n)}$ obtained by
replacing $M$ by $\hat{M}$ in (\ref{eq:def:betaMn}). The situation we shall
consider in the present paper is related to \cite{berk13valid}, but is
different in several aspects: Consider a fixed (real) $p\times 1$ vector $%
x_{0}$ and suppose we want to predict $y_{0}$ which is distributed as $%
N\left( \nu ,\sigma ^{2}\right) $, independently of $Y$. If one is forced to
use a fixed model $M$ for prediction, i.e., to use predictors of the form $%
x_{0}^{\prime }{[M]\gamma }$, the predictor that would then typically be
used is $x_{0}^{\prime }{[M]}\hat{\beta}_{M}$, which can be viewed as an
estimator of the infeasible predictor $x_{0}^{\prime }{[M]}\beta _{M}^{(n)}$%
. Of course, for this predictor to be reasonable there must be some relation
between the training data $(X,Y)$ and $(x_{0},y_{0})$. This is further
discussed in Remark \ref{rem:predictor}. In the presence of model selection
the predictor $x_{0}^{\prime }{[M]}\hat{\beta}_{M}$ will then typically be
replaced by the post-model-selection predictor $x_{0}^{\prime }{[\hat{M}]}%
\hat{\beta}_{\hat{M}}$ which can in turn be seen as a feasible counterpart
to the infeasible predictor%
\begin{equation}
x_{0}^{\prime }{[\hat{M}]}\beta _{\hat{M}}^{(n)}.  \label{eq:def:Beta:hatM:n}
\end{equation}%
The quantity in (\ref{eq:def:Beta:hatM:n}) will be our target for inference
throughout Section \ref{section:design:dependent} and will be called the 
\emph{design-dependent (non-standard) target} (to emphasize that it depends
on the design matrix $X$ apart from its dependence on ${\hat{M}}$, cf. (\ref%
{eq:def:betaMn})). A discussion of the merits of this target and its
interpretation is postponed to Remarks \ref{rem:interpretation} and \ref%
{rem:predictor} given below.

Let now $1-\alpha \in (0,1)$ be a nominal confidence level. Throughout
Section \ref{section:design:dependent} we are interested in confidence
intervals for the design-dependent target $x_{0}^{\prime }{[\hat{M}]}\beta _{%
\hat{M}}^{(n)}$ that are of the form%
\begin{equation}
CI(x_{0})=x_{0}^{\prime }{[\hat{M}]}\hat{\beta}_{\hat{M}}\pm K(x_{0},\hat{M}%
)||s_{{\hat{M}}}||\hat{\sigma},  \label{eq:general:form:CI}
\end{equation}%
where $\left\Vert \cdot \right\Vert $ denotes the Euclidean norm ($\hat{%
\sigma}$ of course representing the nonnegative square root of $\hat{\sigma}%
^{2}$), where%
\begin{equation}
s_{M}^{\prime }=x_{0}^{\prime }{[M]}\left( X{[M]}^{\prime }{X[M]}\right)
^{-1}X{[M]}^{\prime }{,}  \label{eq:def:tM}
\end{equation}%
where $s_{M}=0\in \mathbb{R}^{n}$ for $M=\varnothing $ by our conventions,
and where $K(x_{0},M)=K(x_{0},M,r)=K(x_{0},M,r,X,\alpha ,\mathcal{M})$
denotes a non-negative constant which may depend on $x_{0}$, $M$, $r$, $X$, $%
\alpha $, and $\mathcal{M}$, but does not depend on the observations on $Y$
and $\hat{\sigma}^{2}$. Here we have used the notation $a\pm b$ for the
interval $\left[ a-b,a+b\right] $ ($a\in \mathbb{R}$, $b\geq 0$). The
motivation for the form of the confidence interval stems from the
observation that for \emph{fixed} $M$ the interval $x_{0}^{\prime }{[M]}\hat{%
\beta}_{M}\pm q_{r,1-\alpha /2}||s_{{M}}||\hat{\sigma}$ is a valid $1-\alpha 
$ confidence interval for $x_{0}^{\prime }{[M]}\beta _{M}^{(n)}$, where $%
q_{r,1-\alpha /2}$ is the $(1-\alpha /2)$-quantile of Student's
t-distribution with $r$ degrees of freedom. Furthermore note that on the
event ${\hat{M}=}\varnothing $ the target is equal to zero and the
confidence interval reduces to $\{0\}$, thus always containing the target on
this event. Finally note that $CI(x_{0})$ constitutes a confidence interval
for the predictor $x_{0}^{\prime }{[\hat{M}]}\beta _{\hat{M}}^{(n)}$, and
should not be mistaken for a prediction interval for a new response $y_{0}$.

We aim at finding quantities $K(x_{0},M)$ such that the confidence intervals 
$CI(x_{0})$ satisfy%
\begin{equation}
\inf_{\mu \in \mathbb{R}^{n},\sigma >0}P_{n,\mu ,\sigma }\left(
x_{0}^{\prime }{[\hat{M}]}\beta _{\hat{M}}^{(n)}\in CI(x_{0})\right) \geq
1-\alpha .  \label{eq:PRED:PoSI:coverage:aim}
\end{equation}%
Note that if one replaces $K(x_{0},\hat{M})$ in (\ref{eq:general:form:CI})
by $K_{naive}=q_{r,1-\alpha /2}$, then the confidence interval (\ref%
{eq:general:form:CI}) reduces to the so-called `naive' confidence interval
which is constructed as if $\hat{M}$ were fixed a priori (thus ignoring the
presence of model selection). It does not fulfill (\ref%
{eq:PRED:PoSI:coverage:aim}) as can be seen from the numerical results in
Section \ref{section:simulation:study}, which is in line with the related
results in \cite{leeb13various}.

\subsection{The various confidence intervals \label{CI}}

For the construction of the quantities $K(x_{0},M)$ we distinguish two cases
regarding the observation on $x_{0}$: (i) The vector $x_{0}$ is observed in
its entirety (regardless of which model $\hat{M}$ is selected), or (ii) only
the subvector $x_{0}{[\hat{M}]}$ of $x_{0}$ is observed (note that only this
subvector is needed for the computation of the post-model-selection
predictor $x_{0}^{\prime }{[\hat{M}]}\hat{\beta}_{\hat{M}}$). The former
case will arise if measuring all the components of $x_{0}$ is not too
costly, whereas the latter case will be relevant in practical situations
where the selected model is determined first and then only observations for $%
x_{0}{[\hat{M}]}$ (and not for the other components of $x_{0}$) are
collected, e.g., out of cost considerations. For example, in a medical
application one may want to avoid measuring prognostic variables that
require invasive procedures or that incur high monetary costs, see, e.g., 
\cite{Hepatol}. Cost considerations in the context of model selection or
prediction are also common in fields such as industrial process control or
engineering (\cite{jaupi}, \cite{souders}).

For the case (i), where $x_{0}$ is entirely observed, the following
straightforward adaptation of the approach in \cite{berk13valid} yields a
constant $K_{1}(x_{0})=K_{1}(x_{0},r)=K_{1}(x_{0},r,X,\alpha ,\mathcal{M})$
(not depending on ${M}$) such that the resulting confidence interval (\ref%
{eq:general:form:CI}) satisfies (\ref{eq:PRED:PoSI:coverage:aim}): Observe
that%
\begin{equation}
x_{0}^{\prime }{[\hat{M}]}\hat{\beta}_{\hat{M}}-x_{0}^{\prime }{[\hat{M}]}%
\beta _{\hat{M}}^{(n)}=s_{\hat{M}}^{\prime }\left( Y-\mu \right) ,
\label{eq:betaM:moins:betaMn}
\end{equation}%
define $\bar{s}_{M}=s_{M}/\left\Vert s_{M}\right\Vert $ if $s_{M}\neq 0$,
and set $\bar{s}_{M}=0\in \mathbb{R}^{n}$ if $s_{M}=0$. Then obviously we
have the upper bound%
\begin{equation}
\left\vert \bar{s}_{\hat{M}}^{\prime }\left( Y-\mu \right) \right\vert /\hat{%
\sigma}\leq \max_{M\in \mathcal{M}}\left\vert \bar{s}_{M}^{\prime }\left(
Y-\mu \right) \right\vert /\hat{\sigma}.
\label{eq:principle:upper:bound:PoSI}
\end{equation}%
Define $K_{1}(x_{0})$ to be the smallest constant satisfying%
\begin{equation}
P_{n,\mu ,\sigma }\left( \max_{M\in \mathcal{M}}\left\vert \bar{s}%
_{M}^{\prime }\left( Y-\mu \right) \right\vert /\hat{\sigma}\leq
K_{1}(x_{0})\right) \geq 1-\alpha .  \label{eq:def:posi:xzero}
\end{equation}%
It is important to note that the probability on the left-hand side of the
preceding display neither depends on $\mu $ nor on $\sigma $; it also
depends on the estimator $\hat{\sigma}$ only through the `degrees of
freedom' parameter $r$: To see this note that $\bar{s}_{M}^{\prime }\left(
Y-\mu \right) =\bar{s}_{M}^{\prime }P_{X}\left( Y-\mu \right) $, since $\bar{%
s}_{M}$ belongs to the column space of $X$. Consequently, the collection of
all the quantities $\bar{s}_{M}^{\prime }\left( Y-\mu \right) $ is jointly
distributed as $N(0,\sigma ^{2}C)$, independently of $\hat{\sigma}^{2}\sim
\left( \sigma ^{2}/r\right) \chi ^{2}\left( r\right) $, where the covariance
matrix $C$ depends only on $x_{0}$ and $X$. Hence the joint distribution of
the collection of ratios $\left\vert \bar{s}_{M}^{\prime }\left( Y-\mu
\right) \right\vert /\hat{\sigma}$ does neither depend on $\mu $ nor $\sigma 
$, and depends on the estimator $\hat{\sigma}$ only through $r$. It is now
plain that $K_{1}(x_{0})$ only depends on $x_{0}$, $r$, $X$, $\alpha $, and $%
\mathcal{M}$. Furthermore note that $K_{1}(x_{0})=0$ in case $x_{0}=0$;
otherwise, $K_{1}(x_{0})$ is positive, equality holds in (\ref%
{eq:def:posi:xzero}), and $K_{1}(x_{0})$ is the unique $\left( 1-\alpha
\right) $-quantile of the distribution of the upper bound in (\ref%
{eq:principle:upper:bound:PoSI}). [This follows from Lemma \ref{lem:a1} in
Appendix \ref{app A} and from the observation that, in view of our
assumptions on $\mathcal{M}$, $\bar{s}_{M}^{\prime }=0$ for all $M\in 
\mathcal{M}$ holds if and only if $x_{0}=0$.] Furthermore, observe that $%
K_{1}(x_{0})$ coincides with a PoSI1 constant of \cite{berk13valid} in case $%
x_{0}$ is one of the standard basis vectors $e_{i}$. [This can be seen by
comparison with (4.14) in \cite{berk13valid} and noting that the maximum
inside the probability in (\ref{eq:def:posi:xzero}) effectively extends only
over models satisfying $i\in M$, since $\bar{s}_{M}=0$ holds for models $M$
with $i\notin M$ if $x_{0}=e_{i}$.] Finally, $K_{naive}\leq K_{1}\left(
x_{0}\right) $ clearly holds provided $x_{0}\neq 0$ (since $\bar{s}%
_{M}^{\prime }\left( Y-\mu \right) /\hat{\sigma}$ follows Student's
t-distribution with $r$ degrees of freedom if $s_{M}\neq 0$).

As a consequence of (\ref{eq:principle:upper:bound:PoSI}) and the discussion
in the preceding paragraph we thus immediately obtain the following
proposition.

\begin{prop}
\label{prop:PoSI:coverage} Let $\hat{M}$ be an arbitrary model selection
procedure with values in $\mathcal{M}$, let $x_{0}\in \mathbb{R}^{p}$ be
arbitrary, and let $K_{1}(x_{0})$ be defined by (\ref{eq:def:posi:xzero}).
Then the confidence interval (\ref{eq:general:form:CI}) with $K(x_{0},\hat{M}%
)$ replaced by $K_{1}(x_{0})$ satisfies the coverage property (\ref%
{eq:PRED:PoSI:coverage:aim}).
\end{prop}

The coverage in Proposition \ref{prop:PoSI:coverage} is guaranteed for 
\textit{all} model selection procedures with values in $\mathcal{M}$, and
thus leads to `universally valid post-selection inference' in case $\mathcal{%
M}$ is chosen to be the set of all full-rank submodels obtainable from $X$
(enlarged by the empty set and provided $X$ does not have a zero column);
cf. \cite{berk13valid}, where similar guarantees are obtained for the
components of $\beta _{\hat{M}}^{(n)}$. [In fact, the construction of $%
K_{1}(x_{0})$ implies that the collection of intervals $x_{0}^{\prime }{[M]}%
\hat{\beta}_{M}\pm K_{1}(x_{0})||s_{{M}}||\hat{\sigma}$ with $M\in \mathcal{M%
}$ provides a simultaneous confidence band for $x_{0}^{\prime }{[M]}\beta
_{M}^{(n)}$.]

Consider next case (ii) where only the components of $x_{0}{[\hat{M}]}$ are
observed. In this case, the confidence interval of Proposition \ref%
{prop:PoSI:coverage} is not feasible in that it cannot be computed in
general, because $K_{1}(x_{0})$ will depend on \emph{all} components of $%
x_{0}$ (and not only on those appearing in $x_{0}{[\hat{M}]}$) due the
maximum figuring in (\ref{eq:def:posi:xzero}) and our assumptions on $%
\mathcal{M}$. A first solution is to define 
\begin{equation}
K_{2}(x_{0}{[M]},{M)}=\sup \left\{ K_{1}(x):x{[M]=}x_{0}{[M]}\right\} ,
\label{eq:max:of:PoSI}
\end{equation}%
and then to use the confidence interval (\ref{eq:general:form:CI}) with $%
K(x_{0},\hat{M})$ replaced by $K_{2}(x_{0}{[\hat{M}],\hat{M})}$. Note that $%
K_{2}(x_{0}{[M]},{M)}$, and hence the corresponding confidence interval,
depends on $x_{0}$ only via $x_{0}{[M]}$, and thus can be computed in case
(ii). Of course, $K_{2}(x_{0}{[M]},{M)}$ also depends on $r$, $X$, $\alpha $%
, and $\mathcal{M}$, and we shall write $K_{2}(x_{0}{[M],M,r)}$ if we want
to stress dependence on $r$. It is easy to see that $K_{2}(x_{0}{[M],M)}$ is
finite (as it is not larger than the Scheff\'{e} constant as we shall see
below). Because $K_{2}(x_{0}{[M],M)}$ is never smaller than $K_{1}(x_{0})$,
we have the following corollary to Proposition \ref{prop:PoSI:coverage}.

\begin{cor}
\label{cor:max:of:PoSI:coverage} Let $\hat{M}$ be an arbitrary model
selection procedure with values in $\mathcal{M}$, let $x_{0}\in \mathbb{R}%
^{p}$ be arbitrary, and let $K_{2}(x_{0}{[M],M)}$ be defined by (\ref%
{eq:max:of:PoSI}). Then the confidence interval (\ref{eq:general:form:CI})
with $K(x_{0},\hat{M})$ replaced by $K_{2}(x_{0}{[\hat{M}],\hat{M})}$
satisfies the coverage property (\ref{eq:PRED:PoSI:coverage:aim}).
\end{cor}

The computation of $K_{2}(x_{0}{[\hat{M}],\hat{M})}$ is more costly than
that of $K_{1}(x_{0})$. Indeed, it requires to embed the algorithm for
computing $K_{1}(x_{0})$ in an optimization procedure. Thus, for the cases
where the resulting computational cost is prohibitive, we present in the
subsequent proposition larger constants $K_{3}(x_{0}{[\hat{M}],\hat{M})}$, $%
K_{4}$, and $K_{5}$ that are simpler to compute. Algorithms for computing
these constants are discussed in Appendix \ref{section:practical:algorithms}%
. The constant $K_{4}$ is obtained by applying a union bound to (\ref%
{eq:def:posi:xzero}), whereas $K_{3}$ is obtained by applying a more refined
`partial' union bound. [More precisely, for $M\in \mathcal{M}$ the
complement of the probability in (\ref{eq:def:posi:xzero}) (with $%
K_{1}(x_{0})$ replaced by a generic variable $t$) can be expressed as in (%
\ref{expression}) in Appendix \ref{app A}. For given $M\in \mathcal{M}$, and
after conditioning on the variance estimator (represented by $G$ there), we
apply a union bound by decomposing the maximum over $\mathcal{M}$ into a
maximum over the submodels of the given $M$ and a maximum over the models
not nested in $M$. A further union bound is applied to the latter group of
models, giving rise to the bound (\ref{bound_2}) in Appendix \ref{app A}.
Inspection of this bound shows that the probability appearing in (\ref{eq:F}%
) below springs from the submodels of $M$, whereas the models not nested in $%
M$ give rise to the term in (\ref{eq:F}) involving the $Beta$-distribution
function.]

For $x_{0}\in \mathbb{R}^{p}$ and $M\in \mathcal{M}$ define now the
distribution function $F_{M,x_{0}}^{\ast }$ for $t\geq 0$ via%
\begin{equation}
F_{M,x_{0}}^{\ast }\left( t\right) =1-\min \left[ 
\begin{array}{c}
1,\Pr \left( \max_{M_{\ast }\in \mathcal{M},M_{\ast }\subseteq M}\left\vert 
\bar{s}_{M_{\ast }}^{\prime }V\right\vert >t\right) \\ 
+c\left( M,\mathcal{M}\right) \left( 1-F_{Beta,1/2,(d-1)/2}\left(
t^{2}\right) \right)%
\end{array}%
\right]  \label{eq:F}
\end{equation}%
and via $F_{M,x_{0}}^{\ast }\left( t\right) =0$ for $t<0$. Here $c\left( M,%
\mathcal{M}\right) $ denotes the number of models $M_{\ast }\in \mathcal{M}$
that satisfy $M_{\ast }\nsubseteq M$, $V$ is a random vector that is
uniformly distributed on the unit sphere in the column space of $X$, and $%
F_{Beta,1/2,(d-1)/2}$ denotes the $Beta\left( 1/2,(d-1)/2\right) $%
-distribution function, with the convention that in case $d=1$ we use $%
F_{Beta,1/2,0}$ to denote the distribution function of pointmass at $1$. In
view of our assumptions on $\mathcal{M}$ it follows that $c\left( M,\mathcal{%
M}\right) \geq 1$ always holds, except in the case where $M=\{1,...,p\}$
(and when this set belongs to $\mathcal{M}$). Next define the distribution
function $F_{M,x_{0}}$ via%
\begin{equation}
F_{M,x_{0}}\left( t\right) =\mathbb{E}_{G}F_{M,x_{0}}^{\ast }\left(
t/G\right) ,  \label{eq:int_G}
\end{equation}%
where $G$ denotes a nonnegative random variable such that $G^{2}/d$ follows
an $F$-distribution with $(d,r)$-degrees of freedom and $\mathbb{E}_{G}$
represents expectation w.r.t. the distribution of $G$. We stress that $%
F_{M,x_{0}}$ depends on $x_{0}$ only through $x_{0}[M]$, and hence the same
is true for the constant $K_{3}(x_{0}[M],M)$ we define next: For any $%
x_{0}\in \mathbb{R}^{p}$ and any $M\in \mathcal{M}$ define $%
K_{3}(x_{0}[M],M) $ to be the smallest constant $K$ satisfying%
\begin{equation}
F_{M,x_{0}}\left( K\right) \geq 1-\alpha .  \label{eq:def:K3}
\end{equation}%
Furthermore, set $K_{4}=K_{3}(x_{0}{[\varnothing ],\varnothing )}$. Finally, 
$K_{5}$ is the Scheff\'{e} constant, i.e., the $\left( 1-\alpha \right) $%
-quantile of $G$ (\cite{scheffe59analysis}); see the corresponding
discussion in Section 4.8 of \cite{berk13valid}. Recall that $1-\alpha \in
(0,1)$ has been assumed.

\begin{prop}
\label{prop:PoSI:upper:bound} Let $x_{0}\in \mathbb{R}^{p}$ be arbitrary.
Then we have the following:

(a) $K_{3}(x_{0}[M],M)$ exists and is well-defined. If $M=\{1,...,p\}\in 
\mathcal{M}$ and $x_{0}=0$, then $K_{3}(x_{0}[M],M)={0}$ (and $F_{M,x_{0}}$
is the c.d.f. of pointmass at zero). If $M\neq \{1,...,p\}$ or $x_{0}\neq 0$%
, then (i) $0{<K_{3}(x_{0}[M],M)<\infty }$ holds, and (ii) equality holds in
(\ref{eq:def:K3}) if and only if $K=K_{3}(x_{0}[M],M)$.

(b) For every $M\in \mathcal{M}$ we have%
\begin{equation}
K_{2}(x_{0}[M],M)\leq K_{3}(x_{0}[M],M)\leq K_{4}\leq K_{5}.
\label{eq:comparison:length}
\end{equation}%
Furthermore, 
\begin{equation}
{K_{2}(x_{0}[M}_{2}{],M}_{2}{)\leq K_{2}(x_{0}[M}_{1}{],M}_{1}{),}
\label{eq:monoton:1}
\end{equation}%
\begin{equation}
{K_{3}(x_{0}[M}_{2}{],M}_{2}{)\leq K_{3}(x_{0}[M}_{1}{],M}_{1}{)}
\label{eq:monoton:2}
\end{equation}%
hold whenever $M_{1}\subseteq M_{2}$, $M_{i}\in \mathcal{M}$.
\end{prop}

It is obvious that $K_{3}(x_{0}[M],M)$ depends, besides $x_{0}[M]$ and $M$,
only on $r$, $X$, $\alpha $, and $\mathcal{M}$, whereas $K_{4}$ only depends
on $r$, $d$, $\alpha $, and $\mathcal{M}$, and $K_{5}$ depends only on $r$, $%
d$, and $\alpha $. [Like with $K_{1}\left( x_{0}\right) $, also the other
constants introduced depend on the estimator $\hat{\sigma}$ only through $r$%
.] We shall write $K_{3}(x_{0}[M],M,r)$, $K_{4}\left( r\right) $, and $%
K_{5}\left( r\right) $ if we want to stress dependence on $r$. Note that $%
K_{1}\left( x_{0}\right)
=K_{3}(x_{0}[M_{full}],M_{full})=K_{3}(x_{0},M_{full})$, provided $%
M_{full}:=\{1,...,p\}$ belongs to $\mathcal{M}$, and that $%
K_{3}(x_{0}[M],M)=K_{4}$ holds for any $M\in \mathcal{M}$ satisfying $%
\left\vert M\right\vert =1$ and $\bar{s}_{M}\neq 0$. [Indeed, in this case,
the probability appearing in (\ref{eq:F}) equals $1-F_{Beta,1/2,(d-1)/2}%
\left( t^{2}\right) $ as can be seen from the proof of Proposition \ref%
{prop:PoSI:upper:bound}.] Similarly, $K_{3}(x_{0}[M],M)=K_{4}$ holds for any 
$M\in \mathcal{M}$ in case $d=1$ as is not difficult to see. The proof of
the inequalities involving the constants $K_{3}$ and $K_{4}$ in the above
proposition is an extension of an argument in \cite{berk13validold} (not
contained in the published version \cite{berk13valid}) to find -- in the
case $p=d$ -- an upper-bound for their PoSI constant that does not depend on 
$X$, but only on $d$. [Note that $K_{4}$ is a counterpart to $K_{univ}$ in 
\cite{berk13validold}.]\ Inequalities (\ref{eq:monoton:1}) and (\ref%
{eq:monoton:2}) simply reflect the fact that observing only $x_{0}[M]$
implies that fewer information about $x_{0}$ is provided for smaller models $%
M$. As a consequence of these inequalities it is possible that, on the event
where a small model $M_{1}$ is selected, the resulting confidence interval
is larger than it is on the event where a larger model $M_{2}$ is selected.
Again, this simply reflects the fact that less information on $x_{0}$ is
available under the smaller model. Note, however, that the just discussed
phenomenon is counteracted by the fact that the length of the confidence
interval also depends on $||s_{M}||$ and that we have $||s_{M_{1}}||\leq
||s_{M_{2}}||$ for $M_{1}\subseteq M_{2}$; cf. Figure 1 in Section \ref%
{section:simulation:study}.

Proposition \ref{prop:PoSI:upper:bound} implies that (\ref%
{eq:comparison:length}) holds with $\hat{M}$ replacing $M$, which together
with Corollary \ref{cor:max:of:PoSI:coverage} immediately implies the
following result. We stress that the confidence intervals figuring in the
subsequent corollary depend on $x_{0}$ only through $x_{0}{[\hat{M}]}$ and
thus are feasible in case (ii) discussed at the beginning of Section \ref{CI}%
.

\begin{cor}
\label{cor:PoSI:upper:bound}Let $\hat{M}$ be an arbitrary model selection
procedure with values in $\mathcal{M}$, and let $x_{0}\in \mathbb{R}^{p}$ be
arbitrary. Then the confidence interval (\ref{eq:general:form:CI}) with $%
K(x_{0},\hat{M})$ replaced by $K_{3}(x_{0}{[\hat{M}],\hat{M})}$ ($K_{4}$, or 
$K_{5}$, respectively) satisfies the coverage property (\ref%
{eq:PRED:PoSI:coverage:aim}).
\end{cor}

We conclude this section with a few remarks regarding extensions.

\begin{rem}
\label{rem:kmown:variance:case}\normalfont\emph{(Infeasible variance
estimators) }(i) For later use we note that all results derived in Section %
\ref{section:design:dependent} continue to hold if $\hat{\sigma}^{2}$ is
allowed to also depend on $\sigma $ but otherwise satisfies the assumptions
made earlier (e.g., if $\hat{\sigma}^{2}=\sigma ^{2}Z/r$ where $Z$ is an
observable chi-square distributed random variable with $r$ degrees of
freedom that is independent of $P_{X}Y$).

(ii) If we set $\hat{\sigma}^{2}=\sigma ^{2}$ and $r=\infty $, all of the
results derived in Section \ref{section:design:dependent} continue to hold
with obvious modifications. In particular, in Proposition \ref%
{prop:PoSI:upper:bound} the random variable $G^{2}$ then follows a
chi-squared distribution with $d$ degrees of freedom. We shall denote the
constants corresponding to $K_{1}(x_{0})$, $K_{2}(x_{0}[M],M)$, $%
K_{3}(x_{0}[M],M)$, $K_{4}$, and $K_{5}$ obtained by setting $\hat{\sigma}%
^{2}=\sigma ^{2}$ and $r=\infty $ by $K_{1}(x_{0},\infty )$, etc. We stress
that these constants do \emph{not} depend on $\sigma $.
\end{rem}

\begin{rem}
\normalfont\label{rem_2.6}(i) All results carry over immediately to the case
where $\mu $ can vary only in a subset $\mathfrak{M}$ of $\mathbb{R}^{n}$.

(ii) We have assumed that any non-empty $M\in \mathcal{M}$ is of full-rank.
This assumption could easily be dropped, but this would lead to more
unwieldy results.

(iii) Since the development in Section \ref{section:design:dependent} is
based on the bound (\ref{eq:principle:upper:bound:PoSI}), it is obvious that
all results in Section \ref{section:design:dependent} also hold if $\hat{M}=%
\hat{M}(X,Y,\bar{\sigma}^{2})$ for some arbitrary estimator $\bar{\sigma}%
^{2} $, that may differ from the estimator $\hat{\sigma}^{2}$ that governs
the length of the confidence intervals considered.
\end{rem}

\subsection{On the merits of the non-standard targets}

\begin{rem}
\label{rem:interpretation}\normalfont\emph{\ }(i) As already noted, the
(non-standard) coverage target in \cite{berk13valid} is $\beta _{\hat{M}%
}^{(n)}$ (where these authors choose to represent it in what they call `full
model indexing'). While $\beta _{\hat{M}}^{(n)}$ has a clear technical
meaning as the coefficient vector that provides the best approximation of $%
\mu $ by elements of the form $X[\hat{M}]\gamma $ w.r.t. the Euclidean
distance, adopting this quantity as the target for inference confronts one
with the fact that the target then depends on the data $Y$ via $\hat{M}$
(implying that the target as well as its dimension are random); furthermore,
different model selection procedures give rise to different targets $\beta _{%
\hat{M}}^{(n)}$. Also note that, e.g., the meaning of the first component of
the target $\beta _{\hat{M}}^{(n)}$ depends on the selected model ${\hat{M}}$%
. The target $x_{0}^{\prime }{[\hat{M}]}\beta _{\hat{M}}^{(n)}$ considered
in this paper, while again being random and sharing many of the properties
of $\beta _{\hat{M}}^{(n)}$ just mentioned, seems to be somewhat more
amenable to interpretation since it is simply the random convex combination $%
\sum_{M}x_{0}^{\prime }{[M]}\beta _{M}^{(n)}\boldsymbol{1}({\hat{M}=M})$ of
the (infeasible) predictors $x_{0}^{\prime }{[M]}\beta _{M}^{(n)}$ (which
one would typically use if model $M$ is forced upon one for prediction and
which all have one and the same dimension, not depending on the data).

(ii) In the classical case, i.e., when $\mu =X\beta $ and $d=p\leq n$, one
can justly argue that the target for inference should be $x_{0}^{\prime
}\beta $ rather than $x_{0}^{\prime }{[\hat{M}]}\beta _{\hat{M}}^{(n)}$
because $x_{0}^{\prime }\beta $ is a better (infeasible) predictor in the
mean-squared error sense than is $x_{0}^{\prime }{[\hat{M}]}\beta _{\hat{M}%
}^{(n)}$ provided $y_{0}$ is independent of ${\hat{M}}$ (which will
certainly be the case if $y_{0}$ is independent of $Y$ and $\hat{\sigma}^{2}$%
, or if $y_{0}$ is independent of $Y$ and ${\hat{M}}$ is only a function of $%
X$ and $Y$). [This is so since the mean-squared error of prediction of $%
x_{0}^{\prime }\beta $ is not larger than the one of $x_{0}^{\prime }{[M]}%
\beta _{M}^{(n)}$ for every $M$ and since ${\hat{M}}$ is independent of $%
y_{0}$.] However, this argument does not apply if $x_{0}$ is not observed in
its entirety, but only $x_{0}{[\hat{M}]}$ is observed, because then $%
x_{0}^{\prime }\beta $ is not available. In this case we thus indeed have
some justification for the target $x_{0}^{\prime }{[\hat{M}]}\beta _{\hat{M}%
}^{(n)}$ even in the classical case. This is in contrast with the situation
when, as in \cite{berk13valid}, one's interest focusses on parameters rather
than predictors: Similar as before one can argue that in the classical case
the true parameter $\beta $ should be the target rather than $\beta _{\hat{M}%
}^{(n)}$ but there seems now to be little to justify the non-standard target 
$\beta _{\hat{M}}^{(n)}$ (as the preceding argument justifying the target $%
x_{0}^{\prime }{[\hat{M}]}\beta _{\hat{M}}^{(n)}$ even in the classical case
is obviously not applicable to the target $\beta _{\hat{M}}^{(n)}$).

(iii) In view of the preceding discussion it seems that the non-standard
target $\beta _{\hat{M}}^{(n)}$ of \cite{berk13valid} mainly has a
justification in a non-classical setting where $\mu $ is not assumed to
belong to the column space of $X$ (implying $d<n$), or where $d<p$ holds
(subsuming in particular the important case $p>n=d$), because in these cases 
$\beta $ is no longer available as a target (being not defined or not
uniquely defined). However, in a setting, where $\mu $ is not assumed to
belong to the column space of $X$ or where $p>n=d$ holds, the assumption on
the variance estimator $\hat{\sigma}^{2}$ made in \cite{berk13valid} (as
well as in the present paper) becomes problematic and quite restrictive; see
Remark 2.1(ii) in \cite{leeb13various} as well as Appendix \ref%
{Section_on_sigma}. Hence, there is some advantage in considering the
targets $x_{0}^{\prime }{[\hat{M}]}\beta _{\hat{M}}^{(n)}$ rather than $%
\beta _{\hat{M}}^{(n)}$ as the former has a justification in the classical
as well as in the non-classical framework.

(iv) We note the obvious fact that if the target of inference is the
standard target $x_{0}^{\prime }\beta $ (assuming the classical case) then
the reasoning underlying Proposition \ref{prop:PoSI:coverage} does not apply
since the difference between the post-model-selection predictor and the
standard target is not independent of $\beta $. For the same reason the
approach in \cite{berk13valid} cannot provide a solution to the problem of
constructing confidence sets for the standard target $\beta $.
\end{rem}

\begin{rem}
\label{rem:predictor}\normalfont\emph{(On the optimality of the
design-dependent target) }(i) The infeasible predictor $x_{0}^{\prime }{[M]}%
\beta _{M}^{(n)}$ (for fixed $M$) is the best predictor for $y_{0}$ in the
mean-squared error sense among all predictors of the form $x_{0}^{\prime }{%
[M]}\gamma $ in case $\left. y_{0}\right\vert \nu ,x_{0}\sim N(\nu ,\sigma
^{2})$ and $\left( \nu ,x_{0}^{\prime }\right) $ is drawn from the empirical
distribution of $\left( \mu _{i},x_{i}^{\prime }\right) $ where $%
x_{i}^{\prime }$ denotes the $i$-th row of $X$ (`in-sample prediction').
[More generally, this is so if $\left( \nu ,x_{0}^{\prime }\right) $ is
drawn from the empirical distribution of $\left( \mu
_{i}+a_{i},x_{i}^{\prime }\right) $ where $a$ is a fixed vector orthogonal
to the column space of $X$.] Otherwise, it does in general not have this
optimality property (but nevertheless its feasible counterpart $%
x_{0}^{\prime }{[M]}\hat{\beta}_{M}$ would typically be used if one is
forced to base prediction on model $M$).

(ii) The optimality property in (i) carries over to the design-dependent
target $x_{0}^{\prime }{[\hat{M}]}\beta _{\hat{M}}^{(n)}$ provided $%
(y_{0},x_{0}^{\prime })^{\prime }$ is independent of ${\hat{M}}$.
\end{rem}

\subsection{Behavior of the constants $K_{i}$ as a function of $p$ \label%
{sub2.3}}

In this section we provide some results on the size of the constants $K_{i}$
that govern the length of the confidence intervals. In particular, these
results help in answering the question how tight a bound for $K_{1}$ and $%
K_{2}$ is provided by $K_{3}$ or $K_{4}$.

\subsubsection{Orthogonal designs}

\cite{berk13valid} show that in the case $p=d\leq n$ their PoSI constant
becomes smallest for the case of orthogonal design (provided the model
universe $\mathcal{M}$ is sufficiently rich, e.g., $\mathcal{M}$ contains
all submodels) and then has rate $\sqrt{\log p}$ as $p\rightarrow \infty $,
at least in the known-variance case; cf. Proposition 5.5 in \cite%
{berk13valid} (where the error term $o(d)$ given in this result should read $%
o(1)$). In the next proposition we study the order of magnitude of $%
K_{1}(x_{0})$, the analogue of the PoSI constant and of the closely related
constant $K_{2}(x_{0}[M],M)$ in the case of orthogonal design. Recall that $%
K_{1}(x_{0})$ is only feasible if $x_{0}$ is observed in its entirety, while 
$K_{2}(x_{0}[M],M)$ is the ideal bound for $K_{1}(x_{0})$ given only
knowledge of $x_{0}[M]$. Note that in the following result some of the
objects depend on $p$, but we do not always show this in the notation.
Furthermore, $\phi $ and $\Phi $ denote the p.d.f. and c.d.f. of a standard
normal variable, respectively, and $\left\Vert x\right\Vert _{0}$ denotes
the $l_{0}$-norm.

\begin{prop}
\label{constants_1}Consider the known-variance case (i.e., $r=\infty $ and $%
\hat{\sigma}^{2}=\sigma ^{2}$) and assume that for every $p\geq 1$ the model
universe $\mathcal{M}$ used is the power set of $\{1,...,p\}$. Let $\alpha $%
, $0<\alpha <1$, be given, not depending on $p$.

(a) For any $p\geq 1$ let $X=X(p)$ be an $n(p)\times p$ matrix with
(non-zero) orthogonal columns. For any such sequence $X$ one can find a
corresponding sequence of (non-zero) $p\times 1$ vectors $x_{0}$ such that $%
K_{1}(x_{0},\infty )=K_{1}(x_{0},\infty ,X,\alpha ,\mathcal{M})$ satisfies%
\begin{equation*}
\liminf_{p\rightarrow \infty }K_{1}(x_{0},\infty )/\sqrt{p}\geq \xi
\end{equation*}%
where $\xi =\sup_{b>0}\phi (b)/\sqrt{1-\Phi (b)}\approx 0.6363$.
Furthermore, for any sequence $X$ as above one can find another sequence of
(non-zero) $p\times 1$ vectors $x_{0}$ such that $K_{1}(x_{0},\infty )=O(1)$
(for example, any sequence of (non-zero) $p\times 1$ vectors $x_{0}$
satisfying $\sup_{p}\left\Vert x_{0}\right\Vert _{0}<\infty $ will do).

(b) Let $\gamma \in \lbrack 0,1)$ be given. Then $K_{2}(x_{0}[M],M,\infty
)=K_{2}(x_{0}[M],M,\infty ,X,\alpha ,\mathcal{M})$ satisfies%
\begin{equation*}
\liminf_{p\rightarrow \infty }\inf_{x_{0}\in \mathbb{R}^{p}}\inf_{X\in 
\mathsf{X}(p)}\inf_{M\in \mathcal{M},|M|\leq \gamma
p}K_{2}(x_{0}[M],M,\infty )/\sqrt{p}\geq \xi \sqrt{1-\gamma },
\end{equation*}%
where $\mathsf{X}(p)=\bigcup_{n\geq p}\left\{ X:X\text{ is }n\times p\text{
with non-zero orthogonal columns}\right\} $.
\end{prop}

The lower bounds given in the preceding proposition clearly also apply to $%
K_{3}(x_{0}[M],M,\infty )$ and $K_{4}(\infty )$ a fortiori. Part (a) of the
above proposition shows that, even in the orthogonal case, the growth of $%
K_{1}(x_{0},\infty )$ is -- in the worst-case w.r.t. $x_{0}$ -- of the order 
$\sqrt{p}$. This is in contrast to the above mentioned result of \cite%
{berk13valid} for the PoSI constant. Part (a) also shows that there are
other choices for $x_{0}$ such that $K_{1}(x_{0},\infty )$ stays bounded. In
this context also recall that $K_{1}(x_{0},\infty )$ with $x_{0}$ equal to a 
$p\times 1$ standard basis vector coincides with a PoSI1 constant and thus
equals the $(1-\alpha )$-quantile of the distribution of the absolute value
of a standard normal variable in the orthogonal case. Part (b) goes on to
show that regardless of $x_{0}$ and $X$ the growth of the constants $%
K_{2}(x_{0}[M],M,\infty )$ is of the order $\sqrt{p}$ (except perhaps for
very large submodels $M$).

\subsubsection{Order of magnitude of $K_{3}$ and $K_{4}$}

The next proposition, which exploits results in \cite{zhang13rank}, shows
that $K_{4}(\infty )$ is a tight upper bound for $K_{3}(x_{0}[M],M,\infty )$
at least if $p$ is large. It also provides the growth rates for $%
K_{4}(\infty )$ and $K_{3}(x_{0}[M],M,\infty )$. As before, the dependence
of several objects on $p$ (or $n$) will not always be shown in the notation.
For the following recall the constants $c\left( M,\mathcal{M}\right) $
defined after (\ref{eq:F}).

\begin{prop}
\label{constants_2}Consider the known-variance case (i.e., $r=\infty $ and $%
\hat{\sigma}^{2}=\sigma ^{2}$) and assume that for every $p\geq 1$ a
(non-empty) model universe $\mathcal{M}=\mathcal{M}_{p}$ is given that
satisfies (i) $\bigcup \left\{ M:M\in \mathcal{M}\right\} =\left\{
1,...,p\right\} $, (ii) $\varnothing \in \mathcal{M}$, (iii) $c\left( M,%
\mathcal{M}\right) \geq \tau \left\vert \mathcal{M}\right\vert $ for every $%
M\in \mathcal{M}$ with $M\neq \left\{ 1,\ldots ,p\right\} $, where $\tau >0$
is a given number (neither depending on $M$, $\mathcal{M}$, nor $p$), and
(iv) $\left\vert \mathcal{M}\right\vert \rightarrow \infty $ as $%
p\rightarrow \infty $. For $n\in \mathbb{N}$, the set of positive integers,
let $\mathsf{X}_{n,p}(\mathcal{M})$ denote the set of all $n\times p$
matrices of rank $\min (n,p)$ with the property that $X[M]$ has full
column-rank for every $\varnothing \neq M\in \mathcal{M}$. Furthermore, let $%
\alpha $, $0<\alpha <1$, be given (neither depending on $p$ nor $n$). Let $%
n(p)\in \mathbb{N}$ be a sequence such that $n(p)\rightarrow \infty $ for $%
p\rightarrow \infty $ and such that $\mathsf{X}_{n(p),p}(\mathcal{M})\neq
\varnothing $ for every $p\geq 1$. Then we have%
\begin{equation}
\lim_{p\rightarrow \infty }\sup_{M\in \mathcal{M},M\neq
\{1,...,p\}}\sup_{x_{0}\in \mathbb{R}^{p}}\sup_{X\in \mathsf{X}_{n(p),p}(%
\mathcal{M})}|1-(K_{3}(x_{0}[M],M,\infty )/K_{4}(\infty ))|=0,
\label{Closeness of K3 and K4}
\end{equation}%
where $K_{3}(x_{0}[M],M,\infty )=K_{3}(x_{0}[M],M,\infty ,X,\alpha ,\mathcal{%
M})$ and $K_{4}(\infty )=K_{4}(\infty ,\min (n(p),p),\alpha ,\mathcal{M})$.
Furthermore, 
\begin{equation*}
K_{4}(\infty )/\sqrt{\min (n(p),p)\left( 1-\left\vert \mathcal{M}\right\vert
^{-2/\left( \min (n(p),p)-1\right) }\right) }\rightarrow 1
\end{equation*}%
as $p\rightarrow \infty $.
\end{prop}

\begin{rem}
\normalfont(i) $\mathsf{X}_{n(p),p}(\mathcal{M})\neq \varnothing $ implies $%
\mathsf{X}_{n,p}(\mathcal{M})\neq \varnothing $ for $n\geq n(p)$.

(ii) $\mathsf{X}_{n(p),p}(\mathcal{M})$ is certainly non-empty for $n(p)\geq
p$, but -- depending on $\mathcal{M}$ -- this can already be true for $n(p)$
much smaller than $p$.
\end{rem}

The assumptions (i)-(iv) on $\mathcal{M}$ in the preceding proposition are
shown in the next corollary to be always satisfied in the important case
where $\mathcal{M}$ is of the form $\left\{ M\subseteq \left\{ 1,\ldots
,p\right\} :\left\vert M\right\vert \leq m_{p}\right\} $. Furthermore, in
the special case where $\mathcal{M}$ is the universe of all submodels, a
simple formula for the growth rate of $K_{4}(\infty )$ is found.

\begin{cor}
\label{constants_3}Consider the known-variance case (i.e., $r=\infty $ and $%
\hat{\sigma}^{2}=\sigma ^{2}$) and let $\alpha $, $0<\alpha <1$, be given
(neither depending on $p$ nor $n$). Let $m_{p}\in \mathbb{N}$ satisfy $1\leq
m_{p}\leq p$ for every $p\geq 1$ and define the set $\mathcal{M}%
(m_{p})=\left\{ M\subseteq \left\{ 1,\ldots ,p\right\} :\left\vert
M\right\vert \leq m_{p}\right\} $. Then $\mathcal{M}(m_{p})$ satisfies
(i)-(iv) in Proposition \ref{constants_2} with $\tau =1/3$. Consequently,
for $n(p)$ as in Proposition \ref{constants_2}, (\ref{Closeness of K3 and K4}%
) holds with $\mathcal{M}$ replaced by $\mathcal{M}(m_{p})$ and 
\begin{equation*}
K_{4}(\infty )/\sqrt{\min (n(p),p)\left( 1-\left( \sum_{k=0}^{m_{p}}\binom{p%
}{k}\right) ^{-2/\left( \min (n(p),p)-1\right) }\right) }\rightarrow 1
\end{equation*}%
as $p\rightarrow \infty $. In particular, if $m_{p}=p$ for all $p\geq 1$, we
necessarily have $n(p)\geq p$ and%
\begin{equation*}
K_{4}(\infty )/\sqrt{p}\rightarrow \sqrt{3}/2
\end{equation*}%
as $p\rightarrow \infty $.
\end{cor}

In the important case, where $p=d\leq n$ and $\mathcal{M}$ is the entire
power set of $\{1,...,p\}$, the preceding corollary shows that $K_{4}(\infty
)$ (and hence a fortiori all the constants $K_{1}(x_{0},\infty )$,..., $%
K_{3}(x_{0}[M],M,\infty )$) are `bounded away' from the Scheff\'{e} constant 
$K_{5}$ which clearly satisfies $K_{5}/\sqrt{p}\rightarrow 1$ for $%
p\rightarrow \infty $. This is in line with a similar finding in \cite%
{berk13valid}, Section 6.3, for their PoSI constant.

\begin{rem}
\label{rem:sharp:union:bound}\normalfont In the proof of Proposition \ref%
{prop:PoSI:upper:bound} union bounds were used to obtain the results for $%
K_{3}(x_{0}[M],M)$ and $K_{4}$. Hence, one might ask whether or not these
constants as bounds for $K_{2}(x_{0}[M],M)$ are overly conservative. We now
collect evidence showing that improving $K_{3}(x_{0}[M],M)$ and $K_{4}$ will
not be easy and is sometimes impossible: First, Lemma \ref{lem:case:p:egal:2}
in Appendix \ref{app A} shows that there exist $n\times p$ design matrices $%
X $ with $p=d=2$ and vectors $x_{0}$ such that $K_{4}=K_{1}\left(
x_{0}\right) $ in case $\mathcal{M}$ is the universe of all submodels.
Hence, in this case the union bounds used in the proof of Proposition \ref%
{prop:PoSI:upper:bound} are all exact. Furthermore, in the known-variance
case with $p=d\leq n$ and where $\mathcal{M}$ again is the universe of all
submodels, the propositions given above entail that $K_{4}(\infty )\sim 
\sqrt{p}\sqrt{3}/2\approx 0.866\sqrt{p}$ while $K_{1}(x_{0},\infty )\succeq
\xi \sqrt{p}$ with $\xi \approx 0.6363$ is possible; e.g., as the worst-case
behavior in the orthogonal case, or with $x_{0}=e_{i}$ and the design
matrices constructed in the proof of Theorem 6.2 in \cite{berk13valid}
(recall that $K_{1}(e_{i},\infty )$ coincides with a PoSI1 constant). This
again shows that there is little room for improving $K_{3}$ and $K_{4}$.
[Further evidence in that direction is provided by the observation that the
proof of Theorem 6.3 in \cite{berk13valid} implies that $K_{1}^{\ast }/\sqrt{%
p}$ tends to $\sqrt{3}/2$ in probability as $p\rightarrow \infty $, where $%
K_{1}^{\ast }$ is an analogue of $K_{1}\left( x_{0},\infty \right) $ that is
obtained from (\ref{eq:def:posi:xzero}) (with $r=\infty $) after replacing
the vectors $\bar{s}_{M}$ by $2^{p}$ independent random vectors, each of
which is uniformly distributed on the unit sphere of the column space of $X$
(and these vectors being independent of $Y$). In other words, if one ignores
the particular structure of the vectors $\bar{s}_{M}$, then the bound $%
K_{4}\left( \infty \right) $ is close to being sharp for large values of $p$%
.]
\end{rem}

\begin{rem}
\label{nonuniformity}\normalfont The results for $p\rightarrow \infty $ in
this subsection as well as the related results in \cite{berk13valid} should
be taken with a grain of salt as they obviously are highly non-uniform
w.r.t. $\alpha $: Note that -- for fixed $n$ and $p$ -- any one of the
constants $K_{i}$ will vary in the entire interval $(0,\infty )$ as $\alpha $
varies in $(0,1)$ (except for degenerate cases), while the limits in the
results in question do not depend on $\alpha $ at all.
\end{rem}

\section{Confidence intervals for the design-independent non-standard target 
\label{section:distribution:dependent}}

In this section we again consider the model (\ref{eq:full:linear:model}),
but now assume that $\mu =X\beta $ for some unknown $\beta \in \mathbb{R}^{p}
$ holds and that the $n\times p$ matrix $X$ is random, with $X$ independent
of $U$, where $U$ again follows an $N\left( 0,\sigma ^{2}I_{n}\right) $%
-distribution with $0<\sigma <\infty $. We also assume that $X$ has full
column rank almost surely (implying $p\leq n$) and that each row of $X$ is
distributed according to a common $p$-dimensional distribution $\mathcal{L}$
(not depending on $n$) with a finite and positive definite matrix of
(uncentered) second moments, which we denote by $\Sigma $. [We shall refer
to the preceding assumptions as the maintained model assumptions of this
section.] Furthermore, we assume again that we have available an estimator $%
\hat{\sigma}^{2}$ such that, conditionally on $X$, $\hat{\sigma}^{2}$ is
independent of $P_{X}Y$ (or, equivalently, of $\hat{\beta}=(X^{\prime
}X)^{-1}X^{\prime }Y$) and is distributed as $\sigma ^{2}/r$ times a
chi-squared distributed random variable with $r$ degrees of freedom ($1\leq
r<\infty $). The collection $\mathcal{M}$ of admissible models will be
assumed to be the power set of $\{1,\ldots ,p\}$ in this section for
convenience, but see Remark \ref{restr_univ} for possible extensions.
Observe that all the results of Section \ref{section:design:dependent}
remain valid in the setup of the present section if formulated conditionally
on $X$ (and if $x_{0}$ is treated as fixed). [Alternatively, if $x_{0}$ is
random but independent of $X$, $U$, and $\hat{\sigma}^{2}$, the same is true
if the results in Section \ref{section:design:dependent} are then
interpreted conditionally on $X$ and $x_{0}$.] The joint distribution of $Y$%
, $X$, and $\hat{\sigma}^{2}$ (and of $\tilde{\sigma}$ appearing below) will
be denoted by $P_{n,\beta ,\sigma }$ (see also Appendix \ref{notation}).

In this section we shall consider asymptotic results for $n\rightarrow
\infty $ but where $p$ is held constant (for an extension to the case where $%
p$ is allowed to diverge with $n$ see Appendix \ref{app ext}). It is thus
important to recall that all estimators, estimated models, etc. depend on
sample size $n$. Also note that $r$ may depend on sample size $n$. We shall
typically suppress these dependencies on $n$ in the notation. Furthermore,
we note that, while not explicitly shown in the notation, the rows of $X$
and $U$ (and thus of $Y$) may depend on $n$. [As the results in Section \ref%
{section:design:dependent} are results for fixed $n$, this trivially also
applies to the results in that section.] However, recall that $\mathcal{L}$,
and hence $\Sigma $, are not allowed to depend on $n$.

If $M_{1}$ and $M_{2}$ are subsets of $\{1,...,p\}$ and if $Q$ is a $p\times
p$ matrix we shall denote by $Q[M_{1},M_{2}]$ the matrix that is obtained
from $Q$ by deleting all rows $i$ with $i\notin M_{1}$ as well as all
columns $j$ with $j\notin M_{2}$; if $M_{1}$ is empty but $M_{2}$ is not, we
define $Q[M_{1},M_{2}]$ to be the $1\times \left\vert M_{2}\right\vert $
zero vector; if $M_{2}$ is empty but $M_{1}$ is not, we define $%
Q[M_{1},M_{2}]$ to be the $\left\vert M_{1}\right\vert \times 1$ zero
vector; and if $M_{1}=M_{2}=\varnothing $ we set $Q[M_{1},M_{2}]=0\in 
\mathbb{R}$.

To motivate the target studied in this section, consider now the problem of
predicting a new variable $y_{0}=x_{0}^{\prime }\beta +u_{0}$ where $x_{0}$, 
$u_{0}$, $X$, and $U$ are independent and $u_{0}\sim N\left( 0,\sigma
^{2}\right) $. For a given model $M\subseteq \{1,...,p\}$ we consider the
(infeasible) predictor $x_{0}^{\prime }[M]\beta _{M}^{(\star )}$ where%
\begin{equation*}
\beta _{M}^{(\star )}=\beta \lbrack M]+\left( \Sigma \lbrack M,M]\right)
^{-1}\Sigma \lbrack M,M^{c}]\beta \lbrack M^{c}],
\end{equation*}%
with the convention that the inverse is to be interpreted as the
Moore-Penrose inverse in case $M=\varnothing $. Note that $x_{0}^{\prime
}[M]\beta _{M}^{(\star )}=0$ if $M=\varnothing $ and that $x_{0}^{\prime
}[M]\beta _{M}^{(\star )}=x_{0}^{\prime }\beta $ if $M=\{1,\ldots ,p\}$. A
justification for considering this infeasible predictor is given in Remark %
\ref{rem:predictor_2} below. For purpose of comparison we point out that,
under the assumption $\mu =X\beta $ maintained in the present section, $%
\beta _{M}^{(n)}$ defined in (\ref{eq:def:betaMn}) can be rewritten as $%
\beta _{M}^{(n)}=\beta \lbrack M]+\left( X{[M]}^{\prime }{X[M]}\right) ^{-1}X%
{[M]}^{\prime }{X[M}^{c}{]}\beta \lbrack M^{c}]$. Given a model selection
procedure $\hat{M}=\hat{M}(X,Y,\hat{\sigma}^{2})$ we define now the
(infeasible) predictor%
\begin{equation*}
x_{0}^{\prime }[\hat{M}]\beta _{\hat{M}}^{(\star )}
\end{equation*}%
as our new target for inference. We call this target the \emph{%
design-independent (non-standard) target} as it does not depend on the
design matrix $X$ beyond its dependence on $\hat{M}$. We discuss its merits
in the subsequent remarks.

\begin{rem}
\label{rem:justification}\normalfont As in Remark \ref{rem:interpretation}%
(ii) one can argue that the target for inference should be $x_{0}^{\prime
}\beta $ rather than $x_{0}^{\prime }{[\hat{M}]}\beta _{\hat{M}}^{(\star )}$
because again $x_{0}^{\prime }\beta $ is a better (infeasible) predictor
than $x_{0}^{\prime }{[\hat{M}]}\beta _{\hat{M}}^{(\star )}$ provided that $%
\left( x_{0}^{\prime },u_{0}\right) $ is independent of ${\hat{M}}$ (which,
in particular, will be the case if $\left( x_{0}^{\prime },u_{0}\right) $ is
independent of $X$, $U$, and $\hat{\sigma}$, or if $\left( x_{0}^{\prime
},u_{0}\right) $ is independent of $X$, $U$ and ${\hat{M}}$ is only a
function of $X$ and $Y$). But again, this argument does not apply if $x_{0}$
is not observed in its entirety, but only $x_{0}{[\hat{M}]}$ is observed.
\end{rem}

\begin{rem}
\label{rem:predictor_2}\normalfont\emph{(On the optimality of the
design-independent target)\ }(i) Assume that additionally $x_{0}^{\prime
}\sim \mathcal{L}$. If we are forced to use the (theoretical) predictors of
the form $x_{0}^{\prime }[M]\gamma $, then straightforward computation shows
that $x_{0}^{\prime }[M]\beta _{M}^{(\star )}$ provides the smallest
mean-squared error of prediction among all the linear predictors $%
x_{0}^{\prime }[M]\gamma $. [Note that this result corresponds to the
observation made in Remark \ref{rem:predictor} with $\mathcal{L}$
corresponding to the empirical distribution of the rows of $X$.] If,
furthermore, $x_{0}$ is normally distributed, then $x_{0}$ and $u_{0}$ are
jointly normal and thus $x_{0}^{\prime }[M]\beta _{M}^{(\star )}$ is the
conditional expectation of $y_{0}$ given $x_{0}[M]$ and hence is also the
best predictor in the class of all predictors depending only on $x_{0}[M]$.

(ii) Again assume that $x_{0}^{\prime }\sim \mathcal{L}$. The discussion in
(i) implies that $x_{0}^{\prime }{[\hat{M}]}\beta _{\hat{M}}^{(\star )}$ has
a mean-squared error of prediction not larger than the one of $x_{0}^{\prime
}{[\hat{M}]}\gamma (\hat{M})$ for any choice of $\gamma (\hat{M})$, provided 
$\left( x_{0}^{\prime },u_{0}\right) $ is independent of ${\hat{M}}$. If,
additionally, $x_{0}$ is normally distributed, then $x_{0}^{\prime }{[\hat{M}%
]}\beta _{\hat{M}}^{(\star )}$ is also the best predictor in the class of
all predictors depending only on $x_{0}^{\prime }{[\hat{M}]}$ and ${\hat{M}}$%
.
\end{rem}

After having motivated the design-independent target, we shall, in the
remainder of this section, treat $x_{0}$ as fixed (but see Remark \ref%
{rem:xzero:random} in Appendix \ref{rem_on_theorem_3.6} for the case where $%
x_{0}$ is random). We now proceed to show that the confidence intervals
constructed in Section \ref{section:design:dependent} are also valid as
confidence intervals for the design-independent target $x_{0}^{\prime }[\hat{%
M}]\beta _{\hat{M}}^{(\star )}$ in an asymptotic sense under some mild
conditions. While the results in Section \ref{section:design:dependent}
apply to \emph{any} model selection procedure whatsoever (in case that $%
\mathcal{M}$ is the power set of $\left\{ 1,\ldots ,p\right\} $ as is the
case in the present section), we need here to make the following mild
assumption on the model selection procedure.

\begin{cond}
\label{cond:model:selection:procedure} The model selection procedure
satisfies: For any $M\subseteq \left\{ 1,\ldots ,p\right\} $ with $|M|<p$
and for any $\delta >0$, 
\begin{equation*}
\sup \left\{ P_{n,\beta ,\sigma }(\hat{M}=M|X):\beta \in \mathbb{R}%
^{p},\sigma >0,\left\Vert \beta \lbrack M^{c}]\right\Vert /\sigma \geq
\delta \right\} \rightarrow 0
\end{equation*}%
in probability as $n\rightarrow \infty $.
\end{cond}

Condition \ref{cond:model:selection:procedure} is very mild and typically
holds for model selection procedures such as AIC- and BIC-based procedures
as well as Lasso-type procedures. [This can be established along the lines
of the proof of Corollary 5.4(a) in \cite{leeb03finite}.] In addition, we
assume the following condition on the behavior of the design matrix.

\begin{cond}
\label{cond:L} The sequence of random matrices $\sqrt{n}\left[ \left(
X^{\prime }X/n\right) -\Sigma \right] $ is bounded in probability.
\end{cond}

Condition \ref{cond:L} holds, for example, when the rows of $X$ are
independent, or weakly dependent, and when the distribution $\mathcal{L}$
has finite fourth moments for all its components. We also introduce the
following condition.

\begin{cond}
\label{cond:sigma} The degrees of freedom parameters $r$ of the sequence of
estimators $\hat{\sigma}^{2}$ satisfy $r\rightarrow \infty $ as $%
n\rightarrow \infty $.
\end{cond}

Of course, if we choose for $\hat{\sigma}^{2}$ the usual variance estimator $%
\hat{\sigma}_{OLS}^{2}$ then this condition is certainly satisfied with $%
r=n-p$. We are now in the position to present the asymptotic coverage
result. Recall that the confidence intervals corresponding to $K_{i}$ with $%
2\leq i\leq 5$ depend on $x_{0}$ only through $x_{0}{[\hat{M}]}$ (or not on $%
x_{0}$ at all).

\begin{theo}
\label{theo:asympt:val:CI} Suppose Conditions \ref%
{cond:model:selection:procedure} and \ref{cond:L} hold.

(a) Suppose also that Condition \ref{cond:sigma} is satisfied. Let $%
CI(x_{0}) $ be the confidence interval (\ref{eq:general:form:CI}) where the
constant $K(x_{0},\hat{M})$ is given by the constant $K_{1}(x_{0},r)$
defined in Section \ref{section:design:dependent}. Then the confidence
interval $CI(x_{0})$ satisfies 
\begin{equation}
\inf_{x_{0}\in \mathbb{R}^{p},\beta \in \mathbb{R}^{p},\sigma >0}P_{n,\beta
,\sigma }\left( \left. x_{0}^{\prime }[\hat{M}]\beta _{\hat{M}}^{(\star
)}\in CI(x_{0})\right\vert X\right) \geq (1-\alpha )+o_{p}(1),
\label{eq:asymptotically:valid:CI}
\end{equation}%
where the $o_{p}(1)$ term above depends only on $X$ and converges to zero in
probability as $n\rightarrow \infty $. Relation (\ref%
{eq:asymptotically:valid:CI}) a fortiori holds if the confidence interval $%
CI(x_{0})$ is based on the constants $K_{2}(x_{0}[\hat{M}],\hat{M},r)$, $%
K_{3}(x_{0}[\hat{M}],\hat{M},r)$, $K_{4}\left( r\right) $, or $K_{5}\left(
r\right) $, respectively.

(b) Let $\tilde{\sigma}$ be an arbitrary estimator satisfying 
\begin{equation}
\sup_{\beta \in \mathbb{R}^{p},\sigma >0}P_{n,\beta ,\sigma }(\left\vert 
\tilde{\sigma}/\sigma -1\right\vert \geq \delta \left\vert X\right. )\overset%
{p}{\rightarrow }0  \label{eq:unif_cons_var}
\end{equation}%
for any $\delta >0$ as $n\rightarrow \infty $. Let further $r^{\ast
}=r_{n}^{\ast }$ be an arbitrary sequence in $\mathbb{N\cup }\left\{ \infty
\right\} $ satisfying $r^{\ast }\rightarrow \infty $ for $n\rightarrow
\infty $. Let $CI^{\ast }(x_{0})$ denote the modified confidence interval
which is obtained by replacing $\hat{\sigma}$ by $\tilde{\sigma}$ and $%
K(x_{0},\hat{M})$ by $K_{1}(x_{0},r^{\ast })$ ($K_{2}(x_{0}[\hat{M}],\hat{M}%
,r^{\ast })$, $K_{3}(x_{0}[\hat{M}],\hat{M},r^{\ast })$, $K_{4}\left(
r^{\ast }\right) $, or $K_{5}\left( r^{\ast }\right) $, respectively) in (%
\ref{eq:general:form:CI}) (while keeping $\hat{M}$ unchanged). Then relation
(\ref{eq:asymptotically:valid:CI}) holds with $CI(x_{0})$ replaced by $%
CI^{\ast }(x_{0})$.
\end{theo}

Theorem \ref{theo:asympt:val:CI}(a) shows that for any $x_{0}\in \mathbb{R}%
^{p}$ the interval $CI(x_{0})$ is an asymptotically valid confidence
interval for the design-independent target and additionally that the lower
bound $(1-\alpha )+o_{p}(1)$ for the minimal (over $\beta $ and $\sigma $)
coverage probability can be chosen independently of $x_{0}$. Theorem \ref%
{theo:asympt:val:CI}(b) extends this result to a larger class of intervals.
[Note that Part (a) is in fact a special case of Part (b) obtained by
setting $\tilde{\sigma}=\hat{\sigma}$ and $r^{\ast }=r$ and observing that $%
\hat{\sigma}$ clearly satisfies the condition on $\tilde{\sigma}$ in Part
(b) under Condition \ref{cond:sigma}.] We note that applying Theorem \ref%
{theo:asympt:val:CI}(b) with $\tilde{\sigma}=\hat{\sigma}$ and $r^{\ast
}=\infty $ shows that Theorem \ref{theo:asympt:val:CI}(a) also continues to
hold for the confidence interval that is obtained by replacing the constants 
$K_{1}(x_{0},r)$ ($K_{2}(x_{0}[\hat{M}],\hat{M},r)$, $K_{3}(x_{0}[\hat{M}],%
\hat{M},r)$, $K_{4}\left( r\right) $, or $K_{5}\left( r\right) $,
respectively) by the constants $K_{1}(x_{0},\infty )$ ($K_{2}(x_{0}[\hat{M}],%
\hat{M},\infty )$, $K_{3}(x_{0}[\hat{M}],\hat{M},\infty )$, $K_{4}\left(
\infty \right) $, or $K_{5}\left( \infty \right) $, respectively).
Measurability issues regarding Theorem \ref{theo:asympt:val:CI} are
discussed in Appendix \ref{measurability}.

Condition (\ref{eq:unif_cons_var}) is a uniform consistency property. It is
clearly satisfied by $\hat{\sigma}_{OLS}^{2}$ (and more generally by the
estimator $\hat{\sigma}^{2}$ under Condition \ref{cond:sigma} as already
noted above), but it is also satisfied by the post-model-selection estimator 
$\hat{\sigma}_{\hat{M}}^{2}=||Y-X[\hat{M}]\hat{\beta}_{\hat{M}}||^{2}/(n-|%
\hat{M}|)$ provided the model selection procedure satisfies Condition \ref%
{cond:model:selection:procedure}, see Lemma \ref{lem: PMS-var} in Appendix %
\ref{app:B} for a precise result. As a consequence, Theorem \ref%
{theo:asympt:val:CI}(b) shows that the post-model-selection estimator $\hat{%
\sigma}_{\hat{M}}^{2}$ can be used instead of $\hat{\sigma}^{2}$ in the
construction of the confidence interval.

\begin{rem}
\label{rem:kmown:variance:case_2}\normalfont\emph{(Infeasible variance
estimators) }Theorem \ref{theo:asympt:val:CI}(a) remains valid if $\hat{%
\sigma}^{2}$ is allowed to depend also on $\sigma $ but otherwise satisfies
the assumptions made earlier or if $\hat{\sigma}^{2}=\sigma ^{2}$ and $%
r=\infty $. Similarly, Theorem \ref{theo:asympt:val:CI}(b) remains valid if $%
\tilde{\sigma}^{2}$ is allowed to be infeasible. Furthermore, a remark
similar to Remark \ref{rem_2.6}(iii) also applies here.
\end{rem}

\begin{rem}
\label{restr_univ}\normalfont\emph{(Restricted universe of selected models) }%
Theorem \ref{theo:asympt:val:CI} can easily be generalized to the case where
a universe $\mathcal{M}$ different from the power set of $\left\{ 1,\ldots
,p\right\} $ is employed, provided the full model $\left\{ 1,\ldots
,p\right\} $ belongs to $\mathcal{M}$ (and $\mathcal{M}$ satisfies the basic
assumptions made in Section \ref{section:design:dependent}).
\end{rem}

\section{Numerical study\label{section:simulation:study}}

We next present a numerical study of the lengths and the minimal coverage
probabilities of various confidence intervals. We begin, in Section \ref%
{length_numerical}, with an investigation of the length of the confidence
intervals introduced in Section \ref{section:design:dependent}, including
the `naive' confidence interval that ignores the model selection step, as a
function of the selected model. In Section \ref{minimal_cov_prob_numerical}
we then evaluate numerically the minimal coverage probabilities of these
confidence intervals. As model selectors we consider here AIC, BIC, LASSO,
SCAD (\cite{fan01variable}), and MCP (\cite{zhang10nearly}). Finally, in
Section \ref{taylor} we compare the confidence intervals introduced in
Section \ref{section:design:dependent} with the confidence interval proposed
recently in \cite{lee15exact}, which is specific to the LASSO model
selector. Code for the computations in this section is available from the
first author.

\subsection{Lengths of confidence intervals\label{length_numerical}}

We consider the lengths of the confidence intervals obtained from (\ref%
{eq:general:form:CI}) standardized by $\hat{\sigma}$, i.e., we consider $%
2K(x_{0},\hat{M})\Vert s_{\hat{M}}\Vert $ for the six cases where $K(x_{0},%
\hat{M})$ is replaced by either one of the five constants $K_{1}(x_{0})$, $%
K_{2}(x_{0}[\hat{M}],\hat{M})$, $K_{3}(x_{0}[\hat{M}],\hat{M})$, $K_{4}$, $%
K_{5}$ of Section \ref{section:design:dependent} or by the constant $%
K_{naive}=q_{r,1-\alpha /2}$, the $(1-\alpha /2)$-quantile of Student's
t-distribution with $r$ degrees of freedom. We recall that the constant $%
K_{naive}$ yields the `naive' confidence interval that ignores the model
selection step and that we have $K_{naive}\leq K_{1}(x_{0})\leq ...\leq
K_{5} $ (the first inequality holding provided $x_{0}\neq 0$).

For computing the standardized length, we set $\alpha =0.05$, $n=29$, $%
d=p=10 $, $r=n-p$, $\sigma =1$, and obtain $X$ and $x_{0}$ from a data set
of \cite{rawlings98applied} concerning the peak flow rate of watersheds.
This data set contains a $30\times 10$ design matrix $X_{Raw}$ corresponding
to ten explanatory variables. For a description of these variables see
Appendix \ref{app D}. This data set is also studied in \cite{kabaila06large}
and \cite{leeb13various}. We refer to it as the watershed data set, and $%
x_{0}$ and $X $ are chosen such that $(x_{0},X^{\prime })^{\prime }$ is
equal to the watershed design matrix $X_{Raw}$. It is easily checked that
the so-obtained matrix $X$ is indeed of full column rank (and $x_{0}\neq 0$%
). Furthermore, the model universe $\mathcal{M}$ is chosen to be the power
set of $\{1,...,p\}$.

For the so chosen values of $\alpha $, $n$, $p$, $r$, $\sigma $, $X$, $x_{0}$%
, and $\mathcal{M}$, we compute the standardized lengths $2K(x_{0},M)\Vert
s_{M}\Vert $ of the confidence intervals obtained by replacing $K(x_{0},M)$
by $K_{naive}$, $K_{1}(x_{0})$, $K_{2}(x_{0}[M],M)$, $K_{3}(x_{0}[M],M)$, $%
K_{4}$, and $K_{5}$, respectively. To ease the computational burden and to
enable a simple presentation as in Figure 1 below, we compute the
standardized lengths of the confidence intervals only for $M$ belonging to
the family $\{\{1\},...,\{1,...,10\}\}$ consisting of ten nested submodels.
[This does \emph{not} mean that we compute the constants $K_{i}$ under the
assumption of a restricted universe of models; recall that we use $\mathcal{M%
}$ equal to the power set of $\{1,...,p\}$.] The computation of $K_{naive}$, 
$K_{1}(x_{0})$, $K_{3}(x_{0}[M],M)$, $K_{4}$, and $K_{5}$ is either
straightforward or is obtained from the algorithms described in Appendix \ref%
{section:practical:algorithms}. However, computing $K_{2}(x_{0}[M],M)$ for $%
M\neq \{1,...,10\}$ necessitates to compute $\sup \{K_{1}(x):x[M]=x_{0}[M]\}$%
. We approximate this supremum by using a three-step Monte Carlo procedure
described in Appendix \ref{app D}.

\begin{center}
\begin{tabular}{c}
\includegraphics[scale=.7]{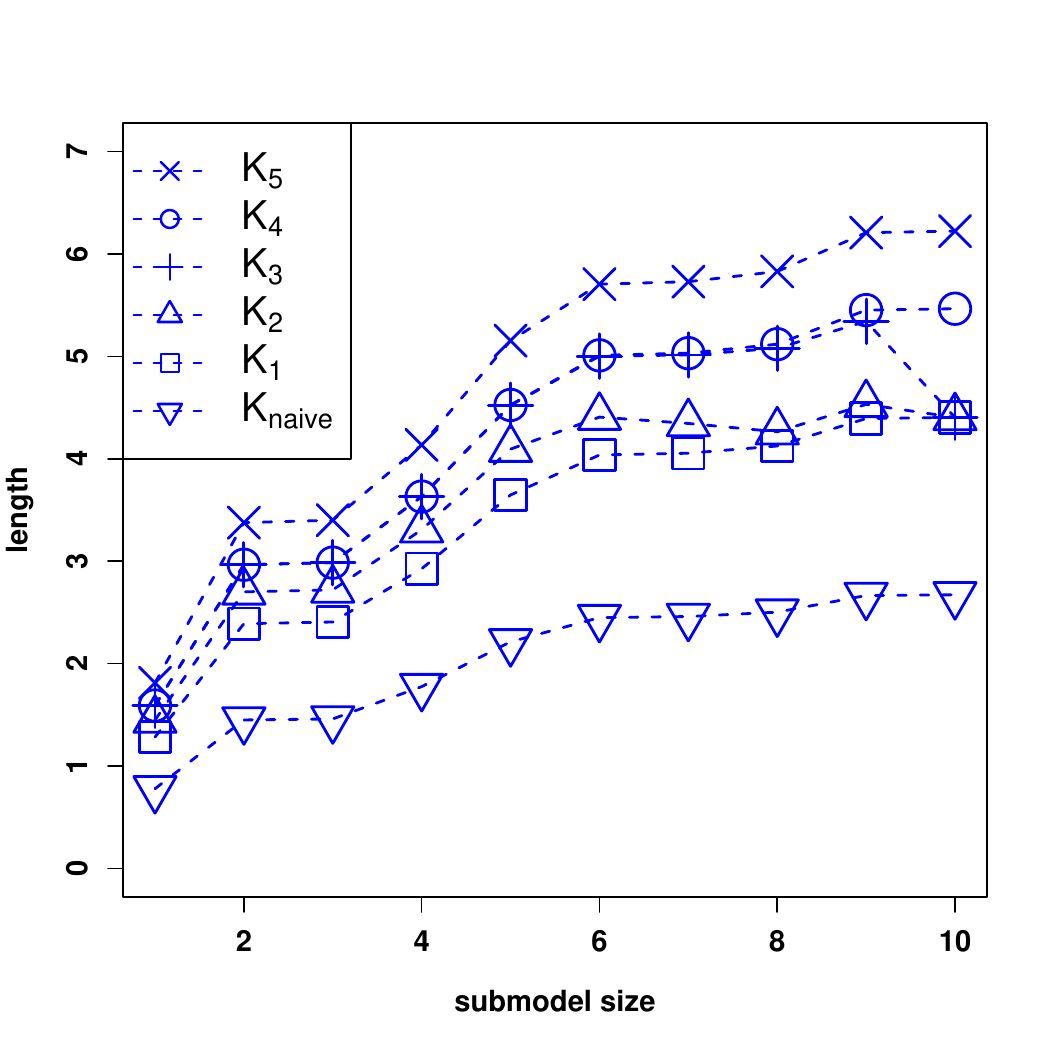}%
\end{tabular}
\begin{quote}
\label{figure1} Figure 1: Standardized lengths of various confidence
intervals as function of model size. Dashed lines are added to improve
readability.
\end{quote}
\end{center}

The standardized lengths of the confidence intervals corresponding to the
constants $K_{naive}$, $K_{1}$,\ldots , $K_{5}$ are reported in Figure 1 for
the ten nested submodels mentioned before. We first see that, for each of
the constants $K_{naive}$, $K_{1}$, $K_{4}$, and $K_{5}$, the standardized
length of the confidence interval increases with submodel size, which must
hold since these constants do not depend on the submodel $M$ and since the
term $||s_{M}||$ increases with submodel size (for nested submodels as
considered in Figure 1). However, as discussed after Proposition \ref%
{prop:PoSI:upper:bound}, the values of $K_{2}$ and $K_{3}$ decrease with
increasing submodel size for nested submodels. Figure 1 shows that the
combined effect of the increase of $||s_{M}||$ and the decrease of $K_{2}$
and $K_{3}$ with submodel size can be an increase or a decrease of the
standardized lengths of the confidence intervals. Indeed, the standardized
lengths increase globally (i.e., from submodel size $1$ to $10$), but can
decrease locally (for example, the standardized length of the confidence
interval obtained from $K_{2}$ decreases from submodel size $6$ to submodel
size $8$; for the interval obtained from $K_{3}$ the standardized length
decreases from submodel size $9$ to submodel size $10$). In Figure 1 the
decreases of the standardized lengths occur only between submodel sizes for
which $||s_{M}||$ is almost constant with $M$ (which can be seen from the
standardized lengths obtained from, say, $K_{5}$, since they are
proportional to $||s_{M}||$). We also see from Figure 1 that the `naive'
interval is much shorter than the other intervals (at the price of typically
not having the correct minimal coverage probability). The difference in
standardized length between the intervals based on $K_{1}$ and $K_{2}$,
respectively, is noticeable but not dramatic. A larger increase in
standardized length is noted when comparing the interval based on the
costly-to-compute constant $K_{2}$ with the one obtained from $K_{3}$,
especially for submodel sizes $6$ to $9$. Furthermore, the standardized
lengths of the confidence intervals obtained from $K_{3}$ are very close to
those obtained from $K_{4}$ for model size $1$ to $9$; cf. (\ref{Closeness
of K3 and K4}). Finally, in Figure 1 we also see that the confidence
intervals obtained from $K_{1}$, $K_{2}$, and $K_{3}$ have the same
standardized length when the model size is $10$, and that the same is true
for the confidence intervals obtained from $K_{3}$ and $K_{4}$ when the
model size is $1$. This, of course, is not a coincidence, but holds
necessarily as has been noted in the discussion of Proposition \ref%
{prop:PoSI:upper:bound}.

Additional computations of confidence interval lengths, with $X$ and $x_{0}$
now randomly generated, yield results very similar to those in Figure 1. For
the sake of brevity, these results are not shown here. We find, in
particular, that the standardized length of the confidence interval obtained
from $K_{3}$ always increases with submodel size when they are averaged with
respect to $X$ and $x_{0}$, but, as in Figure 1, can decrease locally when
not averaged. [In these additional numerical studies we did not consider the
constant $K_{2}$ due to the high computational cost involved in its
evaluation.]

\subsection{Minimal coverage probabilities\label{minimal_cov_prob_numerical}}

In this section we consider the case where $\mu =X\beta $ and $d=p<n$, i.e.,
the case where the given matrix $X$ has full rank less than $n$ and provides
a correct linear model for the data $Y$. We then investigate the minimal
coverage probabilities (the minimum being w.r.t. $\beta \in \mathbb{R}^{p}$
and $\sigma \in (0,\infty )$) of the intervals obtained from the constants $%
K_{naive}$, $K_{1}$, $K_{3}$, and $K_{4}$ when used as confidence intervals
for the target $x_{0}^{\prime }[\hat{M}]\beta _{\hat{M}}^{(n)}$ on the one
hand as well as for the target $x_{0}^{\prime }[\hat{M}]\beta _{\hat{M}%
}^{(\star )}$ on the other hand. The constants $K_{1}$, $K_{3}$, and $K_{4}$
are computed based on $\mathcal{M}$ equal to the power set of $\{1,...,p\}$.
We do not report results for confidence intervals obtained from $K_{2}$,
since the computation of $K_{2}$ is too costly for the study we present
below. The results for confidence intervals obtained from $K_{5}$ would be
qualitatively similar to those for confidence intervals obtained from $K_{4}$%
, so we do not report them for the sake of brevity.

We consider minimal coverage probabilities in the setting where $\alpha
=0.05 $, $p=10$, $n=20$ or $n=100$, and the variance parameter is estimated
by the standard unbiased estimator using the full model, so that $r=n-p$.
For model selection we consider AIC-, BIC-procedures, the LASSO, SCAD (\cite%
{fan01variable}), and MCP (\cite{zhang10nearly}). Tuning parameters of the
latter three procedures are chosen by cross-validation. For all five
procedures we always protect the first explanatory variable (which
corresponds to an intercept term) from selection. However, note that the
information that the first variable is protected is \emph{not} used in
computing the constants $K_{i}$, i.e., we do not use a restricted universe
of models but use $\mathcal{M}$ equal to the power set of $\{1,...,p\}$.
[Additional simulations with no intercept term and no protected explanatory
variable lead to results very similar to the ones given in Table \ref%
{table:coverage} below.] Computational details regarding these procedures
can be found in Appendix \ref{app D}.

The design matrix $X$ and the vector $x_{0}$ are generated in the following
manner: The $10\times 10$ matrix $\Sigma $ of (uncentered) second moments is
chosen to be of the form%
\begin{equation*}
\Sigma =\left( 
\begin{array}{cc}
1 & 0\cdots 0 \\ 
0 & \raisebox{-15pt}{{\huge\mbox{{$\tilde{\Sigma}$}}}} \\[-4ex] 
\vdots &  \\[-0.5ex] 
0 & 
\end{array}%
\right) ,
\end{equation*}%
where we consider three choices for the $9\times 9$ matrix $\tilde{\Sigma}$.
For the first case, $\tilde{\Sigma}$ is obtained by removing the first row
and column of the $10\times 10$ empirical covariance matrix (standardized by 
$30-1=29$) of the variables in the $30\times 10$ watershed design matrix $%
X_{Raw}$. For the second case, we set $\tilde{\Sigma}=I_{\tilde{p}}+(2a+%
\tilde{p}a^{2})E_{\tilde{p}}$ with $\tilde{p}=9$, $a=10$, and with $E_{%
\tilde{p}}$ the $\tilde{p}\times \tilde{p}$ matrix which has all entries
equal to $1$. For the third case $\tilde{\Sigma}$ coincides with the
identity matrix $I_{\tilde{p}}$, except that the zero elements in the last
row and column of $I_{\tilde{p}}$ are replaced by the constant $c=\sqrt{0.8/(%
\tilde{p}-1)}$ where $\tilde{p}=9$. Similar as in \cite{berk13valid} and 
\cite{leeb13various}, we refer to the data set obtained in the second case
as the exchangeable data set (as the covariance matrix $\tilde{\Sigma}$ is
permutation-invariant), and to the one obtained in the third case as the
equicorrelated data set (as $\tilde{\Sigma}$ is the correlation matrix of a
random vector, the last component of which has the same correlation with all
the other components); see Appendix \ref{app D} for more details. For a
given configuration of $n$ and $\Sigma $, we then sample independently $n+1$
vectors of dimension $10\times 1$ such that for each of these vectors the
first component is $1$ and the remaining nine components are jointly
normally distributed with mean zero and covariance matrix $\tilde{\Sigma}$.
The transposes of the first $n$ of theses vectors now form the rows of the $%
n\times p$ design matrix $X$, while the $(n+1)$-th of these vectors is used
for the $p$-dimensional vector $x_{0}$. [It is easy to see that the
mechanism just described generates matrices of full column rank almost
surely. The matrices $X$ actually generated were additionally checked to be
of full column rank.]

Consider now a given configuration of $n$, $\Sigma $, the model selection
procedure, the target (either the design-dependent or the design-independent
target), as well as of a matrix $X$ and a vector $x_{0}$ that have been
obtained in the manner just described. Then we estimate the minimal (over $%
\beta $ and $\sigma $) coverage probabilities (conditional on $X$ and $x_{0}$%
) of the confidence intervals obtained from the constants $K_{naive}$, $%
K_{1} $, $K_{3}$, and $K_{4}$ for the given target under investigation. The
minimal coverage probabilities are estimated by a three-step Monte Carlo
procedure similar to that of \cite{leeb13various}, which is described in
detail in Appendix \ref{app D}. We stress here that the minimal coverage
probabilities found by this Monte Carlo procedure are simulation-based
results obtained by a stochastic search over a $10$-dimensional parameter
space, and thus only provide approximate upper bounds for the true minimal
coverage probabilities.

\begin{table}[tbp]
\begin{center}
\begin{tabular}{c|c|l|cccc|cccc}
Data set & $n$ & Model & \multicolumn{8}{c}{Target} \\ \cline{4-11}
&  & selector & \multicolumn{4}{c|}{design-dependent} & \multicolumn{4}{|c}{
design-independent} \\ 
&  &  & \multicolumn{4}{c|}{$x_{0}[\hat{M}]^{\prime }\beta _{\hat{M}}^{(n)}$}
& \multicolumn{4}{|c}{$x_{0}[\hat{M}]^{\prime }\beta _{\hat{M}}^{(\star )}$}
\\ 
&  &  & $K_{naive}$ & $K_{1}$ & $K_{3}$ & $K_{4}$ & $K_{naive}$ & $K_{1}$ & $%
K_{3}$ & $K_{4}$ \\ \hline
\multirow{4}{*}{Watershed} & $20$ & AIC & 0.84 & 0.99 & 1.00 & 1.00 & 0.79 & 
0.97 & 0.99 & 0.99 \\ 
& $20$ & BIC & 0.84 & 0.99 & 1.00 & 1.00 & 0.74 & 0.96 & 0.98 & 0.98 \\ 
& $20$ & LASSO & 0.90 & 1.00 & 1.00 & 1.00 & 0.18 & 0.48 & 0.61 & 0.61 \\ 
& $20$ & SCAD & 0.90 & 0.99 & 1.00 & 1.00 & 0.45 & 0.77 & 0.84 & 0.84 \\ 
& $20$ & MCP & 0.89 & 0.99 & 1.00 & 1.00 & 0.47 & 0.78 & 0.85 & 0.85 \\ 
\cline{2-11}
& $100$ & AIC & 0.87 & 0.99 & 1.00 & 1.00 & 0.88 & 0.99 & 1.00 & 1.00 \\ 
& $100$ & BIC & 0.88 & 0.99 & 1.00 & 1.00 & 0.87 & 0.99 & 1.00 & 1.00 \\ 
& $100$ & LASSO & 0.88 & 0.99 & 1.00 & 1.00 & 0.87 & 0.99 & 1.00 & 1.00 \\ 
& $100$ & SCAD & 0.88 & 0.99 & 1.00 & 1.00 & 0.88 & 0.99 & 1.00 & 1.00 \\ 
& $100$ & MCP & 0.88 & 0.99 & 1.00 & 1.00 & 0.88 & 0.99 & 1.00 & 1.00 \\ 
\hline
\multirow{4}{*}{Exchangeable} & $20$ & AIC & 0.83 & 0.99 & 1.00 & 1.00 & 0.80
& 0.98 & 0.99 & 0.99 \\ 
& $20$ & BIC & 0.84 & 0.99 & 1.00 & 1.00 & 0.76 & 0.97 & 0.99 & 0.99 \\ 
& $20$ & LASSO & 0.90 & 1.00 & 1.00 & 1.00 & 0.46 & 0.86 & 0.93 & 0.92 \\ 
& $20$ & SCAD & 0.91 & 1.00 & 1.00 & 1.00 & 0.55 & 0.90 & 0.94 & 0.94 \\ 
& $20$ & MCP & 0.91 & 1.00 & 1.00 & 1.00 & 0.54 & 0.89 & 0.94 & 0.94 \\ 
\cline{2-11}
& $100$ & AIC & 0.89 & 0.99 & 1.00 & 1.00 & 0.90 & 0.99 & 1.00 & 1.00 \\ 
& $100$ & BIC & 0.90 & 0.99 & 1.00 & 1.00 & 0.90 & 0.99 & 1.00 & 1.00 \\ 
& $100$ & LASSO & 0.90 & 0.99 & 1.00 & 1.00 & 0.90 & 0.99 & 1.00 & 1.00 \\ 
& $100$ & SCAD & 0.90 & 0.99 & 1.00 & 1.00 & 0.90 & 0.99 & 1.00 & 1.00 \\ 
& $100$ & MCP & 0.90 & 0.99 & 1.00 & 1.00 & 0.90 & 0.99 & 1.00 & 1.00 \\ 
\hline
\multirow{4}{*}{Equicorrelated} & $20$ & AIC & 0.83 & 0.99 & 1.00 & 1.00 & 
0.79 & 0.98 & 0.99 & 0.99 \\ 
& $20$ & BIC & 0.81 & 0.99 & 1.00 & 1.00 & 0.74 & 0.98 & 0.99 & 0.99 \\ 
& $20$ & LASSO & 0.88 & 1.00 & 1.00 & 1.00 & 0.39 & 0.71 & 0.79 & 0.79 \\ 
& $20$ & SCAD & 0.88 & 0.99 & 1.00 & 1.00 & 0.67 & 0.92 & 0.95 & 0.96 \\ 
& $20$ & MCP & 0.86 & 0.99 & 1.00 & 1.00 & 0.66 & 0.93 & 0.96 & 0.96 \\ 
\cline{2-11}
& $100$ & AIC & 0.84 & 0.99 & 1.00 & 1.00 & 0.84 & 0.99 & 1.00 & 1.00 \\ 
& $100$ & BIC & 0.86 & 0.99 & 1.00 & 1.00 & 0.86 & 0.99 & 1.00 & 1.00 \\ 
& $100$ & LASSO & 0.88 & 1.00 & 1.00 & 1.00 & 0.88 & 1.00 & 1.00 & 1.00 \\ 
& $100$ & SCAD & 0.88 & 0.99 & 1.00 & 1.00 & 0.89 & 1.00 & 1.00 & 1.00 \\ 
& $100$ & MCP & 0.88 & 0.99 & 1.00 & 1.00 & 0.89 & 0.99 & 1.00 & 1.00%
\end{tabular}%
\end{center}
\caption{Monte Carlo estimates of the minimal coverage probabilities (w.r.t. 
$\protect\beta $ and $\protect\sigma $) of various confidence intervals. The
nominal coverage probability is $1-\protect\alpha =0.95$ and $p=10$.}
\label{table:coverage}
\end{table}

Table \ref{table:coverage} summarizes the estimated minimal coverage
probabilities for the various confidence sets and targets, and for the
model-selection procedures and data sets considered in the study. The
conclusions are pretty much the same for the three data sets. First, we
observe that, for $n=20$, the differences of minimal coverage probabilities
between the design-dependent and independent targets can be significant,
especially for the `naive' intervals and for the other intervals in case the
LASSO, SCAD, or MCP model selectors are used. However, for $n=100$, these
differences are very small for all the configurations. This is in line with
Lemma \ref{lem:equivalence:beta:betan} in Appendix \ref{app:B}, which
entails that for a large family of model selection procedures, the
difference of coverage probabilities between the two targets vanishes,
uniformly in $\beta $ and $\sigma $, when $n$ increases. For $n=100$, the
results are thus almost identical for the two targets: For the five model
selection procedures, the confidence intervals obtained from the constants $%
K_{1}$, $K_{3}$, and $K_{4}$ are valid, while the `naive' confidence
intervals are moderately too short, so that their minimal coverage
probabilities are below the nominal level, with a minimum of $0.84$.

For $n=20$ and when AIC or BIC is used, the `naive' confidence intervals
fail to have the right coverage probabilities to a somewhat larger extent
than in case $n=100$. Their minimal coverage probabilities can be as small
as $0.81$ for the design-dependent target and $0.74$ for the
design-independent target. [Note that, for the design-dependent target, for $%
n=20$ and $n=100$, the coverage probabilities of the `naive' confidence
interval are generally smaller for the equicorrelated data set than for the
exchangeable data set. This can possibly be explained by the fact that
Theorems 6.1 and 6.2 in \cite{berk13valid} suggest that $K_{1}$ should be
larger for the equicorrelated data set than for the exchangeable data set.
Hence, for the equicorrelated data set, larger confidence intervals seem to
be needed to have the required minimal coverage probability for all model
selection procedures.] Furthermore, again for $n=20$ and when AIC or BIC is
used, the confidence intervals obtained from the constants $K_{1}$, $K_{3}$,
and $K_{4}$ remain valid here for both targets.

However, when $n=20$ and the LASSO model selector is used, the results for
the design-independent target are drastically different from those obtained
with the AIC- or BIC-procedures: All confidence intervals have minimal
coverage probabilities for the design-independent target that are below, and
in most cases significantly below, the nominal level. The failure of all the
confidence intervals is here often more pronounced than the failure of the
`naive' confidence intervals when other model selectors are used. Especially
for the watershed data set, the estimated minimal coverage probability is $%
0.18$ for the `naive' interval and $0.48$ for the confidence interval based
on $K_{1}$. The reason for this phenomenon can be traced to the observation
that the LASSO model selector, as implemented here and for the parameters
used in the stochastic search for the smallest coverage probability, selects
models that are significantly smaller than those AIC and BIC select. In
particular, the LASSO procedure often excludes regressors for which the
corresponding regression coefficients are not small. In our simulation
study, selecting a small model, that excludes regressors with significant
coefficients, makes the difference between the design-dependent and
design-independent targets larger. Since the confidence intervals are
designed to cover the former target, they hence have a hard time to cover
the latter when the two targets are significantly different. In other words,
for $n=20$ the supremum in the display in Condition \ref%
{cond:model:selection:procedure} is not small for the LASSO procedure, so
that the asymptotics in Theorem \ref{theo:asympt:val:CI} does not provide a
good approximation for the finite-sample situation. Finally, for $n=20$ and
for the design-independent target, the results for the SCAD and MCP model
selectors lie somewhere in between those of the AIC and BIC and those of the
LASSO model selectors. Indeed, for SCAD and MCP, the confidence intervals
often fail to have the required minimal coverage probabilities, but less
severely than for the LASSO. We stress that the preceding conclusions hold
for the LASSO, SCAD, and MCP procedures as implemented here where tuning
parameters are chosen by cross-validation. Other implementations of these
procedures may of course give different results.

The results in Table \ref{table:coverage} concern the coverage probabilities
conditional on the design matrix $X$ and on $x_{0}$, and thus depend on the
values of $X$ and $x_{0}$ used. In additional (non-exhaustive) simulations
we have repeated the above analysis for other values of $X$ and $x_{0}$ and
have found similar results.

\subsection{Comparison with the confidence interval of \protect\cite%
{lee15exact} \label{taylor}}

In this section we now compare the confidence intervals of Section \ref%
{section:design:dependent} with a confidence interval recently introduced in 
\cite{lee15exact}. Again, we consider the case where $\mu =X\beta $ and $%
d=p<n$, and we focus on the design-dependent target $x_{0}^{\prime }[\hat{M}%
]\beta _{\hat{M}}^{(n)}$. As in \cite{lee15exact} we consider the
known-variance case and set $\sigma =1$ in this section. The confidence
interval of \cite{lee15exact} is dedicated to the LASSO model selector and
is given in the \verb|R| package accompanying that paper for the case where $%
x_{0}$ is a standard basis vector. We hence assume in the following that $%
x_{0}$ is equal to the first standard basis vector $e_{1}$. The proposed
interval is then \emph{conditionally} valid for the design-dependent target
in the following sense: Consider the model selector $\hat{M}$ obtained by
selecting those explanatory variables for which the LASSO estimator has
non-zero coefficients, with the penalty parameter $\lambda $ in (4.1) of 
\cite{lee15exact} being fixed, independently of $Y$. Then the interval
proposed by \cite{lee15exact}, which we denote by $\bar{CI}$, satisfies, for
any fixed $X$, for $x_{0}=e_{1}$, and for any fixed $M\subseteq \{1,...,p\}$
with $1\in M$, 
\begin{equation}
\inf_{\beta \in \mathbb{R}^{p}}P_{n,\beta ,1}\left( \left. x_{0}^{\prime }[%
\hat{M}]\beta _{\hat{M}}^{(n)}\in \bar{CI}\right\vert \hat{M}=M\right)
=1-\alpha ,  \label{eq:conditionally:valid}
\end{equation}%
with the convention that the probability in the above display is $1$ if $%
P_{n,\beta ,1}(\hat{M}=M)=0$. The computation of $\bar{CI}$ for a given
value of $\hat{M}$ can be carried out without observing $x_{0}[\hat{M}^{c}]$%
, which is also the case for the confidence intervals obtained from $%
K_{2},...,K_{5}$, but not for that obtained from $K_{1}$. Furthermore, the
computation of $\bar{CI}$ (when the conditioning additionally is also on the
signs, see \cite{lee15exact} for details) entails a cost that grows linearly
with $p$. Thus, $\bar{CI}$ can be implemented for significantly larger
values of $p$ than the confidence intervals based on $K_{1},...,K_{3}$
currently can be. We note for later use that, in the case $x_{0}=e_{1}$
considered here, the interval $\bar{CI}$ as given in \cite{lee15exact} is
not defined on the event that a model $\hat{M}$ is selected that does not
contain $1$. Hence, we can not speak about unconditional coverage without
amending the definition in \cite{lee15exact}. [A possible amendment,
consistent with our conventions and maximizing unconditional coverage among
all possible amendments, is to recall that $x_{0}^{\prime }[\hat{M}]\beta _{%
\hat{M}}^{(n)}=0$ if $1\notin \hat{M}$ and to set $\bar{CI}=\{0\}$ on this
event. With such an amendment, $\bar{CI}$ then a fortiori has minimal
unconditional coverage probability not less than $1-\alpha $.]

Despite being specific to the LASSO model selector with fixed $\lambda $, we
nevertheless find below that the confidence interval of \cite{lee15exact} is 
\emph{not} shorter than those based on $K_{1}$, $K_{3}$, and $K_{4}$
(presumably due to the fact that (\ref{eq:conditionally:valid}) imposes a
stricter requirement than requiring only correct unconditional coverage). In
addition, we point out that this confidence interval can be very sensitive
to deviations from the specific model selector it is designed for: In
particular, we show that its coverage can break down, when the LASSO model
selector is used but with a \emph{data-dependent} penalty parameter $\lambda 
$ selected by cross-validation. This is certainly in stark contrast to the
confidence intervals obtained from $K_{1}$, $K_{3}$, and $K_{4}$, that are
valid for \emph{any} model selection procedure whatsoever.

We first present the results on confidence interval lengths. We conduct a
Monte-Carlo study in the case $p=10$, $n=100$, and $\alpha =0.05$. We use
the function \verb|fixedLassoInf| of the \verb|R| package 
\verb|selectiveInference| to construct the confidence interval of \cite%
{lee15exact}. In line with the presentation in that paper, this function is
designed for the cases where $x_{0}$ is one of the standard basis vectors of 
$\mathbb{R}^{p}$ and thus we set $x_{0}^{\prime }=(1,0,...,0)$ as already
mentioned above. We consider two different settings for $X$ and $\beta $. In
the `independent' setting, we sample independently $1000$ values of $X$, $%
\beta $, and $Y$ in the following way: We first sample the (transposes of
the) rows of $X$ and the vector $\beta $ as $n+1$ independent draws from the
standard Gaussian distribution on $\mathbb{R}^{p}$. Given $X$ and $\beta $,
we then sample $Y$ from the $N(X\beta ,I_{n})$-distribution. For each of the 
$1000$ values of $X$, $\beta $, and $Y$ so obtained, we run the LASSO model
selector, with $\lambda $ fixed as a function of $X$ as described at the
beginning of Section 7 of \cite{lee15exact}. We use the function 
\verb|glmnet| of the \verb|R| package \verb|glmnet| to compute the LASSO
model selector. Then, if the first variable is included in the selected
model, we record the lengths of the confidence interval $\bar{CI}$ and of
the confidence intervals (for the same target) obtained from $K_{1}$, $K_{3}$%
, and $K_{4}$ (where these three constants are computed with $r=\infty $ and 
$\mathcal{M}$ equal to the power set of $\left\{ 1,\ldots ,p\right\} $). If
not, we discard the realization of $X$, $\beta $, and $Y$ (this is in line
with the fact that the $\bar{CI}$ is only defined on the event where $1\in 
\hat{M}$ as discussed above). In the `correlated' setting we proceed as just
described, with the only difference that the rows of $X$ are sampled
according to the Gaussian distribution with mean vector $0\in \mathbb{R}^{p}$
and covariance matrix $[\exp (-|i-j|/10)]_{i,j=1}^{10}$ and that then $\beta 
$ is sampled from a random vector $b$, so that $Xb$ follows the standard
Gaussian distribution within the column space of $X$. [We note that this
mechanism almost surely generates matrices $X$ that have full column rank.]

The medians and empirical $90\%$-quantiles of the confidence interval
lengths' distributions obtained that way are reported in Table \ref%
{table:comparison:lengths}. The conclusion is that there is no unilateral
hierarchy of the two methods for confidence interval construction (that of 
\cite{lee15exact} and that using the intervals based on $K_{1}$, $K_{3}$,
and $K_{4}$) in terms of median length. Depending on the situation, any of
the two methods can provide the smallest median length. The $90\%$%
-quantiles, on the other hand, are always larger for the confidence interval 
$\bar{CI}$ of \cite{lee15exact} than for those obtained from $K_{1}$, $K_{3}$%
, and $K_{4}$. [The feature that $\bar{CI}$ can be very long (with small but
non-negligible probability) has also been noted in \cite{lee15exact}, and is
not shared by the confidence intervals obtained from $K_{1}$, $K_{3}$, and $%
K_{4}$.] Note finally that we have obtained the same conclusions in other
length simulations, which we do not report for the sake of brevity.

\begin{table}[tbp]
\begin{center}
\begin{tabular}{c|c|cccc}
Setting & Lengths & \multicolumn{4}{c}{Confidence interval} \\ 
&  & $K_{1}$ & $K_{3}$ & $K_{4}$ & $\bar{CI}$ \\ \hline
\multirow{2}{*}{`Independent'} & Median & 0.46 & 0.78 & 0.78 & 0.43 \\ 
& $90\%$-quantile & 0.51 & 0.85 & 0.85 & 1.06 \\ \hline
\multirow{2}{*}{`Correlated'} & Median & 0.56 & 0.81 & 0.81 & 1.42 \\ 
& $90\%$-quantile & 0.90 & 1.30 & 1.30 & 14.3%
\end{tabular}%
\end{center}
\caption{Medians and empirical quantiles of the lengths of the confidence
intervals $\bar{CI}$ of \protect\cite{lee15exact} and of those obtained from 
$K_{1}$, $K_{3}$, and $K_{4}$. The nominal coverage probability is $1-%
\protect\alpha =0.95$, $n=100$, and $p=10$.}
\label{table:comparison:lengths}
\end{table}

We now demonstrate that the confidence interval $\bar{CI}$ of \cite%
{lee15exact} can have conditional coverage probability considerably smaller
than the nominal one when $\lambda $ is selected by cross-validation. Rather
than evaluating the minimal conditional coverage probabilities conditional
on $\hat{M}=M$ separately for every $M$ satisfying $1\in M$, which would be
quite costly, we evaluate the minimal conditional coverage probability where
conditioning is on the event that $1\in \hat{M}$. We denote this quantity by 
$P_{cond,\min }$. A simple calculation shows that if we find that this
latter minimal conditional coverage probability is smaller than $1-\alpha $,
then it follows that we must have 
\begin{equation*}
\inf_{\beta \in \mathbb{R}^{p}}P_{n,\beta ,1}\left( \left. x_{0}^{\prime }[%
\hat{M}]\beta _{\hat{M}}^{(n)}\in \bar{CI}\right\vert \hat{M}=M\right)
<1-\alpha
\end{equation*}%
for at least some $M$ satisfying $1\in M$, showing that property (\ref%
{eq:conditionally:valid}) is violated. [To see this, note that $P_{n,\beta
,1}(\left. \cdot \right\vert 1\in \hat{M})$ is a convex combination (over
all $M$ with $1\in M$) of the probabilities $P_{n,\beta ,1}(\left. \cdot
\right\vert \hat{M}=M)$, with the (nonnegative) weights summing to $1$.]

In order to numerically evaluate the minimal conditional coverage
probability $P_{cond,\min }$ we proceed as follows: We consider eight
configurations given by all the possible combinations of $n=20,100$, $p=2,10$
and $\alpha =0.05,0.2$. Recall that $x_{0}^{\prime }=(1,0,...,0)$. For each
of these eight configurations, the (transpose of the) rows of $X$ are
sampled once from the $N(0,\Sigma )$-distribution and then remain fixed
throughout the minimal coverage probability evaluation. For $p=2$, we take $%
\Sigma $ to have $1$ as the diagonal and $0.8$ as the off-diagonal elements.
For $p=10$, we take $\Sigma $ to coincide with the identity matrix $I_{p}$,
except that the zero elements in the last row and column of $I_{p}$ are
replaced by the constant $c=\sqrt{0.8/(p-1)}$. For each of these eight
configurations, we carry out a three-step minimal conditional coverage
probability evaluation as described in Appendix \ref{app D}.

The so evaluated minimal conditional coverage probabilities $\hat{P}%
_{cond,\min }$, say, are presented in Table \ref{table:minimal:cond:coverage}%
. For comparison, we also provide similar evaluations of minimal conditional
coverage probabilities using the same procedure as described in Appendix \ref%
{app D}, but now with $\lambda $ fixed as a function of $X$ as in the
beginning of Section 7 of \cite{lee15exact}. When $\lambda $ is fixed, these
minimal conditional coverage probabilities are approximately equal to the
nominal level $1-\alpha $, in agreement with the results of \cite{lee15exact}%
. However, when $\lambda $ is selected by cross-validation, the evaluated
minimal conditional coverage probabilities can be way below the nominal
level. In particular, these probabilities can be equal to $0.31$ for a
nominal level of $0.80$ and to $0.86$ for a nominal level of $0.95$.

In addition, for $\lambda $ selected by cross-validation and in all the
configurations of $n$, $p$, and $\alpha $, for the vector $\beta $ leading
to the minimal conditional coverage probability $\hat{P}_{cond,\min }$, we
can also estimate the unconditional coverage probability $P_{n,\beta
,1}(x_{0}^{\prime }[\hat{M}]\beta _{\hat{M}}^{(n)}\in \bar{CI})$ by $\hat{P}%
(1\in \hat{M})\hat{P}_{cond,\min }+\hat{P}(1\not\in \hat{M})$. [Here we make
use of the aforementioned amendment to $\bar{CI}$ in order to allow for a
well-defined unconditional coverage probability.] In this estimate, $\hat{P}%
(1\in \hat{M})$ is the proportion of times the first regressor belongs to $%
\hat{M}$, over the Monte-Carlo samples in the third step (as described in
Appendix \ref{app D}). The so evaluated unconditional coverage probabilities
are $0.60$ for $p=2$, $n=100$, $1-\alpha =0.80$ and $0.91$ for $p=2$, $n=100$%
, $1-\alpha =0.95$, which implies that the confidence intervals of \cite%
{lee15exact} also have minimal \emph{unconditional} coverage probabilities 
\emph{below} the nominal level when $\lambda $ is estimated by
cross-validation.

\begin{table}[tbp]
\begin{center}
\begin{tabular}{c|c|c|cc}
$p$ & $n$ & $1-\alpha $ & \multicolumn{2}{c}{LASSO} \\ 
&  &  & Fixed $\lambda $ & CV-selected $\lambda $ \\ \hline
\multirow{4}{*}{2} & \multirow{2}{*}{$20$} & $0.80$ & $0.80$ & $0.43$ \\ 
&  & $0.95$ & $0.95$ & $0.93$ \\ \cline{2-5}
& \multirow{2}{*}{$100$} & $0.80$ & $0.80$ & $0.31$ \\ 
&  & $0.95$ & $0.95$ & $0.86$ \\ \hline
\multirow{4}{*}{10} & \multirow{2}{*}{$20$} & $0.80$ & $0.79$ & $0.79$ \\ 
&  & $0.95$ & $0.94$ & 0.$93$ \\ \cline{2-5}
& \multirow{2}{*}{$100$} & $0.80$ & $0.79$ & $0.70$ \\ 
&  & $0.95$ & $0.95$ & $0.92$%
\end{tabular}%
\end{center}
\caption{Monte-Carlo estimates of the minimal conditional coverage
probabilities (w.r.t. $\protect\beta $ and $\protect\sigma $) of the
confidence intervals of \protect\cite{lee15exact}, with the LASSO model
selector where $\protect\lambda $ is either fixed or selected by cross
validation.}
\label{table:minimal:cond:coverage}
\end{table}

The conclusion of this comparison, and particularly of the evaluations of
minimal coverage probabilities, is that, although the confidence intervals
of \cite{lee15exact} are conditionally valid and convenient to compute,
their current applicability appears to be restricted to the case where the
tuning parameter $\lambda $ is fixed. [An extension of the strategy of \cite%
{lee15exact} to cross-validated versions of LASSO has recently been studied
in \cite{LoftTay15} and \cite{Loftus15}. This extension, however, comes with
significantly higher computational cost.] This also highlights the benefit
of the confidence intervals introduced in Section \ref%
{section:design:dependent}, which are intrinsically designed to be valid for
any model selection procedure whatsoever.

\section{Conclusion\label{conclusio}}

We have extended the PoSI confidence intervals of \cite{berk13valid} to PoSI
intervals for predictors. The coverage targets of our intervals, i.e., $%
x_{0}^{\prime }[\hat{M}]\beta _{\hat{M}}^{(n)}$ and $x_{0}^{\prime }[\hat{M}%
]\beta _{\hat{M}}^{(\star )}$, minimize a certain in-sample prediction error
and, under additional assumptions relating the training period to the
prediction period, a certain out-of-sample prediction error, respectively.
For in-sample prediction, i.e., for the target $x_{0}^{\prime }[\hat{M}%
]\beta _{\hat{M}}^{(n)}$, our intervals are valid, in finite samples,
irrespective of the model selection procedure that is being used. For
out-of-sample prediction, i.e., for the target $x_{0}^{\prime }[\hat{M}%
]\beta _{\hat{M}}^{(\star )}$, the same is true asymptotically under very
mild assumptions on the underlying model selector. See also \cite{Gree04a}
for optimality results related to the latter target and for its feasible
counterpart, under appropriate sparsity conditions.

Two types of confidence intervals were studied here: The first one
(corresponding to the constant $K_{1}(x_{0},\hat{M})$) depends on all
components of the vector $x_{0}$ (even if only a subset of these components
is `active' in the selected model $\hat{M}$) and thus is feasible only if $%
x_{0}$ is observed completely. The intervals of the second type
(corresponding to the constants $K_{2}(x_{0}[\hat{M}],\hat{M})$, $%
K_{3}(x_{0}[\hat{M}],\hat{M})$, and $K_{4}$) depend only on the active
components in the selected model, i.e., on $x_{0}[\hat{M}]$. The constants $%
K_{2}$, $K_{3}$, and $K_{4}$ correspond to successively larger confidence
intervals.

Computing the constant $K_{2}$ was found to be quite expensive in practice.
For computing the remaining constants, simple algorithms were presented. The
computational complexity of our algorithms for computing $K_{1}$ and $K_{3}$
is governed by the number of candidate models under consideration, limiting
computations to a few million candidate models in practice. Computation of $%
K_{4}$ is easy and not limited by complexity constraints (see, however, the
warning about numerical stability in Remark~\ref{stability} in Appendix \ref%
{section:practical:algorithms}). Our algorithms are of similar computational
complexity as those proposed in \cite{berk13valid}.

We furthermore have studied the behavior of the constants $K_{i}$ and of the
corresponding confidence intervals through analytic results in a setting
where model dimension is allowed to grow with sample size, and also through
simulations. These results provide evidence that $K_{4}$, which is
relatively cheap to compute, is a reasonably tight bound for the
computationally more expensive constants $K_{1}$ to $K_{3}$. Furthermore,
these results show that all the constants $K_{1}$ to $K_{4}$ are `bounded
away' from the Scheff\'{e} constant.

We have also provided simulation results regarding the coverage
probabilities of the various intervals introduced in the paper. We find that
the asymptotic results in Section \ref{section:distribution:dependent}
regarding the design-independent target already `kick-in' at moderate sample
sizes, and these results demonstrate that the PoSI confidence intervals for
the predictors maintain the desired minimal coverage probability. The
simulation study also shows that `naive' confidence intervals, which ignore
the data-driven model selection step and which use standard confidence
procedures as if the selected model were correct and given a priori, are
invalid also in the setting considered here (which is in line with earlier
findings in \cite{leeb13various}, where inter alia `naive' confidence
intervals for components of $\beta _{\hat{M}}^{(n)}$ were studied).
Furthermore, studying the confidence intervals developed for model selection
with the LASSO by \cite{lee15exact}, and others, we find that these
intervals are invalid if the LASSO penalty is chosen by cross-validation.
This contrasts the established fact that these intervals are valid
(conditionally on the event that a given model is selected), if the penalty
is fixed in advance.

\textbf{Acknowledgement. }We thank the referees and an Associate Editor for
thoughtful feedback and constructive comments. The first author acknowledges
constructive discussions with Lukas Steinberger and Nina Senitschnig on the
topic of the paper. The second author acknowledges partial support from the
Austrian Science Fund (FWF) projects P 28233-N32 and P 26354-N26.

\appendix

\section{Appendix: On the assumptions on $\hat{\protect\sigma}^{2}$\label%
{Section_on_sigma}}

In line with \cite{berk13valid} we have postulated the existence of an
estimator $\hat{\sigma}^{2}$ that is independent of $P_{X}Y$ and is
distributed as $\sigma ^{2}/r$ times a chi-square distributed random
variable with $r$ degrees of freedom ($1\leq r<\infty $). As already noted
in Section \ref{section:design:dependent}, if we assume that $d<n$ \emph{and}
$\mu =X\beta $ hold, such an estimator always exists and is given by the
usual residual variance estimator obtained from the residuals $Y-P_{X}Y$.
However, if $d=n$ holds (which typically is the case if $p>n$) or if $\mu $
is not known to belong to the column space of $X$, such an estimator is much
harder to come by.

Consider first the case where $d=n$ holds. Then it is plain that such an
estimator does not exist if it is to be only a function of $Y$ (and $X$):
Since here $P_{X}Y=Y$ holds, such an estimator would have to be independent
of $Y$ and thus constant with probability one, contradicting the requirement
to be distributed as a positive multiple of a chi-square. In order to
nevertheless be able to come up with an estimator $\hat{\sigma}^{2}$ with
the desired properties, one is hence forced to assume that one has access to 
\emph{additional} data beyond $Y$ that are related to $Y$ in an appropriate
way. A prototypical situation where such a construction is possible is as
follows: Assume that one has available additional data $Y^{\ast }$
distributed as $N(\mu ^{\ast },\sigma ^{2}I_{n^{\ast }})$, independently of $%
Y$ (for example, $Y^{\ast }$ might have been obtained from splitting the
original larger sample into $Y$ and $Y^{\ast }$). Assume further that for $%
Y^{\ast }$ one has available a (non-trivial) correct regression model (i.e., 
$\mu ^{\ast }=X^{\ast }\delta $ with $X^{\ast }$ of full column rank less
than $n^{\ast }$). Obviously, then an estimator $\hat{\sigma}^{2}$
satisfying all the required properties can be constructed from this correct
regression model for $Y^{\ast }$. However, this raises the question why one
would be willing to assume a correct regression model for one part of the
data, but would refuse to do so for the other part. [This might be defended
by reference to a structural break in the mean, which then however would beg
the question why the structural break would not also affect the variance $%
\sigma ^{2}$.] Alternatively to the assumption $\mu ^{\ast }=X^{\ast }\delta 
$, one could assume some `smoothness' in $\mu ^{\ast }$ and then use
nonparametric estimators to produce $\hat{\sigma}^{2}$. Again the question
arises why one would then not make a similar assumption for $\mu $ and use
the nonparametric method also for the first (or the entire) sample. In the
quite special situation where one has replicated observations in $Y^{\ast }$
available, one can abandon the dependence on a correct model (or on
smoothness assumptions) and nevertheless produce an estimator $\hat{\sigma}%
^{2}$ with the desired properties. All this granted, it seems that the
desired assumptions on $\hat{\sigma}^{2}$ and the desire to treat the case $%
d=n$ are not completely at ease.

Second, in case $d<n$ , but it is not assumed that $\mu =X\beta $ holds, it
is not obvious how an estimator $\hat{\sigma}^{2}$ with the desired
properties can be constructed without further assumptions (note that the
residual variance estimator obtained from $Y-P_{X}Y$ while being independent
from $P_{X}Y$ will in general \emph{not} be guaranteed to follow the
required distribution). One such assumption could be that we have available
a correct model $\mu =Z\gamma $, where the column space of $Z$ contains the
column space of $X$ with the rank of $Z$ still less than $n$; we could then
compute $\hat{\sigma}^{2}$ from this larger model (i.e., from $Y-P_{Z}Y$),
the resulting estimator having the desired properties. While this assumption
solves the existence problem for $\hat{\sigma}^{2}$, it raises the question
why one would then still want to keep the model selection exercise
restricted to submodels defined by the columns of $X$, when it is known that
the correct, larger, model $\mu =Z\gamma $ holds (and $Z$ is available).
Hence, we are led back essentially to the classical case with $Z$ playing
the role of $X$. Alternatively, the same constructions as in the preceding
paragraph relying on an independent sample $Y^{\ast }$ are available, but
they again suffer from the limitations pointed out before.

The discussion in this section shows that outside of the framework $d<n$ and 
$\mu =X\beta $ (in which case $p=d$ can be assumed with little loss of
generality) the assumptions on $\hat{\sigma}^{2}$ made in \cite{berk13valid}
as well as in the present paper are less than innocuous and will be
satisfied only in quite special situations; cf. also the discussion in
Remark 2.1(ii) in \cite{leeb13various}. [For this very reason, the first
version of this paper was set in the classical framework.] For ways to work
around this assumption on $\hat{\sigma}^{2}$ by using \textquotedblleft
conservative\textquotedblright\ estimators for the variance see \cite{BPS16}.

\section{Appendix: Proofs for Section \protect\ref{section:design:dependent} 
\label{app A}}

\begin{lem}
\label{lem:a1} Suppose $W$ is a random $m\times 1$ vector that has a density
that is positive almost everywhere. Let $a_{1},\ldots ,a_{L}$, for some $%
L\in \mathbb{N}$, be elements of $\mathbb{R}^{m}$, not all of which are
zero. Define $h(w)=\max_{l=1,\ldots ,L}\left\vert a_{l}^{\prime
}w\right\vert $, and set $H\left( t\right) =\Pr \left( h\left( W\right) \leq
t\right) $ for $t\in \mathbb{R}$. Then $H$ is continuous on $\mathbb{R}$,
satisfies $H\left( t\right) =0$ for $t\leq 0$, and is strictly increasing on 
$\left[ 0,\infty \right) $.
\end{lem}

\textbf{Proof:} For $t<0$ the event $\left\{ h\left( W\right) \leq t\right\} 
$ is empty; for $t=0$ this event is an intersection of the sets $\left\{
a_{l}^{\prime }W=0\right\} $ where at least one of these sets has
probability zero because $W$ possesses a density and not all $a_{l}$ are
zero. Consequently, $H\left( t\right) =0$ for $t\leq 0$ follows. Because $H$
is a distribution function, continuity of $H$ on $\mathbb{R}$ will follow if
we can establish continuity on $\left( 0,\infty \right) $. Now, for every $%
t>0$ the event $\left\{ h\left( W\right) =t\right\} $ is contained in the
union of the events $\left\{ \left\vert a_{l}^{\prime }W\right\vert
=t\right\} $ for which $a_{l}\neq 0$ holds. Since any of these events has
probability zero, it follows that $\Pr \left( \left\{ h\left( W\right)
=t\right\} \right) =0$ and consequently $H$ is continuous on $\left(
0,\infty \right) $. It remains to establish the claim regarding strict
monotonicity: For $t>0$ the set $A(t)=\left\{ w:h\left( w\right) \leq
t\right\} $ contains a sufficiently small ball centered at the origin
because $h\left( 0\right) =0$ and $h$ is continuous, and consequently $%
H\left( t\right) >0$ follows by the assumption on the density of $W$. It
hence suffices to show that $0<t_{1}<t_{2}$ implies $H\left( t_{1}\right)
<H\left( t_{2}\right) $. Because not all $a_{l}$ are zero and $h$ is
positively homogeneous of degree one, we can find an element $w_{1}\in
A(t_{1})$ such that $h\left( w_{1}\right) =t_{1}$ holds. But then there
exists an $l_{1}$ such that $\left\vert a_{l_{1}}^{\prime }w_{1}\right\vert
=t_{1}$ and $\left\vert a_{l}^{\prime }w_{1}\right\vert \leq t_{1}$ for all $%
l$ hold. In fact, we may assume that $a_{l_{1}}^{\prime }w_{1}=t_{1}$ holds
(otherwise we change the sign of $w_{1}$). Consider the set $B$ consisting
of all $w\in \mathbb{R}^{m}$ such that $a_{l_{1}}^{\prime }\left(
w-w_{1}\right) >0$ and such that $\left\vert a_{l}^{\prime }\left(
w-w_{1}\right) \right\vert <\left( t_{2}-t_{1}\right) /2$ for every $l$.
Then $B\subseteq A(t_{2})\backslash A(t_{1})$ holds, since for $w\in B$%
\begin{equation*}
h\left( w\right) \leq \max_{l=1,\ldots ,L}\left\vert a_{l}^{\prime }\left(
w-w_{1}\right) \right\vert +h\left( w_{1}\right) <\left( t_{2}-t_{1}\right)
/2+t_{1}=\left( t_{1}+t_{2}\right) /2<t_{2},
\end{equation*}%
\begin{equation*}
a_{l_{1}}^{\prime }w>a_{l_{1}}^{\prime }w_{1}=t_{1}>0,
\end{equation*}%
and hence also $h\left( w\right) \geq \left\vert a_{l_{1}}^{\prime
}w\right\vert >t_{1}$ hold. But $B$ obviously has positive Lebesgue measure,
implying that $H\left( t_{2}\right) -H\left( t_{1}\right) =\Pr \left(
A(t_{2})\backslash A(t_{1})\right) >0$. $\blacksquare $

\begin{rem}
\normalfont In the special case where $W=W_{1}W_{2}$ with $W_{1}$ a random $%
m\times 1$ vector having a density that is positive almost everywhere, with $%
W_{2}$ a random variable that is independent of $W_{1}$, is positive almost
surely, and has a density that is almost everywhere positive on $\left(
0,\infty \right) $, an alternative, and perhaps simpler, proof is as
follows: Set $H^{\ast }\left( t\right) =\Pr \left( h\left( W_{1}\right) \leq
t\right) $. We conclude that $H^{\ast }$ is continuous on $\mathbb{R}$ and
satisfies $H^{\ast }\left( t\right) =0$ for $t\leq 0$ by repeating the
corresponding arguments in the preceding proof. The same properties for $%
H\left( t\right) =\mathbb{E}_{W_{2}}H^{\ast }\left( t/W_{2}\right) $ then
follow immediately. To establish strict monotonicity of $H$ on $\left[
0,\infty \right) $ consider $0\leq t_{1}<t_{2}$. It is not difficult to see
that we can then find $w_{2}>0$ such that $H^{\ast }\left(
t_{1}/w_{2}\right) <H^{\ast }\left( t_{2}/w_{2}\right) $ holds since
otherwise $H^{\ast }$ would have to be constant on $\left[ 0,\infty \right) $
which is impossible since $H^{\ast }\left( 0\right) =0$ and $H^{\ast }$ is a
distribution function. By continuity of $H^{\ast }$ then also $H^{\ast
}\left( t_{1}/w_{2}^{\prime }\right) <H^{\ast }\left( t_{2}/w_{2}^{\prime
}\right) $ must hold for every $w_{2}^{\prime }$ in a sufficiently small
neighborhood of $w_{2}$. Since $H^{\ast }$ is nondecreasing and since the
distribution of $W_{2}$ puts positive mass on the aforementioned
neighborhood, we can conclude that $\mathbb{E}_{W_{2}}H^{\ast }\left(
t_{1}/W_{2}\right) <\mathbb{E}_{W_{2}}H^{\ast }\left( t_{2}/W_{2}\right) $,
i.e., that $H\left( t_{1}\right) <H\left( t_{2}\right) $ holds.
\end{rem}

The following lemma will be used in the proof of Proposition \ref%
{prop:PoSI:upper:bound} below.

\begin{lem}
\label{lem:a2} Suppose $F^{\ast }$ is a distribution function on $\mathbb{R}$
that is continuous at zero. Let $S$ be a random variable that is positive
with probability one and has a continuous distribution function. Then $%
F\left( t\right) =\mathbb{E}_{S}F^{\ast }\left( t/S\right) $ is continuous
on $\mathbb{R}$.
\end{lem}

\textbf{Proof: }Let $S^{\ast }$ be a random variable which is independent of 
$S$ and which has distribution function $F^{\ast }$. Then $F(t)=\mathbb{E}%
_{S}\mathbb{E}_{S^{\ast }}\boldsymbol{1}\left( S^{\ast }\leq t/S\right) =\Pr
\left( SS^{\ast }\leq t\right) $. Because $S^{\ast }\neq 0$ holds almost
surely by the assumption on $F^{\ast }$, we have $\Pr \left( SS^{\ast
}=t\right) =\mathbb{E}_{S^{\ast }}\mathbb{E}_{S}\boldsymbol{1}\left(
S=t/S^{\ast }\right) $. Since $S$ has a continuous distribution function, we
have $\mathbb{E}_{S}\boldsymbol{1}\left( S=t/S^{\ast }\right) =0$ almost
surely, implying that $\Pr \left( SS^{\ast }=t\right) =0$. $\blacksquare $

\textbf{Proof of Proposition \ref{prop:PoSI:upper:bound}:} (a) Observe that
in case $M=\left\{ 1,\ldots ,p\right\} $ we have $c(M,\mathcal{M})=0$ and
that $x_{0}=0$ implies $\bar{s}_{M_{\ast }}=0$ for every $M_{\ast }$; thus
it is obvious that $F_{M,x_{0}}^{\ast }$, and hence also $F_{M,x_{0}}$, is
the indicator function of $\left[ 0,\infty \right) $, which then implies
that $K_{3}$ exists and is well-defined in this case and that $K_{3}(x_{0}{%
[M],M)=0}$. If $M=\left\{ 1,\ldots ,p\right\} $ but $x_{0}\neq 0$, then $%
F_{M,x_{0}}$ is continuous on $\mathbb{R}$, satisfies $F_{M,x_{0}}\left(
t\right) =0$ for $t\leq 0$, and is strictly increasing on $\left[ 0,\infty
\right) $ in view of Lemma \ref{lem:a1}, since in this case $%
F_{M,x_{0}}\left( t\right) $ reduces to $P_{n,\mu ,\sigma }\left(
\max_{M_{\ast }\in \mathcal{M}}\left\vert \bar{s}_{M_{\ast }}^{\prime
}\left( Y-\mu \right) \right\vert /\hat{\sigma}\leq t\right) $ (see (\ref%
{bound_1}) and (\ref{expression}) below) and since not all $\bar{s}_{M_{\ast
}}$ can be zero (in view of our assumptions on $\mathcal{M}$). Again the
claims then clearly follow in this case.

In case $M\in \mathcal{M}$ is a proper subset of $\left\{ 1,\ldots
,p\right\} $ and $p>1$ holds we argue as follows: Note that then $%
F_{M,x_{0}}^{\ast }\left( 0\right) =0$ holds since now $c(M,\mathcal{M})\geq
1$ holds and since $F_{Beta,1/2,(d-1)/2}\left( 0\right) =0$. Hence $%
F_{M,x_{0}}^{\ast }$ is continuous at $t=0$. We may apply Lemma \ref{lem:a2}
to conclude that $F_{M,x_{0}}$ is continuous on $\mathbb{R}$ and thus
satisfies $F_{M,x_{0}}\left( 0\right) =0$ (since $F_{M,x_{0}}\left( t\right)
=0$ for $t<0$ by its definition). Next let $0\leq t_{1}<t_{2}$. Because $%
F_{M,x_{0}}^{\ast }\left( 0\right) =0$ as noted before and because $%
F_{M,x_{0}}^{\ast }\left( 1\right) =1$ (since $\Pr \left( \max_{M_{\ast
}\subseteq M}\left\vert \bar{s}_{M_{\ast }}^{\prime }V\right\vert >1\right)
=0$ and $F_{Beta,1/2,(d-1)/2}\left( 1\right) =1$) we thus can find a
positive $g_{0}$ such that $F_{M,x_{0}}^{\ast }\left( t_{1}/g_{0}\right)
<F_{M,x_{0}}^{\ast }\left( t_{2}/g_{0}\right) $ holds (if not, constancy of $%
F_{M,x_{0}}^{\ast }$ on $\left[ 0,\infty \right) $ would have to follow).
Because of continuity from the right at $t_{1}/g_{0}$ it follows that $%
F_{M,x_{0}}^{\ast }\left( t_{1}/g\right) <F_{M,x_{0}}^{\ast }\left(
t_{2}/g\right) $ also holds for all $g<g_{0}$ in a sufficiently small
neighborhood of $g_{0}$ that is contained in $\left( 0,\infty \right) $.
Because $F_{M,x_{0}}^{\ast }\left( t_{1}/g\right) \leq F_{M,x_{0}}^{\ast
}\left( t_{2}/g\right) $ holds for every $g>0$ and because $G$ has a density
that is positive everywhere on $\left( 0,\infty \right) $, the strict
inequality $\mathbb{E}_{G}F_{M,x_{0}}^{\ast }\left( t_{1}/G\right) <\mathbb{E%
}_{G}F_{M,x_{0}}^{\ast }\left( t_{2}/G\right) $ follows. This establishes
strict monotonicity of $F_{M,x_{0}}$ on $\left[ 0,\infty \right) $ also in
this case, which proves the claims.

Finally, if $M\in \mathcal{M}$ is a proper subset of $\left\{ 1,\ldots
,p\right\} $ and $p=1$ holds, then $M$ is empty and $d=1$ must hold, and
hence $F_{M,x_{0}}^{\ast }$ reduces to the indicator function of $\left[
1,\infty \right) $. But then $F_{M,x_{0}}\left( t\right) =\Pr \left( G\leq
t\right) $ which obviously is continuous on $\mathbb{R}$, takes the value
zero at $t=0$, and is strictly increasing on $\left[ 0,\infty \right) $,
again implying the claims.

(b) Observe that $\bar{s}_{M_{\ast }}$ belongs to the column space of $X$
for every $M_{\ast }\in \mathcal{M}$ and hence we have%
\begin{eqnarray}
&&P_{n,\mu ,\sigma }\left( \max_{M_{\ast }\in \mathcal{M}}\left\vert \bar{s}%
_{M_{\ast }}^{\prime }\left( Y-\mu \right) \right\vert /\hat{\sigma}>t\right)
\notag \\
&=&P_{n,\mu ,\sigma }\left( \max_{M_{\ast }\in \mathcal{M}}\left\vert \bar{s}%
_{M_{\ast }}^{\prime }P_{X}\left( Y-\mu \right) /\left\Vert P_{X}\left(
Y-\mu \right) \right\Vert \right\vert >\left( \hat{\sigma}/\left\Vert
P_{X}\left( Y-\mu \right) \right\Vert \right) t\right) ,\text{ \ \ \ \ \ \ }
\label{bound_1}
\end{eqnarray}%
where $P_{X}\left( Y-\mu \right) /\left\Vert P_{X}\left( Y-\mu \right)
\right\Vert $ and $\left\Vert P_{X}\left( Y-\mu \right) \right\Vert /\hat{%
\sigma}$ are independent since the random variables $P_{X}\left( Y-\mu
\right) /\left\Vert P_{X}\left( Y-\mu \right) \right\Vert $, $\left\Vert
P_{X}\left( Y-\mu \right) \right\Vert $, and $\hat{\sigma}$ are mutually
independent. [Observe that $P_{X}\left( Y-\mu \right) $ is nonzero with
probability $1$ since $d\geq 1$ holds.] Consequently, the probability given
above can be represented as 
\begin{equation}
\Pr \left( \max_{M_{\ast }\in \mathcal{M}}\left\vert \bar{s}_{M_{\ast
}}^{\prime }V\right\vert >t/G\right)  \label{expression}
\end{equation}%
where $V$ and $G$ are independent and otherwise are as in the definition of $%
F_{M,x_{0}}^{\ast }$ and $F_{M,x_{0}}$. Now, using first independence of $V$
and $G$ and then a union bound twice we have for $M\in \mathcal{M}$ and $%
t\geq 0$%
\begin{eqnarray}
&&\Pr \left( \max_{M_{\ast }\in \mathcal{M}}\left\vert \bar{s}_{M_{\ast
}}^{\prime }V\right\vert >t/G\right) =\int \Pr \left( \max_{M_{\ast }\in 
\mathcal{M}}\left\vert \bar{s}_{M_{\ast }}^{\prime }V\right\vert >t/g\right)
dF_{G}(g)  \notag \\
&\leq &\int \min \left[ 1,\Pr \left( \max_{M_{\ast }\in \mathcal{M},M_{\ast
}\subseteq M}\left\vert \bar{s}_{M_{\ast }}^{\prime }V\right\vert
>t/g\right) +\Pr \left( \max_{M_{\ast }\in \mathcal{M},M_{\ast }\nsubseteqq
M}\left\vert \bar{s}_{M_{\ast }}^{\prime }V\right\vert >t/g\right) \right]
dF_{G}(g)  \notag \\
&\leq &\int \min \left[ 1,\Pr \left( \max_{M_{\ast }\in \mathcal{M},M_{\ast
}\subseteq M}\left\vert \bar{s}_{M_{\ast }}^{\prime }V\right\vert
>t/g\right) +\sum_{M_{\ast }\in \mathcal{M},M_{\ast }\nsubseteqq M}\Pr
\left( \left\vert \bar{s}_{M_{\ast }}^{\prime }V\right\vert >t/g\right) %
\right] dF_{G}(g)  \notag \\
&=&\int \min \left[ 1,\Pr \left( \max_{M_{\ast }\in \mathcal{M},M_{\ast
}\subseteq M}\left\vert \bar{s}_{M_{\ast }}^{\prime }V\right\vert
>t/g\right) +\sum_{M_{\ast }\in \mathcal{M},M_{\ast }\nsubseteqq M}\Pr
\left( \left( \bar{s}_{M_{\ast }}^{\prime }V\right) ^{2}>t^{2}/g^{2}\right) %
\right] dF_{G}(g)  \notag \\
&\leq &\int \left( 1-F_{M,x_{0}}^{\ast }\left( t/g\right) \right) dF_{G}(g)=%
\mathbb{E}_{G}\left( 1-F_{M,x_{0}}^{\ast }\left( t/G\right) \right)
=1-F_{M,x_{0}}\left( t\right) ,  \label{bound_2}
\end{eqnarray}%
where $F_{G}$ here denotes the c.d.f. of $G$. The last inequality follows
from the fact that $\Pr \left( \left( \bar{s}_{M_{\ast }}^{\prime }V\right)
^{2}>t^{2}/g^{2}\right) $ is either equal to zero (if $\bar{s}_{M_{\ast }}=0$%
) or is equal to $1-F_{Beta,1/2,(d-1)/2}\left( t^{2}/g^{2}\right) $ (if $%
\bar{s}_{M_{\ast }}\neq 0$) as is easy to see; for the case where $M$ is the
empty set also observe that $\Pr \left( \max_{M_{\ast }\in \mathcal{M}%
,M_{\ast }\subseteq M}\left\vert \bar{s}_{M_{\ast }}^{\prime }V\right\vert
>t/g\right) =0$ for $t\geq 0$ because $\bar{s}_{\varnothing }=0$. In view of
(\ref{eq:def:posi:xzero}) the chain of inequalities in (\ref{bound_1})-(\ref%
{bound_2}) establishes $K_{1}(x_{0})\leq K_{3}(x_{0}[M],M)$. It follows that 
$K_{1}(x)\leq K_{3}(x[M],M)=K_{3}(x_{0}[M],M)$ for every $x$ satisfying $%
x[M]=x_{0}[M]$, implying $K_{2}(x_{0}[M],M)\leq K_{3}(x_{0}[M],M)$. The
inequality (\ref{eq:monoton:1}) is obvious and inequality (\ref{eq:monoton:2}%
) follows since for $t\geq 0$ we have (again noting that expressions like $%
\Pr \left( \max_{M_{\ast }\in \mathcal{M},M_{\ast }\subseteq
M_{1}}\left\vert \bar{s}_{M_{\ast }}^{\prime }V\right\vert >t\right) $ for $%
t\geq 0$ are equal to zero if $M_{1}$ is empty)%
\begin{eqnarray*}
&&\Pr \left( \max_{M_{\ast }\in \mathcal{M},M_{\ast }\subseteq
M_{2}}\left\vert \bar{s}_{M_{\ast }}^{\prime }V\right\vert >t\right)
+c(M_{2},\mathcal{M})\left( 1-F_{Beta,1/2,(d-1)/2}\left( t^{2}\right) \right)
\\
&\leq &\Pr \left( \max_{M_{\ast }\in \mathcal{M},M_{\ast }\subseteq
M_{1}}\left\vert \bar{s}_{M_{\ast }}^{\prime }V\right\vert >t\right) +\Pr
\left( \max_{M_{\ast }\in \mathcal{M},M_{\ast }\subseteq M_{2},M_{\ast
}\nsubseteqq M_{1}}\left\vert \bar{s}_{M_{\ast }}^{\prime }V\right\vert
>t\right) \\
&&+c(M_{2},\mathcal{M})\left( 1-F_{Beta,1/2,(d-1)/2}\left( t^{2}\right)
\right) \\
&\leq &\Pr \left( \max_{M_{\ast }\in \mathcal{M},M_{\ast }\subseteq
M_{1}}\left\vert \bar{s}_{M_{\ast }}^{\prime }V\right\vert >t\right)
+\sum_{M_{\ast }\in \mathcal{M},M_{\ast }\subseteq M_{2},M_{\ast
}\nsubseteqq M_{1}}\Pr \left( \left\vert \bar{s}_{M_{\ast }}^{\prime
}V\right\vert >t\right) \\
&&+c(M_{2},\mathcal{M})\left( 1-F_{Beta,1/2,(d-1)/2}\left( t^{2}\right)
\right) \\
&\leq &\Pr \left( \max_{M_{\ast }\in \mathcal{M},M_{\ast }\subseteq
M_{1}}\left\vert \bar{s}_{M_{\ast }}^{\prime }V\right\vert >t\right)
+c(M_{1},\mathcal{M})\left( 1-F_{Beta,1/2,(d-1)/2}\left( t^{2}\right)
\right) .
\end{eqnarray*}%
The relation $K_{3}(x_{0}[M],M)\leq K_{4}$ is now immediate. Finally, $%
1-F_{\varnothing ,x_{0}}^{\ast }\left( t\right) \leq 1$ for all $t\in 
\mathbb{R}$ and $1-F_{\varnothing ,x_{0}}^{\ast }\left( t\right) =0$ for $%
t>1 $ lead to%
\begin{eqnarray*}
1-F_{\varnothing ,x_{0}}\left( t\right) &=&\mathbb{E}_{G}\left(
1-F_{\varnothing ,x_{0}}^{\ast }\left( t/G\right) \right) =\mathbb{E}%
_{G}\left( \left( 1-F_{\varnothing ,x_{0}}^{\ast }\left( t/G\right) \right) 
\boldsymbol{1}\left( t\leq G\right) \right) \\
&\leq &\mathbb{E}_{G}\boldsymbol{1}\left( t\leq G\right) =1-\Pr \left( G\leq
t\right) ,
\end{eqnarray*}%
which proves $K_{4}\leq K_{5}$. $\blacksquare $

\begin{lem}
\label{lem:case:p:egal:2} Assume $p=2$ and $n\geq 2$. Then there exists a
design matrix $X$ with full column rank and a vector $x_{0}$ such that $%
K_{4}=K_{1}(x_{0})$ for $\mathcal{M}$ the power set of $\left\{ 1,2\right\} $%
.
\end{lem}

\textbf{Proof:} Assume first that $n=2$. In view of the definition of $%
K_{4}=K_{3}\left( x_{0}[\varnothing ],\varnothing \right) $ it suffices to
exhibit a $2\times 2$ matrix $X$ and a $2\times 1$ vector $x_{0}$ such that
equality holds between the far l.h.s. and the far r.h.s. of (\ref{bound_2})
for $M=\varnothing $ and all $t\geq 0$. Inspection of (\ref{bound_2}) shows
that for this it suffices to find $X$ and $x_{0}$ such that%
\begin{equation*}
\Pr \left( \max_{\varnothing \neq M_{\ast }\subseteq \left\{ 1,2\right\}
}\left\vert \bar{s}_{M_{\ast }}^{\prime }V\right\vert >c\right) =\min \left(
1,\sum_{\varnothing \neq M_{\ast }\subseteq \left\{ 1,2\right\} }\Pr \left(
\left\vert \bar{s}_{M_{\ast }}^{\prime }V\right\vert >c\right) \right)
\end{equation*}%
holds for every $c\geq 0$ and that $\bar{s}_{M_{\ast }}^{\prime }\neq 0$ for
every $\varnothing \neq M_{\ast }\subseteq \left\{ 1,2\right\} $. This is
achieved for 
\begin{equation*}
X=X^{(2)}=\left[ 
\begin{array}{cc}
1 & \cos \left( 2\pi /3\right) \\ 
0 & \sin \left( 2\pi /3\right)%
\end{array}%
\right]
\end{equation*}%
and $x_{0}^{\prime }=x_{0}^{(2)\prime }=\left( \cos \left( 4\pi /3\right)
,\sin \left( 4\pi /3\right) \right) X^{(2)}$: Then $\bar{s}_{\left\{
1\right\} }^{\prime }=-\left( 1,0\right) $, $\bar{s}_{\left\{ 2\right\}
}^{\prime }=-\left( \cos \left( 2\pi /3\right) ,\sin \left( 2\pi /3\right)
\right) $, and $\bar{s}_{\left\{ 1,2\right\} }^{\prime }=\left( \cos \left(
4\pi /3\right) ,\sin \left( 4\pi /3\right) \right) $. Consequently, the
event $\left\{ \max_{\varnothing \neq M_{\ast }\subseteq \left\{ 1,2\right\}
}\left\vert \bar{s}_{M_{\ast }}^{\prime }V\right\vert >c\right\} $ is either
the entire space or is the disjoint union of the events $\{\left\vert \bar{s}%
_{\left\{ 1\right\} }^{\prime }V\right\vert >c\}$, $\{\left\vert \bar{s}%
_{\left\{ 2\right\} }^{\prime }V\right\vert >c\}$ and $\{\left\vert \bar{s}%
_{\left\{ 1,2\right\} }^{\prime }V\right\vert >c\}$. In the case $n>2$
simply set 
\begin{equation*}
X=\left( X^{(2)\prime },0,\ldots ,0\right) ^{\prime }
\end{equation*}%
and $x_{0}^{\prime }=\left( \cos \left( 4\pi /3\right) ,\sin \left( 4\pi
/3\right) ,0,\ldots ,0\right) X$. $\blacksquare $

\begin{rem}
\normalfont Further examples of pairs $X$, $x_{0}$ satisfying the above
lemma can be generated from the matrix $X$ constructed in the proof by
premultiplying $X$ by an orthogonal matrix and leaving $x_{0}$ unchanged.
\end{rem}

\textbf{Proof of Proposition \ref{constants_1}: }(a) The distribution of $%
\omega =\max_{M\in \mathcal{M}}\left\vert \bar{s}_{M}^{\prime }\left( Y-\mu
\right) \right\vert /\sigma $ clearly does not change if $X$ is replaced by $%
AX$, where $A$ is an orthogonal $n\times n$ matrix. Furthermore, scaling the
columns of $X$ and the corresponding columns of $x_{0}^{\prime }$ by the
same (column-specific) positive constants does not alter $\max_{M\in 
\mathcal{M}}\left\vert \bar{s}_{M}^{\prime }\left( Y-\mu \right) \right\vert
/\sigma $. Hence, we may assume w.l.o.g. that $X$ consist of the first $p$
standard basis vectors of $\mathbb{R}^{n}$. Then choose $x_{0}^{\prime }$ as
the $1\times p$ vector $(1,\ldots ,1)$. It follows that $\omega $ can be
written as%
\begin{equation*}
\max_{M\in \mathcal{M}}\left\vert \sum_{i\in M}Z_{i}\right\vert /\sqrt{%
\left\vert M\right\vert }
\end{equation*}%
where $Z_{i}$ are i.i.d. standard normal and where we use the convention
that the expression in the display is zero if $\left\vert M\right\vert =0$.
For a positive real $b$ define the random set $\tilde{M}=\left\{ i\in
\{1,...,p\}:Z_{i}\geq b\right\} $. Since for any realization of the random
variables $Z_{i}$ we have that $\tilde{M}\in \mathcal{M}$, we must have%
\begin{equation*}
\omega /\sqrt{p}\geq p^{-1/2}\left\vert \sum_{i\in \tilde{M}%
}Z_{i}\right\vert /\sqrt{\left\vert \tilde{M}\right\vert }=\left\vert
p^{-1}\sum_{i=1}^{p}Z_{i}\boldsymbol{1}(Z_{i}\geq b)\right\vert /\sqrt{%
p^{-1}\sum_{i=1}^{p}\boldsymbol{1}(Z_{i}\geq b)}.
\end{equation*}%
By the law of large numbers we obtain that the r.h.s. converges to $\phi (b)/%
\sqrt{1-\Phi (b)}$ almost surely. Because $K_{1}(x_{0},\infty )/\sqrt{p}$ is
the $(1-\alpha )$-quantile of $\omega /\sqrt{p}$ with $\alpha $ independent
of $p$ and since $b>0$ was arbitrary, the first claim follows. The second
claim follows immediately by choosing $x_{0}$ equal to a $p\times 1$
standard basis vector and by noting that then $\omega $ is distributed as
the absolute value of a standard normal variable. [More generally, if $%
\sup_{p}\left\Vert x_{0}\right\Vert _{0}<\infty $ holds, then $\omega $ is a
maximum absolute value of at most $k$ standard normal variables, where $k$
is a finite number not depending on $p$.]

(b) For the same reasons as given at the beginning of the proof of part (a)
we have for every $p$ that 
\begin{eqnarray}
&&\inf_{x_{0}\in \mathbb{R}^{p}}\inf_{X\in \mathsf{X}(p)}\inf_{M\in \mathcal{%
M},|M|\leq \gamma p}K_{2}(x_{0}[M],M,\infty ,X,\alpha ,\mathcal{M})/\sqrt{p}
\notag \\
&=&\inf_{x_{0}\in \mathbb{R}^{p}}\inf_{n\geq p}\inf_{M\in \mathcal{M}%
,|M|\leq \gamma p}K_{2}(x_{0}[M],M,\infty ,(I_{p},0_{p\times (n-p)})^{\prime
},\alpha ,\mathcal{M})/\sqrt{p}  \label{first_display} \\
&=&\inf_{x_{0}\in \mathbb{R}^{p}}\inf_{M\in \mathcal{M},|M|\leq \gamma
p}K_{2}(x_{0}[M],M,\infty ,I_{p},\alpha ,\mathcal{M})/\sqrt{p}  \notag
\end{eqnarray}%
where $I_{p}$ is the identity matrix of dimension $p$. By the monotonicity
property (\ref{eq:monoton:1}) the far r.h.s. of (\ref{first_display}) equals%
\begin{equation}
\inf_{x_{0}\in \mathbb{R}^{p}}\inf_{M\in \mathcal{M},|M|=\left\lfloor \gamma
p\right\rfloor }K_{2}(x_{0}[M],M,\infty ,I_{p},\alpha ,\mathcal{M})/\sqrt{p}.
\label{a_display}
\end{equation}%
Now fix an arbitrary $x_{0}\in \mathbb{R}^{p}$ and $M\in \mathcal{M}$ with $%
|M|=\left\lfloor \gamma p\right\rfloor $. Define $x_{0}^{\ast }$ via $%
x_{0i}^{\ast }=x_{0i}$ for $i\in M$ and set $x_{0i}^{\ast }=1$ else. Then 
\begin{equation}
K_{2}(x_{0}[M],M,\infty ,I_{p},\alpha ,\mathcal{M})\geq K_{1}(x_{0}^{\ast
},\infty ,I_{p},\alpha ,\mathcal{M})  \label{second_display}
\end{equation}%
holds and the latter quantity is the $(1-\alpha )$-quantile of 
\begin{equation*}
\omega ^{\ast }=\max_{M^{\ast }\in \mathcal{M}}\left\vert \sum_{i\in M^{\ast
}}x_{0i}^{\ast }Z_{i}\right\vert /\sqrt{\sum_{i\in M^{\ast }}\left(
x_{0i}^{\ast }\right) ^{2}}.
\end{equation*}%
For a positive real $b$ define now the random set $\breve{M}=\left\{ i\notin
M:Z_{i}\geq b\right\} $. Similar as above we then conclude that 
\begin{equation*}
\omega ^{\ast }/\sqrt{p}\geq \left\vert p^{-1}\sum_{i\notin M}Z_{i}%
\boldsymbol{1}(Z_{i}\geq b)\right\vert /\sqrt{p^{-1}\sum_{i\notin M}%
\boldsymbol{1}(Z_{i}\geq b)}.
\end{equation*}%
While the r.h.s. of the above display depends on $M$, its distribution does
not as it coincides with the distribution of 
\begin{equation*}
A_{p}=\left\vert p^{-1}\sum_{i=\left\lfloor \gamma p\right\rfloor
+1}^{p}Z_{i}\boldsymbol{1}(Z_{i}\geq b)\right\vert /\sqrt{%
p^{-1}\sum_{i=\left\lfloor \gamma p\right\rfloor +1}^{p}\boldsymbol{1}%
(Z_{i}\geq b)}.
\end{equation*}%
Consequently, $K_{1}(x_{0}^{\ast },\infty ,I_{p},\alpha ,\mathcal{M})/\sqrt{p%
}$ as the $(1-\alpha )$-quantile of $\omega ^{\ast }/\sqrt{p}$ is not
smaller than the corresponding quantile of $A_{p}$, and this is true
independently of the choice of $x_{0}$ and of $M$ with $|M|=\left\lfloor
\gamma p\right\rfloor $. Since $A_{p}$ converges to $\sqrt{1-\gamma }\phi
(b)/\sqrt{1-\Phi (b)}$ almost surely and $\alpha $ does not depend on $p$,
we can conclude that%
\begin{equation*}
\liminf_{p\rightarrow \infty }\inf_{x_{0}\in \mathbb{R}^{p}}\inf_{M\in 
\mathcal{M},|M|=\left\lfloor \gamma p\right\rfloor }K_{1}(x_{0}^{\ast
},\infty ,I_{p},\alpha ,\mathcal{M})/\sqrt{p}\geq \sqrt{1-\gamma }\phi (b)/%
\sqrt{1-\Phi (b)}.
\end{equation*}%
Since $b>0$ was arbitrary, the proof is then complete in view of (\ref%
{first_display}), (\ref{a_display}), and (\ref{second_display}). $%
\blacksquare $

\begin{lem}
\label{quantiles}Let $\alpha $, $0<\alpha <1$, be a fixed number, let $N\in 
\mathbb{N}$, and $d\in \mathbb{N}$. Let $K(N,d,\alpha )$ denote the ($%
1-\alpha $)-quantile of the distribution function given by%
\begin{equation*}
1-\mathbb{E}_{G}\left( \min \left[ 1,N\left( 1-F_{Beta,1/2,\left( d-1\right)
/2}\left( t^{2}/G^{2}\right) \right) \right] \right)
\end{equation*}%
for $t\geq 0$ and by $0$ for $t<0$. Here $G$ is a nonnegative random
variable such that $G^{2}$ follows a chi-square distribution with $d$
degrees of freedom. Then as $\min (N,d)\rightarrow \infty $%
\begin{equation*}
K(N,d,\alpha )/\sqrt{d\left( 1-N^{-2/\left( d-1\right) }\right) }\rightarrow
1\text{.}
\end{equation*}
\end{lem}

\textbf{Proof:} The c.d.f. in the lemma is the c.d.f. of $GW$, where $W$ is
independent of $G$, is nonnegative, and has distribution function given by 
\begin{equation*}
1-\min \left[ 1,N\left( 1-F_{Beta,1/2,\left( d-1\right) /2}\left(
t^{2}\right) \right) \right]
\end{equation*}%
for $t\geq 0$. Observe that $G/\sqrt{d}$ converges to $1$ in probability as $%
d\rightarrow \infty $. To complete the proof it thus suffices to show that $%
W/\sqrt{\left( 1-N^{-2/\left( d-1\right) }\right) }$ converges to $1$ in
probability as $\min (N,d)\rightarrow \infty $: For $t>1$ we have that%
\begin{equation*}
\Pr \left( W/\sqrt{\left( 1-N^{-2/\left( d-1\right) }\right) }>t\right) \leq
N\left( 1-F_{Beta,1/2,\left( d-1\right) /2}\left( t^{2}\left( 1-N^{-2/\left(
d-1\right) }\right) \right) \right) .
\end{equation*}%
But the r.h.s. of the preceding display has been shown in \cite{zhang13rank}
to converge to zero as $\min (N,d)\rightarrow \infty $, cf. (A.5) and (A.6)
in that paper. For $t<1$ we have%
\begin{equation*}
\Pr \left( W/\sqrt{\left( 1-N^{-2/\left( d-1\right) }\right) }\leq t\right)
=1-\min \left[ 1,N\left( 1-F_{Beta,1/2,\left( d-1\right) /2,}\left(
t^{2}\left( 1-N^{-2/\left( d-1\right) }\right) \right) \right) \right] ,
\end{equation*}%
and hence it suffices to show that 
\begin{equation*}
N\left( 1-F_{Beta,1/2,\left( d-1\right) /2,}\left( t^{2}\left(
1-N^{-2/\left( d-1\right) }\right) \right) \right) \rightarrow \infty
\end{equation*}%
as $\min (N,d)\rightarrow \infty $. But this has again been established in 
\cite{zhang13rank}, see (C.4) and (C.5) in that paper. $\blacksquare $

\textbf{Proof of Proposition \ref{constants_2}:} Observe that $K_{4}(\infty
) $ is always positive and hence the ratio in (\ref{Closeness of K3 and K4})
is well-defined. In view of (12) and (13) and the assumptions on $\mathcal{M}
$ we have for $M\in \mathcal{M}$ with $M\neq \left\{ 1,\ldots ,p\right\} $,
for $X\in \mathsf{X}_{n(p),p}(\mathcal{M})$, and for $x_{0}\in \mathbb{R}%
^{p} $ that $K_{3}(x_{0}[M],M,\infty )$ is not less than $K(\left\lfloor
\tau \left\vert \mathcal{M}\right\vert \right\rfloor ,\min (n(p),p),\alpha )$
in the notation of Lemma \ref{quantiles}, where we note that $\left\lfloor
\tau \left\vert \mathcal{M}\right\vert \right\rfloor \geq 1$ holds at least
from a certain $p$ onwards. Since $K_{3}(x_{0}[M],M,\infty )\leq
K_{4}(\infty )$ always holds, and since $K_{4}(\infty )=K(\left\vert 
\mathcal{M}\right\vert -1,\min (n(p),p),\alpha )$ in the notation of Lemma %
\ref{quantiles}, it suffices to show that 
\begin{equation*}
K(\left\lfloor \tau \left\vert \mathcal{M}\right\vert \right\rfloor ,\min
(n(p),p),\alpha )/K(\left\vert \mathcal{M}\right\vert -1,\min
(n(p),p),\alpha )\rightarrow 1
\end{equation*}%
as $p\rightarrow \infty $. Note that $d=\min (n(p),p)\rightarrow \infty $ as 
$p\rightarrow \infty $ by the assumption on $n(p)$, and that $\left\vert 
\mathcal{M}\right\vert \rightarrow \infty $ as $p\rightarrow \infty $. By
Lemma \ref{quantiles} we thus need to show that%
\begin{equation*}
A_{d}\left( \left\vert \mathcal{M}\right\vert \right) :=\left( 1-(\left\vert 
\mathcal{M}\right\vert -1)^{-2/\left( d-1\right) }\right) /\left(
1-\left\lfloor \tau \left\vert \mathcal{M}\right\vert \right\rfloor
^{-2/\left( d-1\right) }\right) \rightarrow 1
\end{equation*}%
as $p\rightarrow \infty $. Observe that $\tau <1$ must hold, and thus $%
\left\vert \mathcal{M}\right\vert -1\geq \left\lfloor \tau \left\vert 
\mathcal{M}\right\vert \right\rfloor >1$ holds for large $p$. This, in
particular, implies $A_{d}\left( \left\vert \mathcal{M}\right\vert \right)
\geq 1$ for large $p$. But then for large $p$%
\begin{equation*}
1\leq A_{d}\left( \left\vert \mathcal{M}\right\vert \right) \leq \left(
1-\left\vert \mathcal{M}\right\vert ^{-2/\left( d-1\right) }\right) /\left(
1-\left( (\tau /2)\left\vert \mathcal{M}\right\vert \right) ^{-2/\left(
d-1\right) }\right) =:B_{d}\left( \left\vert \mathcal{M}\right\vert \right)
\end{equation*}%
holds since also $\left\lfloor \tau \left\vert \mathcal{M}\right\vert
\right\rfloor \geq (\tau /2)\left\vert \mathcal{M}\right\vert >1$ is
satisfied for large $p$. It thus suffices to show that $B_{d}\left(
\left\vert \mathcal{M}\right\vert \right) \rightarrow 1$ for $p\rightarrow
\infty $. Let $f>2/\tau $ be a real number. Then $\left\vert \mathcal{M}%
\right\vert \geq f$ holds for large $p$. Because $B_{d}(x)$ is monotone
decreasing in $x$ for $x>2/\tau >1$ and for every $d$ as is easily checked
by inspection of the derivative, we have that $B_{d}\left( f\right) \geq
B_{d}\left( \left\vert \mathcal{M}\right\vert \right) $ holds for large $p$.
But now it is easily checked (H\^{o}pital's rule) that $B_{d}\left( f\right) 
$ converges to $\log f/\left( \log f+\log (\tau /2)\right) $ as $%
p\rightarrow \infty $ (and thus $d\rightarrow \infty $). Making $f$
arbitrarily large, $\log f/\left( \log f+\log (\tau /2)\right) $ approaches $%
1$. This completes the proof of (\ref{Closeness of K3 and K4}). The second
claim follows immediately from $K_{4}=K(\left\vert \mathcal{M}\right\vert
-1,\min (n(p),p),\alpha )$, the preceding lemma, the observation that $%
B_{d}\left( \left\vert \mathcal{M}\right\vert \right) \rightarrow 1$, and
that%
\begin{equation*}
1/B_{d}\left( \left\vert \mathcal{M}\right\vert \right) \leq \left(
1-(\left\vert \mathcal{M}\right\vert -1)^{-2/\left( d-1\right) }\right)
/\left( 1-\left\vert \mathcal{M}\right\vert ^{-2/\left( d-1\right) }\right)
\leq 1
\end{equation*}%
holds. $\blacksquare $

\textbf{Proof of Corollary \ref{constants_3}:} Properties (i), (ii), and
(iv) are obvious. In case $m_{p}=p$, we have that $\mathcal{M}(m_{p})$ is
the power set and hence $\left\vert \mathcal{M}(m_{p})\right\vert
=2^{m_{p}}=2^{p}$. But then we have $\left\vert M\right\vert <m_{p}=p$ for $%
M\neq \left\{ 1,\ldots ,p\right\} $. Consequently, 
\begin{equation*}
c\left( M,\mathcal{M}(m_{p})\right) =\left\vert \mathcal{M}%
(m_{p})\right\vert -2^{\left\vert M\right\vert }\geq
2^{m_{p}}-2^{m_{p}-1}=(1/2)\left\vert \mathcal{M}(m_{p})\right\vert \geq
(1/3)\left\vert \mathcal{M}(m_{p})\right\vert .
\end{equation*}%
Next consider the case where $m_{p}<p$. Then certainly $\left\vert \mathcal{M%
}(m_{p})\right\vert \geq 2^{m_{p}+1}-1$ holds. Now, for $M\in \mathcal{M}%
(m_{p})$ we have $\left\vert M\right\vert \leq m_{p}$ (and $M\neq \left\{
1,\ldots ,p\right\} $) and thus%
\begin{eqnarray*}
c\left( M,\mathcal{M}(m_{p})\right) &=&\left\vert \mathcal{M}%
(m_{p})\right\vert -2^{\left\vert M\right\vert }\geq \left\vert \mathcal{M}%
(m_{p})\right\vert -2^{m_{p}}=\left\vert \mathcal{M}(m_{p})\right\vert
\left( 1-2^{m_{p}}/\left\vert \mathcal{M}(m_{p})\right\vert \right) \\
&\geq &\left\vert \mathcal{M}(m_{p})\right\vert \left(
1-2^{m_{p}}/(2^{m_{p}+1}-1)\right) \geq (1/3)\left\vert \mathcal{M}%
(m_{p})\right\vert .
\end{eqnarray*}%
Thus (iii) with $\tau =1/3$ holds. The next claim now follows from
Proposition \ref{constants_2} since $\left\vert \mathcal{M}%
(m_{p})\right\vert =\sum_{k=0}^{m_{p}}\binom{p}{k}$. The final claim is then
a trivial consequence, since $\left\vert \mathcal{M}(p)\right\vert =2^{p}$
in this case. Note that $\mathsf{X}_{n(p),p}(\mathcal{M}(p))\neq \varnothing 
$ by the assumptions on $n(p)$, implying $n(p)\geq p$. $\blacksquare $

\section{Appendix: Proofs for Section \protect\ref%
{section:distribution:dependent} \label{app:B}}

In the subsequent lemma we assume that $\tilde{\sigma}_{1}$ and $\tilde{%
\sigma}_{2}$ are defined on the same probability space as are $Y$, $X$, and $%
\hat{\sigma}^{2}$. In slight abuse of notation, we shall then denote by $%
P_{n,\beta ,\sigma }$ the joint distribution of $Y$, $X$, $\hat{\sigma}^{2}$%
, $\tilde{\sigma}_{1}$, and $\tilde{\sigma}_{2}$. We note that an argument
corresponding to a special case of this lemma has been used in \cite%
{ewald12influence}.

\begin{lem}
\label{lem:equivalence:beta:betan} Suppose that the maintained model
assumptions of Section \ref{section:distribution:dependent} are satisfied.
Assume further that Conditions \ref{cond:model:selection:procedure} and \ref%
{cond:L} hold. Let $\mathcal{W}$ be the set of all measurable non-negative
functions of the form $W(x_{0},X,M)$. Then, for any two sequences of random
variables $\tilde{\sigma}_{1}=\tilde{\sigma}_{1,n}$ and $\tilde{\sigma}_{2}=%
\tilde{\sigma}_{2,n}$ (which may be functions of $\sigma $) satisfying 
\begin{equation}
\sup_{\beta \in \mathbb{R}^{p},\sigma >0}P_{n,\beta ,\sigma }\left( \left.
\left\vert \left( \tilde{\sigma}_{i}/\sigma \right) -1\right\vert >\delta
\right\vert X\right) \rightarrow 0  \label{eq:var-condition}
\end{equation}%
in probability as $n\rightarrow \infty $ for every $\delta >0$ and for $%
i=1,2 $, we have that%
\begin{eqnarray*}
&&\sup_{x_{0}\in \mathbb{R}^{p},\beta \in \mathbb{R}^{p},\sigma >0,W\in 
\mathcal{W}}\left\vert P_{n,\beta ,\sigma }\left( \left. \left\vert
x_{0}^{\prime }[\hat{M}]\hat{\beta}_{\hat{M}}-x_{0}^{\prime }[\hat{M}]\beta
_{\hat{M}}^{(\star )}\right\vert \leq W(x_{0},X,\hat{M})\tilde{\sigma}%
_{1}\right\vert X\right) \right. \\
&&-\left. P_{n,\beta ,\sigma }\left( \left. \left\vert x_{0}^{\prime }[\hat{M%
}]\hat{\beta}_{\hat{M}}-x_{0}^{\prime }[\hat{M}]\beta _{\hat{M}%
}^{(n)}\right\vert \leq W(x_{0},X,\hat{M})\tilde{\sigma}_{2}\right\vert
X\right) \right\vert
\end{eqnarray*}%
converges to $0$ in probability as $n\rightarrow \infty $.
\end{lem}

\textbf{Proof:} Because the number of variables $p$ is fixed, it suffices to
show for arbitrary but fixed $M\subseteq \{1,\ldots ,p\}$ that%
\begin{eqnarray*}
Q_{n} &=&\sup_{x_{0}\in \mathbb{R}^{p},\beta \in \mathbb{R}^{p},\sigma
>0,W\in \mathcal{W}}\left\vert P_{n,\beta ,\sigma }\left( \left. \left\vert
x_{0}^{\prime }[M]\hat{\beta}_{M}-x_{0}^{\prime }[M]\beta _{M}^{(\star
)}\right\vert \leq W_{M}\tilde{\sigma}_{1};\hat{M}=M\right\vert X\right)
\right. \\
&&-\left. P_{n,\beta ,\sigma }\left( \left. \left\vert x_{0}^{\prime }[M]%
\hat{\beta}_{M}-x_{0}^{\prime }[M]\beta _{M}^{(n)}\right\vert \leq W_{M}%
\tilde{\sigma}_{2};\hat{M}=M\right\vert X\right) \right\vert
\end{eqnarray*}%
goes to $0$ in probability, where we have used the abbreviation $%
W_{M}=W(x_{0},X,M)$. We may assume in what follows that $M\neq \varnothing $
since otherwise $Q_{n}$ is zero. Furthermore, $Q_{n}$ does not change its
value if the supremum is restricted to those $x_{0}$ which have $\left\Vert
x_{0}[M]\right\Vert =1$ (since the expression inside the supremum is
identically zero if $x_{0}$ satisfies $x_{0}[M]=0$ and since otherwise the
norm of $x_{0}[M]$ can be absorbed into $W_{M}$). Hence we have%
\begin{eqnarray}
Q_{n} &=&\sup_{x_{0}\in S(M),\beta \in \mathbb{R}^{p},\sigma >0,W\in 
\mathcal{W}}\left\vert P_{n,\beta ,\sigma }\left( \left. \left\vert \sigma
e_{1}+\sigma e_{2}\right\vert \leq W_{M}\tilde{\sigma}_{1};\hat{M}%
=M\right\vert X\right) \right.  \notag \\
&&-\left. P_{n,\beta ,\sigma }\left( \left. \left\vert \sigma
e_{2}\right\vert \leq W_{M}\tilde{\sigma}_{2};\hat{M}=M\right\vert X\right)
\right\vert  \label{eq:diff_prob}
\end{eqnarray}%
where we have used the abbreviations $S(M)=\left\{ x_{0}\in \mathbb{R}%
^{p}:\left\Vert x_{0}[M]\right\Vert =1\right\} $,%
\begin{equation*}
e_{1}=\sigma ^{-1}n^{1/2}x_{0}^{\prime }[M]\left( \left( X[M]^{\prime
}X[M]\right) ^{-1}X[M]^{\prime }X[M^{c}]-\left( \Sigma \lbrack M,M]\right)
^{-1}\Sigma \lbrack M,M^{c}]\right) \beta \lbrack M^{c}]
\end{equation*}%
and 
\begin{equation*}
e_{2}=\sigma ^{-1}n^{1/2}x_{0}^{\prime }[M]\left( X[M]^{\prime }X[M]\right)
^{-1}X[M]^{\prime }\left( Y-X\beta \right) .
\end{equation*}%
Note that we have also absorbed a factor $n^{1/2}$ into $W_{M}$, which is
possible because of the supremum operation w.r.t. $W_{M}$. Using the
inequality $\left\vert \Pr \left( A\cap C\right) -\Pr \left( B\cap C\right)
\right\vert \leq \Pr \left( A^{c}\cap B\cap C\right) +\Pr \left( A\cap
B^{c}\cap C\right) $ we can bound the absolute value inside the supremum in (%
\ref{eq:diff_prob}) by%
\begin{eqnarray}
&&P_{n,\beta ,\sigma }\left( \left. \left\vert \sigma e_{1}+\sigma
e_{2}\right\vert >W_{M}\tilde{\sigma}_{1};\left\vert \sigma e_{2}\right\vert
\leq W_{M}\tilde{\sigma}_{2};\hat{M}=M\right\vert X\right)  \notag \\
&&+P_{n,\beta ,\sigma }\left( \left. \left\vert \sigma e_{1}+\sigma
e_{2}\right\vert \leq W_{M}\tilde{\sigma}_{1};\left\vert \sigma
e_{2}\right\vert >W_{M}\tilde{\sigma}_{2};\hat{M}=M\right\vert X\right) .
\label{eq:bound_1}
\end{eqnarray}%
Let now $\delta _{n,1}$ be an arbitrary sequence of positive numbers
converging to zero. Then we can further bound the above expression by%
\begin{eqnarray}
&&P_{n,\beta ,\sigma }\left( \left. \tilde{\sigma}_{1}\left( W_{M}-\delta
_{n,1}\right) _{+}\leq \left\vert \sigma e_{2}\right\vert \leq \tilde{\sigma}%
_{2}W_{M}\right\vert X\right)  \notag \\
&&+P_{n,\beta ,\sigma }\left( \left. \tilde{\sigma}_{2}W_{M}\leq \left\vert
\sigma e_{2}\right\vert \leq \tilde{\sigma}_{1}\left( W_{M}+\delta
_{n,1}\right) \right\vert X\right)  \notag \\
&&+2P_{n,\beta ,\sigma }\left( \left. \left\vert \sigma e_{1}\right\vert
\geq \tilde{\sigma}_{1}\delta _{n,1};\hat{M}=M\right\vert X\right) .
\label{eq:bound_2}
\end{eqnarray}%
By the assumption on the estimators $\tilde{\sigma}_{1}$ and $\tilde{\sigma}%
_{2}$ we can find a sequence $\delta _{n,2}<1$ of positive numbers
converging to zero such that%
\begin{equation}
\sup_{\beta ,\sigma }P_{n,\beta ,\sigma }\left( \left.
\max_{i=1,2}\left\vert \left( \tilde{\sigma}_{i}/\sigma \right)
-1\right\vert >\delta _{n,2}\right\vert X\right) \rightarrow 0
\label{eq:sigma}
\end{equation}%
in probability as $n\rightarrow \infty $. This can easily be seen from a
diagonal sequence argument. Now, using (\ref{eq:diff_prob}), (\ref%
{eq:bound_1}), (\ref{eq:bound_2}), and (\ref{eq:sigma}), we have%
\begin{eqnarray*}
Q_{n} &\leq &\sup_{x_{0}\in S(M),\beta \in \mathbb{R}^{p},\sigma >0,W\in 
\mathcal{W}}P_{n,\beta ,\sigma }\left( \left. \left( 1-\delta _{n,2}\right)
\left( W_{M}-\delta _{n,1}\right) _{+}\leq \left\vert e_{2}\right\vert \leq
\left( 1+\delta _{n,2}\right) W_{M}\right\vert X\right) \\
&&+\sup_{x_{0}\in S(M),\beta \in \mathbb{R}^{p},\sigma >0,W\in \mathcal{W}%
}P_{n,\beta ,\sigma }\left( \left. \left( 1-\delta _{n,2}\right) W_{M}\leq
\left\vert e_{2}\right\vert \leq \left( 1+\delta _{n,2}\right) \left(
W_{M}+\delta _{n,1}\right) \right\vert X\right) \\
&&+2\sup_{x_{0}\in S(M),\beta \in \mathbb{R}^{p},\sigma >0,W\in \mathcal{W}%
}P_{n,\beta ,\sigma }\left( \left. \left\vert e_{1}\right\vert \geq \left(
1-\delta _{n,2}\right) \delta _{n,1};\hat{M}=M\right\vert X\right) +o_{p}(1)
\\
&\leq &2\sup_{x_{0}\in S(M),\beta \in \mathbb{R}^{p},\sigma >0,W\in \mathcal{%
W}}P_{n,\beta ,\sigma }\left( \left. \left( 1-\delta _{n,2}\right) \left(
W_{M}-\delta _{n,1}\right) _{+}\leq \left\vert e_{2}\right\vert \leq \left(
1+\delta _{n,2}\right) \left( W_{M}+\delta _{n,1}\right) \right\vert X\right)
\\
&&+2\sup_{x_{0}\in S(M),\beta \in \mathbb{R}^{p},\sigma >0}P_{n,\beta
,\sigma }\left( \left. \left\vert e_{1}\right\vert \geq \left( 1-\delta
_{n,2}\right) \delta _{n,1};\hat{M}=M\right\vert X\right) +o_{p}(1) \\
&=&2Q_{n,1}+2Q_{n,2}+o_{p}(1).
\end{eqnarray*}%
We first bound $Q_{n,1}$ as follows: Observe that, conditionally on $X$, the
quantity $e_{2}$ is normally distributed with mean zero and variance given
by $c_{n}\left( x_{0},X\right) =x_{0}^{\prime }[M]\left( n^{-1}X[M]^{\prime
}X[M]\right) ^{-1}x_{0}[M]$. By Condition \ref{cond:L} the variance $%
c_{n}\left( x_{0},X\right) $ converges to $c(x_{0})=x_{0}^{\prime }[M]\left(
\Sigma \lbrack M,M]\right) ^{-1}x_{0}[M]>0$ in probability, and in fact even
uniformly in $x_{0}\in S(M)$. Since $\Sigma \lbrack M,M]$ is obviously
positive definite, $0<c_{\ast }\leq c(x_{0})\leq c^{\ast }<\infty $ must
hold for all $x_{0}\in S(M)$. Consequently, 
\begin{equation}
\sup_{x_{0}\in S(M)}\left\vert \left( c_{n}^{1/2}\left( x_{0},X\right)
/c^{1/2}\left( x_{0}\right) \right) -1\right\vert  \label{eq:variance}
\end{equation}%
converges to zero in probability. Therefore we can find a sequence $\delta
_{n,3}\in \left( 0,1\right) $ converging to zero for $n\rightarrow \infty $
such that the event $D_{n}$ where (\ref{eq:variance}) is less than $\delta
_{n,3}$ has probability converging to $1$. On this event $\inf_{x_{0}\in
S(M)}c_{n}\left( x_{0},X\right) $ is then positive for sufficiently large $n$
and we have on $D_{n}$ and for sufficiently large $n$%
\begin{eqnarray*}
Q_{n,1} &=&2\sup_{x_{0}\in S(M),W\in \mathcal{W}}\left\{ \Phi \left( \left(
1+\delta _{n,2}\right) \left( W_{M}+\delta _{n,1}\right) /c_{n}^{1/2}\left(
x_{0},X\right) \right) \right. \\
&&\left. -\Phi \left( \left( 1-\delta _{n,2}\right) \left( W_{M}-\delta
_{n,1}\right) _{+}/c_{n}^{1/2}\left( x_{0},X\right) \right) \right\} \\
&\leq &2\sup_{x_{0}\in S(M),W\in \mathcal{W}}\left\{ \Phi \left( \frac{%
1+\delta _{n,2}}{1-\delta _{n,3}}\left( W_{M}+\delta _{n,1}\right)
/c^{1/2}\left( x_{0}\right) \right) \right. \\
&&\left. -\Phi \left( \frac{1-\delta _{n,2}}{1+\delta _{n,3}}\left(
W_{M}-\delta _{n,1}\right) _{+}/c^{1/2}\left( x_{0}\right) \right) \right\}
\\
&\leq &2\sup_{x_{0}\in S(M),W\in \mathcal{W}}\left\{ \Phi \left( \frac{%
1+\delta _{n,2}}{1-\delta _{n,3}}\left( W_{M}+\delta _{n,1}\right)
/c^{1/2}\left( x_{0}\right) \right) \right. \\
&&\left. -\Phi \left( \frac{1-\delta _{n,2}}{1+\delta _{n,3}}\left(
W_{M}-\delta _{n,1}\right) /c^{1/2}\left( x_{0}\right) \right) \right\} \\
&\leq &2\sup_{x_{0}\in S(M),z\geq 0}\left\{ \Phi \left( \frac{1+\delta _{n,2}%
}{1-\delta _{n,3}}\left( z+\delta _{n,1}/c^{1/2}\left( x_{0}\right) \right)
\right) \right. \\
&&\left. -\Phi \left( \frac{1-\delta _{n,2}}{1+\delta _{n,3}}\left( z-\delta
_{n,1}/c^{1/2}\left( x_{0}\right) \right) \right) \right\} \\
&\leq &2\sup_{z\geq 0}\left\{ \Phi \left( \frac{1+\delta _{n,2}}{1-\delta
_{n,3}}\left( z+\delta _{n,1}/c_{\ast }^{1/2}\right) \right) -\Phi \left( 
\frac{1-\delta _{n,2}}{1+\delta _{n,3}}\left( z-\delta _{n,1}/c_{\ast
}^{1/2}\right) \right) \right\} ,
\end{eqnarray*}%
where $\Phi $ denotes the standard normal c.d.f. But the far right-hand side
in the above display obviously converges to zero for $n\rightarrow \infty $
since $\delta _{n,1}$, $\delta _{n,2}$, as well as $\delta _{n,3}$ converge
to zero. We have thus established that $Q_{n,1}$ converges to zero in
probability as $n\rightarrow \infty $.

We next turn to $Q_{n,2}$. In case $M=\left\{ 1,\ldots ,p\right\} $, we have
that $e_{1}=0$, and hence $Q_{n,2}=0$. Otherwise, from Condition \ref%
{cond:model:selection:procedure} we can conclude (from a diagonal sequence
argument) the existence of a sequence of positive numbers $\delta _{n,4}$
that converge to zero for $n\rightarrow \infty $ such that 
\begin{equation*}
\sup \left\{ P_{n,\beta ,\sigma }(\hat{M}=M|X):\beta \in \mathbb{R}%
^{p},\sigma >0,\left\Vert \beta \lbrack M^{c}]\right\Vert /\sigma \geq
\delta _{n,4}\right\} \rightarrow 0
\end{equation*}%
in probability as $n\rightarrow \infty $. Then%
\begin{eqnarray}
Q_{n,2} &\leq &\sup_{x_{0}\in S(M),\left\Vert \beta \lbrack
M^{c}]\right\Vert /\sigma \geq \delta _{n,4}}P_{n,\beta ,\sigma }\left(
\left. \left\vert e_{1}\right\vert \geq \left( 1-\delta _{n,2}\right) \delta
_{n,1};\hat{M}=M\right\vert X\right)  \notag \\
&&+\sup_{x_{0}\in S(M),\left\Vert \beta \lbrack M^{c}]\right\Vert /\sigma
<\delta _{n,4}}P_{n,\beta ,\sigma }\left( \left. \left\vert e_{1}\right\vert
\geq \left( 1-\delta _{n,2}\right) \delta _{n,1};\hat{M}=M\right\vert
X\right)  \notag \\
&\leq &\sup_{\left\Vert \beta \lbrack M^{c}]\right\Vert /\sigma \geq \delta
_{n,4}}P_{n,\beta ,\sigma }\left( \left. \hat{M}=M\right\vert X\right) 
\notag \\
&&+\sup_{x_{0}\in S(M),\left\Vert \beta \lbrack M^{c}]\right\Vert /\sigma
<\delta _{n,4}}P_{n,\beta ,\sigma }\left( \left. \left\vert e_{1}\right\vert
\geq \left( 1-\delta _{n,2}\right) \delta _{n,1}\right\vert X\right)  \notag
\\
&\leq &o_{p}(1)+\sup_{x_{0}\in S(M),\left\Vert \beta \lbrack
M^{c}]\right\Vert /\sigma <\delta _{n,4}}P_{n,\beta ,\sigma }\left( \left.
\left\vert e_{1}\right\vert \geq \left( 1-\delta _{n,2}\right) \delta
_{n,1}\right\vert X\right) .  \label{eq:bound_Qn2}
\end{eqnarray}%
Using the Cauchy-Schwartz inequality we obtain for $x_{0}\in S(M)$%
\begin{eqnarray*}
\left\vert e_{1}\right\vert &\leq &\left\Vert x_{0}^{\prime }[M]\right\Vert
\left\Vert \beta \lbrack M^{c}]/\sigma \right\Vert \left\Vert n^{1/2}\left(
\left( X[M]^{\prime }X[M]\right) ^{-1}X[M]^{\prime }X[M^{c}]-\left( \Sigma
\lbrack M,M]\right) ^{-1}\Sigma \lbrack M,M^{c}]\right) \right\Vert \\
&\leq &\left\Vert \beta \lbrack M^{c}]/\sigma \right\Vert B_{n}\left(
X\right)
\end{eqnarray*}%
where $B_{n}\left( X\right) \geq 0$ is $O_{p}\left( 1\right) $, this
following from Condition \ref{cond:L} and positive definiteness of $\Sigma
\lbrack M,M]$. This shows that the second term on the far right-hand side of
(\ref{eq:bound_Qn2}) is bounded by%
\begin{equation*}
\boldsymbol{1}\left( \delta _{n,4}B_{n}\left( X\right) \geq \left( 1-\delta
_{n,2}\right) \delta _{n,1}\right) .
\end{equation*}%
If we set now, for example, $\delta _{n,1}=\delta _{n,4}^{1/2}$, we see that
the above quantity converges to zero in probability as $n\rightarrow \infty $%
, implying that $Q_{n,2}$ converges to zero in probability as $n\rightarrow
\infty $. This completes the proof. $\blacksquare $

\textbf{Proof of Theorem \ref{theo:asympt:val:CI}:} (a) Use Lemma \ref%
{lem:equivalence:beta:betan} with $W(x_{0},X,M)$ equal to $K_{1}\left(
x_{0},r\right) \left\Vert s_{M}\right\Vert $ ($K_{2}(x_{0}[M],M,r)\left\Vert
s_{M}\right\Vert $, $K_{3}(x_{0}[M],M,r)\left\Vert s_{M}\right\Vert $, $%
K_{4}\left( r\right) \left\Vert s_{M}\right\Vert $, or $K_{5}\left( r\right)
\left\Vert s_{M}\right\Vert $, respectively) and $\tilde{\sigma}_{1}=\tilde{%
\sigma}_{2}=\hat{\sigma}$ and combine this with Proposition \ref%
{prop:PoSI:coverage} (Corollaries \ref{cor:max:of:PoSI:coverage}, \ref%
{cor:PoSI:upper:bound}, respectively). Note that $r=r_{n}\rightarrow \infty $
because of Condition \ref{cond:sigma}, and hence $\hat{\sigma}$ satisfies (%
\ref{eq:var-condition}).

(b) Let $\tilde{\sigma}_{2}$ be a sequence of random variables such that,
conditionally on $X$, $\tilde{\sigma}_{2}^{2}$ is independent of $\hat{\beta}
$ and is distributed as $\sigma ^{2}/r^{\ast }$ times a chi-squared
distributed random variable with $r^{\ast }$ degrees of freedom with the
convention that $\tilde{\sigma}_{2}=\sigma $ in case $r^{\ast }=\infty $.
[Such a sequence exists: Possibly after redefining the relevant random
variables on a sufficiently rich probability space we may find a sequence $%
\left( Z_{i}\right) _{i\in \mathbb{N}}$ of i.i.d. standard Gaussian random
variables that is independent of $Y$ and $X$. Then define $\tilde{\sigma}%
_{2}^{2}=\sigma ^{2}\sum_{i=1}^{r^{\ast }}Z_{i}^{2}/r^{\ast }$ if $r^{\ast
}<\infty $ and set $\tilde{\sigma}_{2}^{2}=\sigma ^{2}$ otherwise.] In view
of Remarks \ref{rem:kmown:variance:case} and \ref{rem_2.6}(iii) we have that
Proposition \ref{prop:PoSI:coverage} (Corollaries \ref%
{cor:max:of:PoSI:coverage}, \ref{cor:PoSI:upper:bound}, respectively) also
hold if the confidence interval (\ref{eq:general:form:CI}) for the target $%
x_{0}^{\prime }{[\hat{M}]}\beta _{\hat{M}}^{(n)}$ uses $\tilde{\sigma}_{2}$
instead of $\hat{\sigma}$ (but uses ${\hat{M}}$ as in part (a)) and uses the
constants $K_{1}\left( x_{0},r^{\ast }\right) $ ($K_{2}(x_{0}[\hat{M}],\hat{M%
},r^{\ast })$, $K_{3}(x_{0}[\hat{M}],\hat{M},r^{\ast })$, $K_{4}\left(
r^{\ast }\right) $, or $K_{5}\left( r^{\ast }\right) $, respectively). Now
apply Lemma \ref{lem:equivalence:beta:betan} with $W(x_{0},X,M)$ equal to $%
K_{1}\left( x_{0},r^{\ast }\right) \left\Vert s_{M}\right\Vert $ ($%
K_{2}(x_{0}[M],M,r^{\ast })\left\Vert s_{M}\right\Vert $, $%
K_{3}(x_{0}[M],M,r^{\ast })\left\Vert s_{M}\right\Vert $, $K_{4}\left(
r^{\ast }\right) \left\Vert s_{M}\right\Vert $, or $K_{5}\left( r^{\ast
}\right) \left\Vert s_{M}\right\Vert $, respectively) and with $\tilde{\sigma%
}_{1}=\tilde{\sigma}$. Note that $\tilde{\sigma}_{1}$ satisfies (\ref%
{eq:var-condition}) by assumption, while $\tilde{\sigma}_{2}$ satisfies it
because $r^{\ast }\rightarrow \infty $ has been assumed. $\blacksquare $

\begin{lem}
\label{lem: PMS-var} Suppose that the maintained model assumptions of
Section \ref{section:distribution:dependent} are satisfied and that $%
X^{\prime }X/n\rightarrow \Sigma $ in probability for $n\rightarrow \infty $%
. Assume further that Condition \ref{cond:model:selection:procedure} holds
and define $\hat{\sigma}_{\hat{M}}^{2}=||Y-X[\hat{M}]\hat{\beta}_{\hat{M}%
}||^{2}/(n-|\hat{M}|)$ for $n>p$. Then $\hat{\sigma}_{\hat{M}}^{2}$
satisfies condition (\ref{eq:unif_cons_var}).
\end{lem}

\textbf{Proof:} Clearly 
\begin{equation*}
Y-X[\hat{M}]\hat{\beta}_{\hat{M}}=Y-P_{X[\hat{M}]}Y=P_{X[\hat{M}]^{\bot
}}U+P_{X[\hat{M}]^{\bot }}X[\hat{M}^{c}]\beta \lbrack \hat{M}^{c}]=A+B,
\end{equation*}%
where $P_{X[\hat{M}]^{\bot }}$ denotes orthogonal projection on the
orthogonal complement of the column space of $X[\hat{M}]$. By the triangle
inequality we hence have%
\begin{eqnarray*}
\left\vert \left( \hat{\sigma}_{\hat{M}}/\sigma \right) -1\right\vert &\leq
&\left\vert \left( n-|\hat{M}|\right) ^{-1/2}\left\Vert A/\sigma \right\Vert
-1\right\vert +\left( n-|\hat{M}|\right) ^{-1/2}\left\Vert B/\sigma
\right\Vert \\
&\leq &\left\vert \left( n-|\hat{M}|\right) ^{-1/2}\left\Vert A/\sigma
\right\Vert -1\right\vert +\left( n-|\hat{M}|\right) ^{-1/2}\left\Vert X[%
\hat{M}^{c}]\beta \lbrack \hat{M}^{c}]/\sigma \right\Vert .
\end{eqnarray*}%
We now bound the probability in (\ref{eq:unif_cons_var}) by the sum of the
probabilities that the first and second term on the far r.h.s. of the
preceding display, respectively, exceed $\delta /2$. Because $p$ is fixed
there is a fixed finite number of possible models $\hat{M}$ and thus for $%
\delta >0$ we have the bound for the first term 
\begin{eqnarray*}
&&\sup_{\beta \in \mathbb{R}^{p},\sigma >0}P_{n,\beta ,\sigma }\left( \left.
\left\vert \left( n-|\hat{M}|\right) ^{-1/2}\left\Vert A/\sigma \right\Vert
-1\right\vert \geq \delta /2\right\vert X\right) \\
&=&\sup_{\beta \in \mathbb{R}^{p},\sigma >0}\sum_{M}P_{n,\beta ,\sigma
}\left( \left. \left\vert \left( n-|M|\right) ^{-1/2}\left\Vert
P_{X[M]^{\bot }}U/\sigma \right\Vert -1\right\vert \geq \delta /2,\hat{M}%
=M\right\vert X\right) \\
&\leq &\sum_{M}\sup_{\beta \in \mathbb{R}^{p},\sigma >0}P_{n,\beta ,\sigma
}\left( \left. \left\vert \left( n-|M|\right) ^{-1/2}\left\Vert
P_{X[M]^{\bot }}U/\sigma \right\Vert -1\right\vert \geq \delta /2\right\vert
X\right) .
\end{eqnarray*}%
Note that the probabilities in the upper bound on the far r.h.s. of the
preceding display do actually neither depend on $\beta $ nor $\sigma $ and
are each of the form $\Pr \left( \left\vert W/w-1\right\vert \geq \delta
\right) $ where $W$ is distributed as the square root of a chi-squared
random variable with $w^{2}$ degrees of freedom. Since $w^{2}=n-\left\vert
M\right\vert $ goes to infinity for $n\rightarrow \infty $ and any fixed $M$%
, and since the sum has a fixed finite number of terms, we can conclude that
the upper bound converges to zero in probability as $n\rightarrow \infty $.

Turning to the second term we have, letting $\lambda _{\max }$ denote the
largest eigenvalue of a symmetric matrix, 
\begin{eqnarray*}
&&\sup_{\beta \in \mathbb{R}^{p},\sigma >0}P_{n,\beta ,\sigma }\left( \left.
\left( n-|\hat{M}|\right) ^{-1/2}\left\Vert X[\hat{M}^{c}]\beta \lbrack \hat{%
M}^{c}]/\sigma \right\Vert \geq \delta /2\right\vert X\right) \\
&\leq &\sup_{\beta \in \mathbb{R}^{p},\sigma >0}P_{n,\beta ,\sigma }\left(
\left. \left( n-|\hat{M}|\right) ^{-1/2}\lambda _{\max }^{1/2}\left( X[\hat{M%
}^{c}]^{\prime }X[\hat{M}^{c}]\right) \left\Vert \beta \lbrack \hat{M}%
^{c}]/\sigma \right\Vert \geq \delta /2\right\vert X\right) \\
&\leq &\sup_{\beta \in \mathbb{R}^{p},\sigma >0}P_{n,\beta ,\sigma }\left(
\left. \lambda _{\max }^{1/2}\left( X^{\prime }X/\left( n-p\right) \right)
\left\Vert \beta \lbrack \hat{M}^{c}]/\sigma \right\Vert \geq \delta
/2\right\vert X\right) \\
&\leq &\sum_{M\neq \left\{ 1,\ldots ,p\right\} }\sup_{\beta \in \mathbb{R}%
^{p},\sigma >0}P_{n,\beta ,\sigma }\left( \hat{M}=M,\left. \lambda _{\max
}^{1/2}\left( X^{\prime }X/\left( n-p\right) \right) \left\Vert \beta
\lbrack M^{c}]/\sigma \right\Vert \geq \delta /2\right\vert X\right) .
\end{eqnarray*}%
Now, since $X^{\prime }X/\left( n-p\right) $ converges to the positive
definite matrix $\Sigma $ in probability, we can find an event $D_{n}$,
which has probability converging to $1$ for $n\rightarrow \infty $, such
that on this event $\lambda _{\max }\left( X^{\prime }X/\left( n-p\right)
\right) $ is not larger than $4\lambda _{\max }\left( \Sigma \right) $.
Hence, on $D_{n}$ we can bound each supremum on the far r.h.s. of the
preceding display by 
\begin{eqnarray*}
&&\sup_{\beta \in \mathbb{R}^{p},\sigma >0}P_{n,\beta ,\sigma }\left( \hat{M}%
=M,\left. \left\Vert \beta \left[ M^{c}\right] /\sigma \right\Vert \geq
\lambda _{\max }^{-1/2}\left( \Sigma \right) \delta /4\right\vert X\right) \\
&=&\sup \left\{ P_{n,\beta ,\sigma }(\hat{M}=M|X):\beta \in \mathbb{R}%
^{p},\sigma >0,\left\Vert \beta \lbrack M^{c}]\right\Vert /\sigma \geq
\lambda _{\max }^{-1/2}\left( \Sigma \right) \delta /4\right\} ,
\end{eqnarray*}%
which goes to zero in probability as $n\rightarrow \infty $ by Condition \ref%
{cond:model:selection:procedure}. $\blacksquare $

\section{Appendix: Comments on and extension of results in Section \protect
\ref{section:distribution:dependent}\label{app_various}}

\subsection{Measurability issues \label{measurability}}

Various statements concerning uncountable suprema (infima) of conditional
probabilities occur in Section \ref{section:distribution:dependent},
Appendix \ref{app:B}, and in Appendix \ref{app ext} below, such as
statements that these quantities converge in probability. It is not
difficult to see that -- in absence of measurability -- all these statements
remain valid if they are properly interpreted as statements referring to
outer probability. This thus relieves one from the need to establish
measurability. For this reason we do not explicitly mention the
measurability issues in the presentation of the results in Section \ref%
{section:distribution:dependent} as well as Appendices \ref{app:B} and \ref%
{app ext}.

\subsection{Some remarks on Theorem \protect\ref{theo:asympt:val:CI}\label%
{rem_on_theorem_3.6}}

\begin{rem}
\normalfont Under the assumptions of Theorem \ref{theo:asympt:val:CI}(b) we
further have that 
\begin{equation*}
\inf_{x_{0}\in \mathbb{R}^{p},\beta \in \mathbb{R}^{p},\sigma >0}P_{n,\beta
,\sigma }\left( \left. x_{0}^{\prime }[\hat{M}]\beta _{\hat{M}}^{(n)}\in
CI^{\ast }(x_{0})\right\vert X\right) \geq (1-\alpha )+o_{p}(1),
\end{equation*}%
holds, where the $o_{p}(1)$ term above depends only on $X$ and converges to
zero in probability as $n\rightarrow \infty $. This follows easily from a
repeated application of Lemma \ref{lem:equivalence:beta:betan} in Appendix %
\ref{app:B}. [Regarding Theorem \ref{theo:asympt:val:CI}(a) recall that the
finite-sample coverage result for the target $x_{0}^{\prime }[\hat{M}]\beta
_{\hat{M}}^{(n)}$ in Section \ref{section:design:dependent} continues to
hold in the context of Section \ref{section:distribution:dependent} if
interpreted conditionally on $X$.]
\end{rem}

\begin{rem}
\label{rem:xzero:random}\normalfont\emph{(Random }$x_{0}$\emph{) }If $x_{0}$
is random and independent of $X$, $U$, and $\hat{\sigma}^{2}$, Theorem \ref%
{theo:asympt:val:CI} continues to hold if the result is then being
interpreted as conditional on $X$ and $x_{0}$. A particular consequence of
this result conditional on $X$ and $x_{0}$ is then that the confidence
interval $CI(x_{0})$ also satisfies%
\begin{equation*}
\inf_{\beta \in \mathbb{R}^{p},\sigma >0}P_{n,\beta ,\sigma }\left( \left.
x_{0}^{\prime }[\hat{M}]\beta _{\hat{M}}^{(\star )}\in CI(x_{0})\right\vert
X\right) \geq (1-\alpha )+o_{p}(1)
\end{equation*}%
where again $o_{p}(1)$ is a function of $X$ only (and $P_{n,\beta ,\sigma }$
here represents the distribution of $Y$, $X$, $\hat{\sigma}^{2}$, and $x_{0}$%
). As noted at the beginning of Section \ref{section:distribution:dependent}%
, the results in Section \ref{section:design:dependent} continue to hold if
interpreted conditionally on $X$ and $x_{0}$. As a consequence, we thus also
have that%
\begin{equation*}
\inf_{\beta \in \mathbb{R}^{p},\sigma >0}P_{n,\beta ,\sigma }\left( \left.
x_{0}^{\prime }[\hat{M}]\beta _{\hat{M}}^{(n)}\in CI(x_{0})\right\vert
X\right) \geq 1-\alpha
\end{equation*}%
holds. See also \cite{Lee09a}, where prediction intervals for $y_{0}$ are
studied in a similar setting.
\end{rem}

\begin{rem}
\normalfont\emph{(Relaxing the assumptions on }$X$\emph{)} The assumption
that the rows of $X$ follow a common distribution $\mathcal{L}$ has been
used only to define the matrix $\Sigma $, which in turn is used in the
definition of $\beta _{M}^{(\star )}$. If this assumption is dropped, but
instead it is assumed that Condition \ref{cond:L} holds for \emph{some}
positive matrix $\Sigma $, which is then used to define $\beta _{M}^{(\star
)}$, Theorem \ref{theo:asympt:val:CI} continues to hold. Note that this
version of Theorem \ref{theo:asympt:val:CI} also covers the case of
nonrandom design matrices for which $n^{-1}X^{\prime }X$ converges to a
positive definite limit at rate $n^{-1/2}$ (or faster).
\end{rem}

\subsection{Extension of Theorem \protect\ref{theo:asympt:val:CI} to the
case $p\rightarrow \infty $\label{app ext}}

We now provide a variant of Theorem \ref{theo:asympt:val:CI} where we
consider the same setup as in Section \ref{section:distribution:dependent}
but now allow $p$ to depend on $n$ such that $p=p(n)\rightarrow \infty $ as $%
n\rightarrow \infty $. The positive definite second moment matrix $\Sigma $
now depends on $n$, i.e., $\Sigma =\Sigma (n)$, and has dimension $%
p(n)\times p(n)$. Note that $p(n)\leq n$ holds. For a symmetric non-negative
definite matrix $A$, we let $A^{1/2}$ denote the unique symmetric
non-negative definite square root of $A$. Furthermore, set $\Sigma
^{(n)}=X^{\prime }X/n$. The theorem below relies on the following two
conditions, the second of which is somewhat intransparent.

\begin{cond}
\label{cond:L:new} We have for $n\rightarrow \infty $ that%
\begin{equation*}
A_{n}(X):=\max_{M\in \mathcal{M},M\neq \varnothing }\left\Vert (\Sigma
\lbrack M,M])^{-1/2}\Sigma ^{(n)}[M,M](\Sigma \lbrack
M,M])^{-1/2}-I_{|M|}\right\Vert =o_{p}(1),
\end{equation*}%
where the norm here is the spectral norm.
\end{cond}

\begin{cond}
\label{cond:model:selection:procedure:new:large:p} The model selection
procedure satisfies for all $\epsilon >0$ that 
\begin{eqnarray*}
&&\sup_{\beta \in \mathbb{R}^{p},\sigma >0,x_{0}\neq 0}P_{n,\beta ,\sigma
}\left( \left\Vert \sqrt{n}(\Sigma \lbrack \hat{M},\hat{M}])^{1/2}\left[
(\Sigma ^{(n)}[\hat{M},\hat{M}])^{-1}\Sigma ^{(n)}[\hat{M},\hat{M}%
^{c}]\right. \right. \right. \\
&&\quad \quad \quad \quad \quad \quad -\left. \left. \left. \left. (\Sigma
\lbrack \hat{M},\hat{M}])^{-1}\Sigma \lbrack \hat{M},\hat{M}^{c}]\right] 
\frac{1}{\sigma }\beta \lbrack \hat{M}^{c}]\right\Vert \geq \epsilon
K_{1}(x_{0},\infty )\right\vert X\right)
\end{eqnarray*}%
converges to $0$ in probability as $n\rightarrow \infty $.
\end{cond}

The confidence interval in the subsequent theorem assumes knowledge of $%
\sigma $ for defining the length of the interval. However, the model
selection procedure $\hat{M}$ may or may not make use of this knowledge.

\begin{theo}
\label{thm:coverage:large:p}Suppose that Conditions \ref{cond:L:new} and \ref%
{cond:model:selection:procedure:new:large:p} hold. Let $\delta >0$ be given.
Denote by $CI^{\delta }(x_{0})$ the interval obtained from (\ref%
{eq:general:form:CI}) by replacing $K(x_{0},\hat{M})$ by $(1+\delta
)K_{1}(x_{0},\infty )$ and $\hat{\sigma}$ by $\sigma $. Then the interval $%
CI^{\delta }(x_{0})$ satisfies 
\begin{equation*}
\inf_{x_{0}\in \mathbb{R}^{p},\beta \in \mathbb{R}^{p},\sigma >0}P_{n,\beta
,\sigma }\left( x_{0}^{\prime }[\hat{M}]\beta _{\hat{M}}^{(\star )}\in
CI^{\delta }(x_{0})|X\right) \geq 1-\alpha +o_{p}(1),
\end{equation*}%
where the $o_{p}(1)$ term above depends only on $X$ (and possibly on the
sequence of second moment matrices $\Sigma =\Sigma (n)$) and converges to
zero in probability as $n\rightarrow \infty $. This result a fortiori holds
if in the definition of $CI^{\delta }(x_{0})$ the constant $%
K_{1}(x_{0},\infty )$ is replaced by any one of the constants $K_{2}(x_{0}[%
\hat{M}],\hat{M},\infty )$, $K_{3}(x_{0}[\hat{M}],\hat{M},\infty )$, $%
K_{4}\left( \infty \right) $, or $K_{5}\left( \infty \right) $, respectively.
\end{theo}

\textbf{Proof of Theorem \ref{thm:coverage:large:p}:} Obviously it suffices
to prove the result for the case where $K_{1}(x_{0},\infty )$ is used. Since 
$x_{0}^{\prime }[\hat{M}]\beta _{\hat{M}}^{(\star )}\in CI^{\delta }(x_{0})$
trivially holds when $x_{0}=0$, it suffices to show that%
\begin{equation*}
\inf_{\beta \in \mathbb{R}^{p},\sigma >0,x_{0}\neq 0}P_{n,\beta ,\sigma
}\left( \left. \left\vert x_{0}^{\prime }[\hat{M}]\beta _{\hat{M}}^{(\star
)}-x_{0}^{\prime }[\hat{M}]\hat{\beta}_{\hat{M}}\right\vert >(1+\delta
)K_{1}(x_{0},\infty )\left\Vert s_{\hat{M}}\right\Vert \sigma \right\vert
X\right) \leq \alpha +o_{p}(1).
\end{equation*}%
Now, we obtain from the triangular inequality, 
\begin{align}
& P_{n,\beta ,\sigma }\left( \left. \left\vert x_{0}^{\prime }[\hat{M}]\beta
_{\hat{M}}^{(\star )}-x_{0}^{\prime }[\hat{M}]\hat{\beta}_{\hat{M}%
}\right\vert >(1+\delta )K_{1}(x_{0},\infty )\left\Vert s_{\hat{M}%
}\right\Vert \sigma \right\vert X\right)  \notag \\
& \leq P_{n,\beta ,\sigma }\left( \left. \left\vert x_{0}^{\prime }[\hat{M}%
]\beta _{\hat{M}}^{(n)}-x_{0}^{\prime }[\hat{M}]\hat{\beta}_{\hat{M}%
}\right\vert >K_{1}(x_{0},\infty )\left\Vert s_{\hat{M}}\right\Vert \sigma
\right\vert X\right)  \notag \\
& +P_{n,\beta ,\sigma }\left( \left. \left\vert x_{0}^{\prime }[\hat{M}%
]\beta _{\hat{M}}^{(\star )}-x_{0}^{\prime }[\hat{M}]\beta _{\hat{M}%
}^{(n)}\right\vert >\delta K_{1}(x_{0},\infty )\left\Vert s_{\hat{M}%
}\right\Vert \sigma \right\vert X\right) .  \label{eq:for:thm:3p6:large:p}
\end{align}%
The first probability on the r.h.s. in the preceding display is not larger
than $\alpha $ as a consequence of Proposition \ref{prop:PoSI:coverage} and
Remarks \ref{rem:kmown:variance:case}, \ref{rem_2.6}(iii). For the second
probability we have 
\begin{eqnarray}
&&P_{n,\beta ,\sigma }\left( \left. \left\vert x_{0}^{\prime }[\hat{M}]\beta
_{\hat{M}}^{(\star )}-x_{0}^{\prime }[\hat{M}]\beta _{\hat{M}%
}^{(n)}\right\vert >\delta K_{1}(x_{0},\infty )\left\Vert s_{\hat{M}%
}\right\Vert \sigma \right\vert X\right)  \notag \\
&=&P_{n,\beta ,\sigma }\left( \left\vert x_{0}^{\prime }[\hat{M}]\left[
(\Sigma ^{(n)}[\hat{M},\hat{M}])^{-1}\Sigma ^{(n)}[\hat{M},\hat{M}%
^{c}]-(\Sigma \lbrack \hat{M},\hat{M}])^{-1}\Sigma \lbrack \hat{M},\hat{M}%
^{c}]\right] \frac{\beta \lbrack \hat{M}^{c}]}{\sigma }\right\vert \right. 
\notag \\
&&\quad \quad \quad \left. >\left. \frac{1}{\sqrt{n}}\left\Vert (\Sigma
^{(n)}[\hat{M},\hat{M}])^{-1/2}x_{0}[\hat{M}]\right\Vert \delta
K_{1}(x_{0},\infty )\right\vert X\right)  \label{event} \\
&\leq &P_{n,\beta ,\sigma }\left( \left\Vert (\Sigma \lbrack \hat{M},\hat{M}%
])^{-1/2}x_{0}[\hat{M}]\right\Vert \left\Vert \sqrt{n}(\Sigma \lbrack \hat{M}%
,\hat{M}])^{1/2}\left[ (\Sigma ^{(n)}[\hat{M},\hat{M}])^{-1}\Sigma ^{(n)}[%
\hat{M},\hat{M}^{c}]\right. \right. \right.  \notag \\
&&\quad \quad \left. \left. -\left. \left. (\Sigma \lbrack \hat{M},\hat{M}%
])^{-1}\Sigma \lbrack \hat{M},\hat{M}^{c}]\right] \frac{\beta \lbrack \hat{M}%
^{c}]}{\sigma }\right\Vert >\left\Vert (\Sigma ^{(n)}[\hat{M},\hat{M}%
])^{-1/2}x_{0}[\hat{M}]\right\Vert \delta K_{1}(x_{0},\infty )\right\vert
X\right) .  \notag
\end{eqnarray}%
With $\lambda _{\min }(.)$ denoting the smallest eigenvalue of a symmetric
matrix, we have 
\begin{align*}
& \left\Vert (\Sigma ^{(n)}[\hat{M},\hat{M}])^{-1/2}x_{0}[\hat{M}%
]\right\Vert ^{2}=\left\Vert (\Sigma ^{(n)}[\hat{M},\hat{M}])^{-1/2}(\Sigma
\lbrack \hat{M},\hat{M}])^{1/2}(\Sigma \lbrack \hat{M},\hat{M}])^{-1/2}x_{0}[%
\hat{M}]\right\Vert ^{2} \\
& \geq \lambda _{\min }\left( (\Sigma \lbrack \hat{M},\hat{M}])^{1/2}(\Sigma
^{(n)}[\hat{M},\hat{M}])^{-1}(\Sigma \lbrack \hat{M},\hat{M}])^{1/2}\right)
\left\Vert (\Sigma \lbrack \hat{M},\hat{M}])^{-1/2}x_{0}[\hat{M}]\right\Vert
^{2}.
\end{align*}%
On the event where $\hat{M}\neq \varnothing $ it holds that%
\begin{equation*}
\lambda _{\min }\left( (\Sigma \lbrack \hat{M},\hat{M}])^{1/2}(\Sigma ^{(n)}[%
\hat{M},\hat{M}])^{-1}(\Sigma \lbrack \hat{M},\hat{M}])^{1/2}\right) \geq
1-A_{n}(X).
\end{equation*}%
This gives%
\begin{equation*}
\left\Vert (\Sigma ^{(n)}[\hat{M},\hat{M}])^{-1/2}x_{0}[\hat{M}]\right\Vert
^{2}\geq (1-A_{n}(X))\left\Vert (\Sigma \lbrack \hat{M},\hat{M}%
])^{-1/2}x_{0}[\hat{M}]\right\Vert ^{2}.
\end{equation*}%
Note that on the event inside the probability on the far r.h.s. of (\ref%
{event}) $\hat{M}\neq \varnothing $ and $x_{0}[\hat{M}]\neq 0$ must hold (in
view of our conventions). Furthermore, $1-A_{n}(X)>1/4$ holds with
probability going to one as $n\rightarrow \infty $ in view of Condition \ref%
{cond:L:new}. Hence, with $o_{p}(1)$ denoting a term that goes to zero in
probability as $n\rightarrow \infty $ and only depends on $X$, (\ref%
{eq:for:thm:3p6:large:p}) is not larger than%
\begin{eqnarray*}
&&\alpha +o_{p}(1)+P_{n,\beta ,\sigma }\left( \left\Vert \sqrt{n}(\Sigma
\lbrack \hat{M},\hat{M}])^{1/2}\left[ (\Sigma ^{(n)}[\hat{M},\hat{M}%
])^{-1}\Sigma ^{(n)}[\hat{M},\hat{M}^{c}]\right. \right. \right. \\
&&\quad \quad \quad \quad \quad \quad \quad \quad \quad -\left. \left.
\left. \left. (\Sigma \lbrack \hat{M},\hat{M}])^{-1}\Sigma \lbrack \hat{M},%
\hat{M}^{c}]\right] \frac{1}{\sigma }\beta \lbrack \hat{M}^{c}]\right\Vert
\geq \delta K_{1}(x_{0},\infty )/2\right\vert X\right)
\end{eqnarray*}%
The result then follows from Condition \ref%
{cond:model:selection:procedure:new:large:p}. $\blacksquare $

\begin{rem}
If $C=C(n)$ is a subset of $\mathbb{R}^{p(n)}\times (0,\infty )$ and only a
weaker version of Condition \ref{cond:model:selection:procedure:new:large:p}
is assumed to hold, where the supremum in this condition only extends over $%
(\beta ,\sigma )\in C$ and $x_{0}\neq 0$, then the same proof delivers a
version of Theorem \ref{thm:coverage:large:p} where the infimum in the
display in that theorem now extends over $(\beta ,\sigma )\in C$ and all $%
x_{0}$.
\end{rem}

\subsection{A comment on notation\label{notation}}

In Section \ref{section:design:dependent} the probability measure $P_{n,\mu
,\sigma }$ denotes the joint distribution of $Y$ and $\hat{\sigma}^{2}$. In
the case where $\hat{\sigma}^{2}$ depends on extraneous data besides $X$ and 
$Y$, the joint distribution of $Y$ and $\hat{\sigma}^{2}$ can in principle
depend not only on $\mu $ and $\sigma $, but also on additional parameters
governing the joint distribution of $Y$ and the extraneous data. [Of course,
this does not occur if $\hat{\sigma}^{2}$ depends only on $X$ and $Y$.] Such
a dependence on additional parameters is, however, not relevant in the
context of Section \ref{section:design:dependent} and thus is not shown in
the notation, since all the results in that section are based on (\ref%
{eq:def:posi:xzero}) and since the probability in (\ref{eq:def:posi:xzero})
is free from any such additional parameters, the reason being that the event
inside this probability can be expressed in terms of $(P_{X}Y,\hat{\sigma})$
only, with $P_{X}Y$ and $\hat{\sigma}$ being independent by assumption and
their marginal distributions only depending on $\mu $ and $\sigma $ (and the
degree of freedom parameter $r$).

In the context of Section \ref{section:distribution:dependent} various
probability statements such as (\ref{eq:asymptotically:valid:CI}) may in
general depend on additional parameters governing the joint distribution of $%
(X,Y)$ and $\hat{\sigma}^{2}$, which, however, are not shown in the
notation. The assumptions and results in this section (e.g., Theorem \ref%
{theo:asympt:val:CI}) are hence to be read as holding for \emph{given}
values (or a \emph{given} sequence of values depending on $n$) of such
additional parameters. As a consequence, the results in Section \ref%
{section:distribution:dependent} do not provide uniformity guarantees w.r.t.
the additional parameters in general. Of course, if all the conditions can
be checked to hold for all possible sequences of values of the additional
parameters in a certain domain, uniformity w.r.t. the additional parameters
is obtained.

\section{Appendix: Algorithms for computing the confidence intervals \label%
{section:practical:algorithms}}

In this appendix we consider the setting of Section \ref%
{section:design:dependent}. In particular, recall that $X$ is a fixed $%
n\times p$ matrix of rank $d\geq 1$. Let $Q$ be a $n\times d$ matrix so that
the columns of $Q$ form an orthonormal basis of the column space of $X$.
Following \cite{berk13valid} we define $\tilde{Y}=Q^{\prime }Y$ and $\tilde{X%
}=Q^{\prime }X$, the so-called canonical coordinates of $Y$ and $X$, cf.
Section 5.1 in \cite{berk13valid}. We then have $\tilde{Y}=\tilde{\mu}+%
\tilde{U}$ with $\tilde{\mu}=Q^{\prime }\mu $ and $\tilde{U}=Q^{\prime
}U\sim N\left( 0,\sigma ^{2}I_{d}\right) $. It is now easy to see that $P_{%
\tilde{X}}\tilde{Y}=Q^{\prime }P_{X}Y$ and $QP_{\tilde{X}}\tilde{Y}=P_{X}Y$
hold. In particular, the independence of the \emph{given} $\hat{\sigma}^{2}$
from the projection of the data vector on the space spanned by the regressor
holds whether we work with the original data or with the data in canonical
coordinates. Furthermore, setting $\tilde{s}_{M}^{\prime }=x_{0}^{\prime }{%
[M]}(\tilde{X}{[M]}^{\prime }{\tilde{X}[M])}^{-1}\tilde{X}{[M]}^{\prime }$
for $\varnothing \neq M\in \mathcal{M}$ with $\mathcal{M}$ as in Section \ref%
{section:design:dependent} and $\tilde{s}_{M}^{\prime }=0\in \mathbb{R}^{d}$
for $M=\varnothing $, it follows that $\left\Vert s_{M}^{\prime }\right\Vert
=\left\Vert \tilde{s}_{M}^{\prime }\right\Vert $ and $s_{M}^{\prime }(Y-\mu
)=\tilde{s}_{M}^{\prime }(\tilde{Y}-\tilde{\mu})$. For later use define $%
\overset{\_}{\tilde{s}}_{M}=\tilde{s}_{M}/\left\Vert \tilde{s}%
_{M}\right\Vert $ if $\left\Vert \tilde{s}_{M}\right\Vert \neq 0$ and define 
$\overset{\_}{\tilde{s}}_{M}=0$ if $\tilde{s}_{M}=0$. Inspection of (\ref%
{eq:def:posi:xzero}) and of the definition of $F_{M,x_{0}}$ now shows that
all the constants $K_{i}$ remain the same whether they are computed from the
original problem using the design matrix $X$ or from the transformed problem
using the canonical coordinates $\tilde{X}$ (but using the originally given $%
\hat{\sigma}^{2}$ in both cases). Hence, in the algorithms below we shall
work with the canonical coordinates as this facilitates computation. Note
that $x_{0}$ is unaffected by this transformation. In the important case $%
d=p $($\leq n$) the matrices $Q$ and $\tilde{X}$ can be obtained, for
example, from a SVD or a QR decomposition of $X$, cf. Section 5.1 in \cite%
{berk13valid}. [In case $p\geq n=d$, one can always set $Q=I_{n}$.]

The following algorithm for computing $K_{1}(x_{0})$ is similar to that of 
\cite{berk13validold} for computing the PoSI constant. We present it here
for completeness. From Proposition \ref{prop:PoSI:upper:bound} and from the
arguments used to prove (\ref{expression}) in Appendix \ref{app A} we see
that, in case $x_{0}\neq 0$, $K_{1}(x_{0})$ is the solution to%
\begin{equation*}
\mathbb{E}_{G}\Pr \left( \left. \max_{M\in \mathcal{M}}\left\vert \overset{\_%
}{\tilde{s}}_{M}^{\prime }V\right\vert \leq t/G\right\vert G\right)
=1-\alpha ,
\end{equation*}%
where $V$ here is uniformly distributed on the unit sphere of $\mathbb{R}%
^{d} $, independently of $G$ (and $G$ is as in Section \ref%
{section:design:dependent}). The algorithm now replaces the (conditional)
probability in the preceding display by a Monte-Carlo estimator,
analytically performs the integration w.r.t. $G$, and then numerically
solves the resulting equation. We note that in this and the other algorithms
to follow there is no need for Monte-Carlo integration w.r.t. $G$. We shall
denote by $F_{d,r}^{\sharp }$ the c.d.f. of $G$; note that then $%
F_{d,r}^{\sharp }\left( t\right) =F_{d,r}\left( t^{2}/d\right) $, where $%
F_{d,r}$ denotes the c.d.f. of an $F$-distribution with $(d,r)$-degrees of
freedom.

\begin{alg}
\label{alg:Kun} In case $x_{0}\neq 0$, choose $I\in \mathbb{N}$ and generate
independent identically distributed random vectors $V_{1},\ldots ,V_{I}$,
where each $V_{i}$ is uniformly distributed on the unit sphere in $\mathbb{R}%
^{d}$. Calculate the quantities $c_{i}=\max_{M\in \mathcal{M}}\left\vert 
\overset{\_}{\tilde{s}}_{M}^{\prime }V_{i}\right\vert $ with $\overset{\_}{%
\tilde{s}}_{M}$ as defined above. A numerical approximation to $K_{1}(x_{0})$
is then obtained by searching for that value of $K$ that solves%
\begin{equation}
\frac{1}{I}\sum_{i=1}^{I}F_{d,r}^{\sharp }\left( \frac{K}{c_{i}}\right)
=1-\alpha .  \label{eq:K}
\end{equation}%
In case $x_{0}=0$, set $K_{1}(x_{0})=0$.
\end{alg}

Note that for $x_{0}\neq 0$ at least one of the vectors $\overset{\_}{\tilde{%
s}}_{M}$, $M\in \mathcal{M}$, is non-zero, implying that the quantities $%
c_{i}$ are all non-zero with probability $1$; hence the terms $%
F_{d,r}^{\sharp }\left( \frac{K}{c_{i}}\right) $ are well-defined with
probability one. It is now obvious that -- on the event where all $c_{i}$
are non-zero -- the solution $K$ of (\ref{eq:K}) exists, is unique and
positive. The costly factor in Algorithm \ref{alg:Kun} is the maximization
involved in the computation of the quantities $c_{i}$, while searching for
the value of $K$ \ that solves (\ref{eq:K}), for example by bisection
searches, incurs only negligible cost. In our simulations, computing $%
K_{1}(x_{0})$ for $p=d=10$ (with $\mathcal{M}$ the power set of $\{1,...,p\}$%
) and $I=10,000$ takes around one second on a personal computer; and around
10 minutes for $p=d=20$ and $I=1,000$. In case $\mathcal{M}$ is the power
set of $\{1,...,p\}$, the complexity of Algorithm \ref{alg:Kun} will be
exponential in $p$ and thus will be feasible only for moderately large
values of $p$. In relation to this we mention that \cite{berk13validold}
found their algorithm (which is similar to Algorithm \ref{alg:Kun} as noted
above) to be tractable for up to about $p=20$ and $I=1,000$, in which case
the elapsed time was around one hour on 2012 desktop computer equipment.
[The longer running time is due to the fact that \cite{berk13validold} have
to search over $p2^{p-1}$ unit vectors, while we have to search only over $%
2^{p}$ unit vectors.]

The algorithm for computing $K_{3}(x_{0}[M],M)$ is given next. We provide
this algorithm only for non-empty $M\neq \{1,...,p\}$ since in case $%
M=\{1,...,p\}$ we have $K_{3}(x_{0}[M],M)=K_{1}(x_{0})$, which can be
computed by Algorithm \ref{alg:Kun}, and in case $M$ is empty we have $%
K_{3}(x_{0}[M],M)=K_{4}$, which can be computed by Algorithm \ref%
{alg:Kquatre} given below. We now search for the solution of the equation%
\begin{equation*}
1-\alpha =\mathbb{E}_{G}\breve{F}_{M,x_{0}}^{\ast }\left( t/G\right)
\end{equation*}%
where $\breve{F}_{M,x_{0}}^{\ast }$ is a Monte-Carlo estimator of $%
F_{M,x_{0}}^{\ast }$ obtained by replacing the probability involving $V$ by
an empirical Monte-Carlo estimator (and where $\overset{\_}{\tilde{s}}_{M}$%
instead of $\bar{s}_{M}$ is being used). Observing that we need only to
integrate over the range where $\breve{F}_{M,x_{0}}^{\ast }$ is positive
(i.e., where $t/G>m_{\ast }$ defined below), the integrand can be additively
decomposed into a `jump' part and a continuous part. The integral over the
jump part can be expressed analytically in terms of the c.d.f. $%
F_{d,r}^{\sharp }$, whereas the integral over the continuous part is
approximated by an integral over a step function, which again can be
expressed in terms of the c.d.f. $F_{d,r}^{\sharp }$. Recall that $c\left( M,%
\mathcal{M}\right) $ has been defined subsequent to (\ref{eq:F}).

\begin{alg}
\label{alg:Ktrois}Suppose that $M\in \mathcal{M}$ satisfies $\varnothing
\neq M\neq \{1,...,p\}$. Choose $I\in \mathbb{N}$, generate independent
identically distributed random vectors $V_{1},\ldots ,V_{I}$, where each $%
V_{i}$ is uniformly distributed on the unit sphere in $\mathbb{R}^{d}$, and
calculate the quantities $c_{i}=\max_{M_{\ast }\in \mathcal{M},M_{\ast
}\subseteq M}\left\vert \overset{\_}{\tilde{s}}_{M_{\ast }}^{\prime
}V_{i}\right\vert $ with $\overset{\_}{\tilde{s}}_{M_{\ast }}$ as defined
above. In case $d>1$, find $m_{\ast }$ as the smallest value such that%
\begin{equation*}
\frac{1}{I}\sum_{i=1}^{I}\boldsymbol{1}\left( c_{i}>t\right) +c\left( M,%
\mathcal{M}\right) \left( 1-F_{Beta,1/2,(d-1)/2}\left( t^{2}\right) \right)
<1
\end{equation*}%
holds for all $t>m_{\ast }$. Next, choose $J\in \mathbb{N}$, $J>1$, and find
the values $m_{1},...,m_{J-1}$ so that, for $j=1,...,J-1$%
\begin{equation}
\left( 1-F_{Beta,1/2,(d-1)/2}\left( m_{\ast }^{2}\right) \right) \frac{j}{J}%
=\left( 1-F_{Beta,1/2,(d-1)/2}\left( m_{j}^{2}\right) \right)
\label{eq:grid}
\end{equation}%
holds. Set $m_{J}=m_{\ast }$. A numerical approximation to $%
K_{3}(x_{0}[M],M) $ is then obtained by searching for that value of $K$ that
solves%
\begin{align}
1-\alpha & =F_{d,r}^{\sharp }\left( \frac{K}{m_{J}}\right) -\frac{1}{I}%
\sum_{i:c_{i}>m_{J}}\left( F_{d,r}^{\sharp }\left( \frac{K}{m_{J}}\right)
-F_{d,r}^{\sharp }\left( \frac{K}{c_{i}}\right) \right)  \notag \\
& +c\left( M,\mathcal{M}\right) \left( 1-F_{Beta,1/2,(d-1)/2}\left(
m_{J}^{2}\right) \right) \frac{1}{J}\sum_{j=1}^{J-1}\left( F_{d,r}^{\sharp
}\left( \frac{K}{m_{j}}\right) -F_{d,r}^{\sharp }\left( \frac{K}{m_{J}}%
\right) \right) .  \label{eq:KK}
\end{align}%
In case $d=1$, $K_{3}(x_{0}[M],M)$ is the (uniquely determined and positive)
constant $K$ that solves 
\begin{equation*}
1-\alpha =F_{1,r}^{\sharp }\left( K\right) \text{.}
\end{equation*}
\end{alg}

Note that $m_{\ast }$ exists, is uniquely determined, is always positive,
and satisfies $m_{\ast }\leq 1$. [In fact, $m_{\ast }<1$ holds, except in
case $c_{i}=1$ for all $i$, which is a probability zero event.] Provided $%
m_{\ast }<1$ holds, the values $m_{j}$ for $j\geq 1$ are uniquely defined
and satisfy $m_{\ast }<m_{J-1}<\ldots <m_{1}<1$. [In case $m_{\ast }=1$,
then any $m_{j}\geq 1$ would solve (\ref{eq:grid}). But in this case the
r.h.s. of (\ref{eq:KK}) reduces to $F_{d,r}^{\sharp }\left( K\right) $
anyway and hence there is no need for solving equation (\ref{eq:grid}).]
Furthermore, note that the r.h.s. of (\ref{eq:KK}) can be written as%
\begin{eqnarray*}
&&F_{d,r}^{\sharp }\left( \frac{K}{m_{J}}\right) \left[ 1-\frac{1}{I}%
\sum_{i=1}^{I}\boldsymbol{1}\left( c_{i}>m_{J}\right) -c\left( M,\mathcal{M}%
\right) \left( 1-F_{Beta,1/2,(d-1)/2}\left( m_{J}^{2}\right) \right) \frac{%
J-1}{J}\right] \\
&&+\frac{1}{I}\sum_{i:c_{i}>m_{J}}F_{d,r}^{\sharp }\left( \frac{K}{c_{i}}%
\right) +c\left( M,\mathcal{M}\right) \left( 1-F_{Beta,1/2,(d-1)/2}\left(
m_{J}^{2}\right) \right) \frac{1}{J}\sum_{j=1}^{J-1}F_{d,r}^{\sharp }\left( 
\frac{K}{m_{j}}\right) .
\end{eqnarray*}%
Observing that the expression in brackets is nonnegative (in fact, positive)
because of the definition of $m_{J}$, we see that the r.h.s. of (\ref{eq:KK}%
) is strictly increasing in $K$. Furthermore, inspection of the r.h.s. of (%
\ref{eq:KK}) shows that it is zero for $K=0$ and converges to one for $%
K\rightarrow \infty $. Consequently, equation (\ref{eq:KK}) has a unique
solution for $K$, which necessarily is positive. We note that in Algorithm %
\ref{alg:Ktrois} the cost of searching for $m_{\ast }$, for the $m_{j}$'s,
and for $K$, for example by bisection searches,\ is negligible compared to
that of computing the quantities $c_{i}$, which is again the limiting factor.

The above algorithm is based on approximating $1-F_{Beta,1/2,(d-1)/2}\left(
t^{2}\right) $ for $t>m_{\ast }$ by a step function from below. If we
approximate by a step function from above, this results in the same
algorithm except that now the second sum on the r.h.s. of equation (\ref%
{eq:KK}) runs from $j=0$ to $j=J-1$ with the convention that $m_{0}=1$. A
similar argument as above shows that the solution to this modification of (%
\ref{eq:KK}) exists, is unique and is positive. Note that the solutions
obtained from running both versions of the algorithm in parallel provide a
lower as well as an upper bound for the solution one would obtain if the
integration of the continuous part could be performed without error. These
lower and upper bounds allow one to gauge whether or not $J$ has been chosen
large enough such that the effect of the numerical integration error on $K$
is negligible. Note that running the two versions of the algorithm in
parallel is not much more costly than running just one version, as only
(bisection) searches are involved once the $c_{i}$'s have been computed.

The following algorithm for computing $K_{4}$ is similar to the algorithm in 
\cite{berk13validold}, Section 7.2, for computing the universal upper-bound
for the PoSI constants. The computational cost of this algorithm is
negligible compared to those of Algorithms \ref{alg:Kun} and \ref{alg:Ktrois}%
.

\begin{alg}
\label{alg:Kquatre} In case $d>1$, choose $J\in \mathbb{N}$, $J>1$, and find
the values $m_{1},...,m_{J}$ so that, for $j=1,...,J$, 
\begin{equation}
c\left( \varnothing ,\mathcal{M}\right) \left( 1-F_{Beta,1/2,(d-1)/2}\left(
m_{j}^{2}\right) \right) =\frac{j}{J}.  \label{eq:quantile:beta}
\end{equation}%
Then, $K_{4}$ is numerically approximated by the (uniquely determined and
positive) constant $K$ that solves%
\begin{equation}
\frac{1}{J}\sum_{m_{j}>0}F_{d,r}^{\sharp }\left( \frac{K}{m_{j}}\right)
=1-\alpha .  \label{eq:KKK}
\end{equation}%
In case $d=1$, $K_{4}$ is the (uniquely determined and positive) constant $K$
that solves%
\begin{equation}
F_{1,r}^{\sharp }\left( K\right) =1-\alpha .  \label{eq:KKKK}
\end{equation}
\end{alg}

Note that in case $d>1$ the constants $m_{j}$ always exist and are unique;
they are all positive in case $\mathcal{M}\neq \left\{ \varnothing ,\left\{
1,\ldots ,p\right\} \right\} $ (as then $c\left( \varnothing ,\mathcal{M}%
\right) >1$ must hold in view of our assumptions on $\mathcal{M}$), and they
are positive for $j=1,\ldots ,J-1$ in case $\mathcal{M}=\left\{ \varnothing
,\left\{ 1,\ldots ,p\right\} \right\} $. Consequently, the solution $K$ of (%
\ref{eq:KKK}) exists, is unique and positive. In case $d=1$ the solution of (%
\ref{eq:KKKK}) also exists, is unique and positive. As before, this
algorithm relies on approximation by a step function from below. A version
of the algorithm that uses a step function that approximates from above is
obtained if equation (\ref{eq:KKK}) is replaced by 
\begin{equation*}
\frac{1}{J}\sum_{j=0}^{J-1}F_{d,r}^{\sharp }\left( \frac{K}{m_{j}}\right)
=1-\alpha
\end{equation*}%
with the convention that $m_{0}=1$.

\begin{rem}
\label{rem:known:variance:case:algos}\normalfont For the computation of the
constants $K_{1}(x_{0},\infty )$, $K_{3}(x_{0}[M],M,\infty )$, and $%
K_{4}\left( \infty \right) $ (cf. Remark \ref{rem:kmown:variance:case}) one
can use the above algorithms with the only modification that the
distribution function $F_{d,r}^{\sharp }$ is replaced by the distribution
function of the square root of a chi-squared-distributed random variable
with $d$ degrees of freedom.
\end{rem}

\begin{rem}
\label{stability}\normalfont When the collection $\mathcal{M}$ becomes large
(e.g., if $\mathcal{M}$ is the power set of $\{1,...,p\}$ in case $p=d\leq n$
and $d$ is larger than $20$), Algorithms \ref{alg:Kun} or \ref{alg:Ktrois}
may not be tractable, but Algorithm \ref{alg:Kquatre} can still be as it
does not require the costly step of searching over the model universe $%
\mathcal{M}$. However, it is reported in \cite{berk13validold} that, for
about $d\geq 40$, it can be problematic to compute the extreme quantiles in (%
\ref{eq:quantile:beta}) with standard routines. In this case, one can of
course always use the Scheff\'{e} constant $K_{5}$. In practice, one may
also consider in such cases (since $p$ is large) to use rule-of-thumb
constants smaller than $K_{5}$ that are based on asymptotic considerations
such as Corollary \ref{constants_3}: For example, if $p=d\leq n$, $\mathcal{M%
}$ is the power set of $\{1,...,p\}$, but $p$ is very large, this corollary
could be read as suggesting to use the constant $K_{6}=0.866K_{5}$ in (\ref%
{eq:general:form:CI}). A similar advice is given in the framework of \cite%
{berk13validold}. [In case $p>n$ and $\mathcal{M}$ is as in Corollary \ref%
{constants_3}, this corollary can be used to provide appropriate substitutes
for $K_{6}$.] However, we would like to issue a warning here: The asymptotic
results for $p\rightarrow \infty $ like Corollary \ref{constants_3} and the
related results in \cite{berk13validold} and \cite{berk13valid} are highly
non-uniform w.r.t. $\alpha $ (cf. Remark \ref{nonuniformity}), showing that
rule-of-thumb approximations such as $K_{6}$ have to be taken with a grain
of salt; see also the warning expressed at the end of Section 5.2 of \cite%
{berk13validold}.
\end{rem}

\section{Appendix: Details for computations in Section \protect\ref%
{section:simulation:study}\label{app D}}

\subsection{Description of the variables in the watershed data set}

The explanatory variables in $X_{Raw}$ are a constant term (to include an
intercept in the model), rainfall (inches), area of watershed (square
miles), area impervious to water (square miles), average slope of watershed
(percent), longest stream flow in watershed (thousands of feet), surface
absorbency index (0= complete absorbency; 100 = no absorbency), estimated
soil storage capacity (inches of water), infiltration rate of water into
soil (inches/hour) and time period during which rainfall exceeded 1/4
inch/hour. Logarithms are taken of the explanatory variables except for the
intercept. [In \cite{rawlings98applied}, the response corresponding to these
explanatory variables is peak flow rate from watersheds.]

\subsection{Three-step Monte Carlo procedure for determining $K_{2}$ in
Section \protect\ref{length_numerical}}

For the given $X$, $x_{0}$, and $M$ we proceed as follows: First, we
randomly sample $100,000$ independent vectors $x\in \mathbb{R}^{10}$, so
that $x[M]=x_{0}[M]$ and $x[M^{c}]$ follows a Gaussian distribution with
mean vector $0\in \mathbb{R}^{10-|M|}$ and covariance matrix $%
(1/n)(X[M^{c}]^{\prime }X[M^{c}])$. For each of these vectors, we evaluate $%
K_{1}(x)$ with Algorithm \ref{alg:Kun} in Appendix \ref%
{section:practical:algorithms}, with $I_{1}=1,000$ Monte Carlo samples. In
the second step, we keep the $1,000$ vectors $x$ corresponding to the
largest evaluations of $K_{1}(x)$ and we reevaluate $K_{1}(x)$ for them,
with a number of Monte Carlo samples equal to $I_{2}=100,000$ in Algorithm %
\ref{alg:Kun}. In the third step, we keep the vector $x$ from the second
step corresponding to the largest value of $K_{1}$ and we reevaluate $%
K_{1}(x)$ for this $x$, but this time with a number of Monte Carlo samples
equal to $I_{3}=1,000,000$ in Algorithm \ref{alg:Kun}.

\subsection{Details for computing AIC, BIC, LASSO, MCP, and SCAD in Section 
\protect\ref{minimal_cov_prob_numerical}}

For the AIC- and BIC-procedures we use the \verb|step()| function in \verb|R|%
, with penalty parameter \verb|k| equal to $2$ for AIC and $\log (n)$ for
BIC. The AIC and BIC objective functions are minimized through a greedy
general-to-specific search over the resulting $2^{p-1}$ candidate models
(recall that the intercept is protected).

For the LASSO, the selected model corresponds to the explanatory variables
for which the LASSO estimator has non-zero coefficients. More precisely, we
use the \verb|lars| package in \verb|R| and follow suggestions outlined in 
\cite{efron04least}: To protect the first regressor, we first compute the
residual of the orthogonal projection of $Y$ on the first regressor; write $%
\tilde{Y}$ for this residual vector, and write $\tilde{X}$ for the design
matrix $X$ with the first column removed. We then compute the
LASSO-estimator for a regression of $\tilde{Y}$ on $\tilde{X}$ using the 
\verb|lars()| function; the LASSO-penalty is chosen by $10$-fold
cross-validation using the \verb|cv.lars()| function. In both functions we
set the \verb|intercept| parameter to \verb|FALSE|, but otherwise use the
default settings. The selected model is comprised of those regressors in $%
\tilde{X}$ for which the corresponding LASSO coefficients are non-zero, plus
the first column of X.

For SCAD and MCP, we use the \verb|ncvreg| package in \verb|R|. With the
function \verb|cv.ncvreg()| (with parameters \verb|SCAD| or \verb|MCP|) the
penalty is selected by $10$-fold cross-validation and the corresponding
estimated regression coefficients are computed. Like for the LASSO, the
function \verb|cv.ncvreg()| is applied for a regression of $\tilde{Y}$ on $%
\tilde{X}$, and the selected model is comprised of those regressors in $%
\tilde{X}$ for which the SCAD (or MCP) coefficients are non-zero, plus the
first column of X.

\subsection{On the exchangeable and the equicorrelated data set}

The matrix $\tilde{\Sigma}$ in the exchangeable case is related to the
design matrix $\mathbf{X}^{(\tilde{p})}(a)$ defined in Section 6.1 of \cite%
{berk13valid} via $\tilde{\Sigma}=(\mathbf{X}^{(\tilde{p})}(a))^{\prime }(%
\mathbf{X}^{(\tilde{p})}(a))$, where $\tilde{p}=9$ and $a=10$. In the
equicorrelated case $\tilde{\Sigma}$ is related to the design matrix $%
\mathbf{X}^{(\tilde{p})}(c)$ defined in Section 6.2 of \cite{berk13valid}
again via $\tilde{\Sigma}=(\mathbf{X}^{(\tilde{p})}(c))^{\prime }(\mathbf{X}%
^{(\tilde{p})}(c))$, where $\tilde{p}=9$ and $c=\sqrt{0.8/(\tilde{p}-1)}$.

\subsection{Three-step Monte Carlo procedure for estimating minimal coverage
probabilities in Section \protect\ref{minimal_cov_prob_numerical}}

For each configuration of $n$, $\Sigma $, the model selection procedure, the
target (either the design-dependent or the design-independent target), as
well as of a matrix $X$ and a vector $x_{0}$ as described in Section \ref%
{minimal_cov_prob_numerical} and for each of the constants $K_{naive}$, $%
K_{1}$, $K_{3}$, and $K_{4}$ we estimate the minimal (over $\beta $ and $%
\sigma $) coverage probabilities (conditional on $X$ and $x_{0}$) of the
confidence intervals by a three-step Monte Carlo procedure as follows: We
first sample independently $m_{1}=1,000$ parameters $\beta $ from a $p$%
-dimensional random vector $b$ where $Xb$ follows a standard Gaussian
distribution within the column-space of $X$. Then, for each of these vectors 
$\beta $, we draw $I_{1}=1000$ Monte Carlo samples from the full model
(i.e., from a $N(X\beta ,\sigma ^{2}I_{n})$-distribution) using $\beta $ and 
$\sigma =1$ as the true parameters. [For invariance reasons it suffices to
consider only the case where $\sigma =1$.] For each Monte Carlo sample, we
use the standard unbiased estimator $\hat{\sigma}^{2}$ of the error variance
(under the full linear model), we carry out the model-selection procedure $%
\hat{M}$, and we record whether or not the target currently under
investigation is covered by the confidence interval obtained from (\ref%
{eq:general:form:CI}) with $K(x_{0},\hat{M})$ replaced by the constant $K$
under investigation. [For $K=K_{i}$, $i=1,3,4$, the value of $K_{i}$ is
obtained from the algorithms described in Appendix \ref%
{section:practical:algorithms}.] For each $\beta $, the $I_{1}$ recorded
results are then averaged, resulting in $m_{1}$ Monte Carlo estimates of the
coverage probabilities depending on the $m_{1}$ sampled vectors $\beta $.
Then for the $m_{2}=100$ vectors $\beta $ corresponding to the smallest
estimated coverage probabilities from the first step, we repeat the Monte
Carlo procedures, but this time with $I_{2}=10,000$ Monte Carlo samples, and
we record the vector $\beta $ that yields the smallest estimate for the
coverage probability in this second step. Performing these two steps for
each of the four constants $K_{naive}$, $K_{1}$, $K_{3}$, and $K_{4}$
results in four vectors $\beta (1)$, $\beta (2)$, $\beta (3)$, and $\beta
(4) $. In a third step, we now reevaluate the coverage probability of any of
the four confidence intervals at each of the vectors $\beta (j)$, $%
j=1,\ldots ,4$, this time now with $I_{3}=100,000$ Monte Carlo samples, and
record, for each of the confidence intervals, the minimum of these four
estimates of the coverage probabilities. This is then used as the final
estimate of the minimal coverage probability of the confidence interval
under consideration.

\subsection{Three-step Monte Carlo procedure for estimating minimal
conditional coverage probabilities in Section \protect\ref{taylor}}

For each of the eight configurations mentioned in Section \ref{taylor}, we
carry out the three-step minimal coverage probability evaluation described
just above, with the same values of $m_{1}$, $m_{2}$ and $I_{1}$, $I_{2}$, $%
I_{3}$ (with the only difference that in the third step the (now
conditional) coverage probability is reevaluated only for one value $\beta
(1)$, say, where $\beta (1)$ corresponds to that value of $\beta $ that
gives the smallest estimate for the coverage probability in the second
step). When we evaluate a conditional coverage probability for a given $%
\beta $ in this process, we proceed as follows: We sample $I_{1}$ (or $I_{2}$
or $I_{3}$) values of $Y$ from the $N(X\beta ,I_{n})$ distribution. For each
value of $Y$ we run the model selection procedure $\hat{M}$, where $\lambda =%
\hat{\lambda}$ is first selected by cross-validation with the 
\verb|cv.glmnet| function of the \verb|R| package \verb|glmnet|, and where
the \verb|glmnet| function is then used to compute the selected model with
the LASSO with penalty parameter $\hat{\lambda}$. Then, if the selected
model does not contain the first explanatory variable, we discard the value
of $Y$, and else, we record whether the design-dependent target belongs to $%
\bar{CI}$ or not. The conditional coverage probability is then obtained by
taking the average number of times this is the case, over all the recorded
events.

\bibliographystyle{ims}
\bibliography{Biblio}

\end{document}